\documentclass[a4paper]{article}

\usepackage[cp1250]{inputenc}
\usepackage[IL2]{fontenc}

\usepackage{a4wide}
\usepackage[english]{babel}
\usepackage{euscript}
\usepackage{amstext,amsbsy,amscd,amssymb}
\usepackage{amsmath,enumitem}
\usepackage{amsfonts}
\usepackage{graphics}
\usepackage{graphicx}
\usepackage{rotating}
\usepackage{float}
\usepackage{refcount,xstring}
\usepackage{array}
\usepackage{microtype}
\usepackage{xcolor}
\usepackage{nicefrac}
\usepackage{dynkin-diagrams}
\usepackage{tikz,caption,subcaption}
\usetikzlibrary{arrows,decorations.pathmorphing,backgrounds,positioning,fit,petri}
\include{diagxy}
\allowdisplaybreaks

\let\rarr=\rightarrow
\let\xrarr=\xrightarrow

\let\veps=\varepsilon
\let\mcal=\mathcal
\let\mfrak=\mathfrak
\let\eus=\EuScript

\def\N{\mathbb{N}}
\def\Z{\mathbb{Z}}
\def\R{\mathbb{R}}
\def\C{\mathbb{C}}

\def\vac{|0\rangle}

\def\normOrd #1{{\mathop{:}\nolimits\!#1\!\mathop{:}\nolimits}}

\DeclareMathSymbol{\squares}{\mathord}{AMSa}{"03}

\newcolumntype{P}[1]{>{\centering\arraybackslash}p{#1}}

\DeclareMathOperator{\Ann}{Ann}

\DeclareMathOperator{\vspan}{span}
\DeclareMathOperator{\diag}{diag}
\DeclareMathOperator{\Spec}{Spec}
\DeclareMathOperator{\Specm}{Specm}
\DeclareMathOperator{\Res}{Res}
\DeclareMathOperator{\cdeg}{cdeg}

\def\Mod{\mathop {\rm Mod} \nolimits}

\def\End{\mathop {\rm End} \nolimits}

\def\ad{\mathop {\rm ad} \nolimits}
\def\gr{\mathop {\rm gr} \nolimits}
\def\im{\mathop {\rm im} \nolimits}

\def\id{{\rm id}}

\def\rank{\mathop {\rm rk} \nolimits}
\def\Ad{\mathop {\rm Ad} \nolimits}

\def\tr{\mathop {\rm tr} \nolimits}

\def\Com{\mathop {\rm Com} \nolimits}

\def\ch{\mathop {\rm ch} \nolimits}

\def\ind{\mathop {\rm ind} \nolimits}

\long\def\proof #1{\noindent \emph{Proof.}\ #1 \hfill $\squares$

\medskip}

\newcounter{num}[section]
\numberwithin{equation}{section}
\numberwithin{num}{section}
\newcounter{thmnum}
\renewcommand{\thethmnum}{\Alph{thmnum}}

\long\def\definition #1 {\refstepcounter{num} \noindent {\bf Definition \thenum.} #1

\medskip}

\long\def\theorem #1{\refstepcounter{num} \noindent {\bf Theorem \thenum.} #1

\medskip}

\long\def\maintheorem #1{\refstepcounter{thmnum} \noindent {\bf Main\,Theorem \thethmnum.} #1

\medskip}

\long\def\lemma #1{\refstepcounter{num}  \noindent {\bf Lemma \thenum.} #1

\medskip}

\long\def\proposition #1{\refstepcounter{num}  \noindent {\bf Proposition \thenum.} #1

\medskip}

\long\def\corollary #1{\refstepcounter{num}  \noindent {\bf Corollary \thenum.} #1

\medskip}

\long\def\remark #1{\noindent {\bf Remark.}\ #1}

\makeatletter
\newcommand{\raisemath}[1]{\mathpalette{\raisem@th{#1}}}
\newcommand{\raisem@th}[3]{\raisebox{#1}{$#2#3$}}
\makeatother

\makeatletter
\newcommand*\if@single[3]{%
  \setbox0\hbox{${\mathaccent"0362{#1}}^H$}%
  \setbox2\hbox{${\mathaccent"0362{\kern0pt#1}}^H$}%
  \ifdim\ht0=\ht2 #3\else #2\fi
  }
\newcommand*\rel@kern[1]{\kern#1\dimexpr\macc@kerna}
\newcommand*\widebar[1]{\@ifnextchar^{{\wide@bar{#1}{0}}}{\wide@bar{#1}{1}}}
\newcommand*\wide@bar[2]{\if@single{#1}{\wide@bar@{#1}{#2}{1}}{\wide@bar@{#1}{#2}{2}}}
\newcommand*\wide@bar@[3]{%
  \begingroup
  \def\mathaccent##1##2{%
    \if#32 \let\macc@nucleus\first@char \fi
    \setbox\z@\hbox{$\macc@style{\macc@nucleus}_{}$}%
    \setbox\tw@\hbox{$\macc@style{\macc@nucleus}{}_{}$}%
    \dimen@\wd\tw@
    \advance\dimen@-\wd\z@
    \divide\dimen@ 3
    \@tempdima\wd\tw@
    \advance\@tempdima-\scriptspace
    \divide\@tempdima 10
    \advance\dimen@-\@tempdima
    \ifdim\dimen@>\z@ \dimen@0pt\fi
    \rel@kern{0.6}\kern-\dimen@
    \if#31
      \overline{\rel@kern{-0.6}\kern\dimen@\macc@nucleus\rel@kern{0.4}\kern\dimen@}%
      \advance\dimen@0.4\dimexpr\macc@kerna
      \let\final@kern#2%
      \ifdim\dimen@<\z@ \let\final@kern1\fi
      \if\final@kern1 \kern-\dimen@\fi
    \else
      \overline{\rel@kern{-0.6}\kern\dimen@#1}%
    \fi
  }%
  \macc@depth\@ne
  \let\math@bgroup\@empty \let\math@egroup\macc@set@skewchar
  \mathsurround\z@ \frozen@everymath{\mathgroup\macc@group\relax}%
  \macc@set@skewchar\relax
  \let\mathaccentV\macc@nested@a
  \if#31
    \macc@nested@a\relax111{#1}%
  \else
    \def\gobble@till@marker##1\endmarker{}%
    \futurelet\first@char\gobble@till@marker#1\endmarker
    \ifcat\noexpand\first@char A\else
      \def\first@char{}%
    \fi
    \macc@nested@a\relax111{\first@char}%
  \fi
  \endgroup
}
\makeatother

\newcommand\rsmraise[1]{%
  \ifx#1\displaystyle .8\else
    \ifx#1\textstyle .8\else
      \ifx#1\scriptstyle .6\else
        .45%
      \fi
    \fi
  \fi}

\title{Generalized Grothendieck's simultaneous resolution
 and associated varieties of simple affine vertex algebras}

\author{Tomoyuki Arakawa, Vyacheslav Futorny, Libor K\v{r}i\v{z}ka}

\AtEndDocument{\bigskip{\footnotesize%
  (T.\,Arakawa) \textsc{Ningbo University, Ningbo, China and Research Institute for Mathematical Sciences, Kyoto University, Kyoto, Japan} \par
  \textit{E-mail address}: \texttt{arakawa@kurims.kyoto-u.ac.jp} \par
  \addvspace{\medskipamount}
  (V.\,Futorny) \textsc{International Center for Mathematics, SUSTech, Shenzhen, China} \par
  \textit{E-mail address}: \texttt{futorny@sustech.edu.cn} \par
  \addvspace{\medskipamount}
  (L.\,K\v{r}i\v{z}ka)  \textsc{International Center for Mathematics, SUSTech, Shenzhen, China} \par
  \textit{E-mail address}: \texttt{krizka.libor@gmail.com} \par
}}

\date{}

\begin{document}

\maketitle

\begin{abstract}
The closure of a Diximier sheet  is the image of a generalized Grothendieck's simultaneous resolution. We show that the associated variety of simple affine vertex algebras is contained in  the closure of the Diximier sheet when a chiralization of generalized Grothendieck's simultaneous resolution exists. This generalizes in a conceptual manner the results obtained by the first named author and Anne Moreau and, in particular, allows to extend them to the types $A_2$, $C_n$, $E_6$, $E_7$.
\end{abstract}

\medskip

\medskip


\thispagestyle{empty}

\tableofcontents


\section*{Introduction}
\addcontentsline{toc}{section}{Introduction}

\emph{Associated variety} \cite{Arakawa2012} of a vertex algebra $\mcal{V}$ is an affine Poisson variety $X_\mcal{V}$, defined as the spectrum of Zhu's $C_2$-algebra \cite{Zhu1996} of $\mcal{V}$.
This notion has turned out to be useful not only in the representation theory of vertex algebras (e.g.\ \cite{Arakawa2015b, Arakawa-Moreau2018, Arakawa-Moreau2017, Arakawa-Kawasetsu2018, Arakawa-Ekeren2023}), but also in the study of four dimensional $\mcal{N}=2$ superconformal field theories via the 4d/2d correspondence discovered in \cite{Beem-Lemos-Liendo2015}, see e.g.\ \cite{Song-Xie-Yan2017, Beem-Rastelli2018, Arakawa2018, Bonetti-Meneghelli-Rastelli2019, Dedushenko2021, Xie-Yan2021}.

In the case when $\mcal{V}$ is the simple affine vertex algebra $\mcal{V}_\kappa(\mfrak{g})$ associated with a finite-dimensional simple Lie algebra $\mfrak{g}$ and a $\mfrak{g}$-invariant symmetric bilinear form $\kappa$ on $\mfrak{g}$, the associated variety $X_{\mcal{V}_\kappa(\mfrak{g})}$ is a conic $G$-invariant subvariety of $\mfrak{g}^*$, as in the case of the associated variety \cite{Vogan1991} of a primitive ideal of $U(\mfrak{g})$, where $G$ is the adjoint group of $\mfrak{g}$. However, the associated variety $X_{\mcal{V}_\kappa(\mfrak{g})}$ need not be contained in the nilpotent cone $\mcal{N}(\mfrak{g})$ of $\mfrak{g}$, as in the case of the associated variety \cite{Vogan1991} of a primitive ideal of $U(\mfrak{g})$. In fact, we have $X_{\mcal{V}_\kappa(\mfrak{g})}=\mfrak{g} \simeq \mfrak{g}^*$ for a generic $\kappa$. On the other hand, there are cases where $X_{\mcal{V}_\kappa(\mfrak{g})}$ is a nilpotent orbit closure as in the case of the associated variety of a primitive ideal of $U(\mfrak{g})$, and this happens for instance when $\mcal{V}_\kappa(\mfrak{g})$ is admissible \cite{Arakawa2015} or when $\mfrak{g}$ is a member of the Deligne exceptional series and $\kappa=-\big(\smash{{h^\vee \over 6}}+1\big)\kappa_0$, where $\kappa_0$ is the normalized $\mfrak{g}$-invariant inner product \cite{Arakawa-Moreau2018}. There are also cases where $X_{\mcal{V}_\kappa(\mfrak{g})}$ is a proper subvariety of $\mfrak{g}$ that is not contained in $\mcal{N}(\mfrak{g})$. It was found in \cite{Arakawa-Moreau2017} that there are cases where $X_{\mcal{V}_\kappa(\mfrak{g})}$ are the closures of some \emph{Diximier sheets}, in types $A$ and $D$. It is natural to expect that there are also cases where $X_{\mcal{V}_\kappa(\mfrak{g})}$ are the closures of Diximier sheets in other types as well. However, it is difficult to generalize the method in \cite{Arakawa-Moreau2017}, because it is based on the explicit computation of singular vectors.

The aim of this paper is to generalize the result of \cite{Arakawa-Moreau2017} in a conceptual manner modeling a result of Borho and Brylinski \cite{Borho-Brylinski1982}, as we explain below in more details.
\smallskip

Let us consider a connected parabolic subgroup $P$ of $G$ with its Lie algebra $\mfrak{p}$. The closure of the Diximier sheet $\mcal{S}_\mfrak{p}$ determined by $\mfrak{p}$ is by definition the image of the \emph{generalized Grothendieck's simultaneous resolution}
\begin{align*}
  \mu_P \colon G \times_P [\mfrak{p},\mfrak{p}]^\perp \rarr \mfrak{g}.
  \end{align*}
Let $Y=G/(P,P)$ be the \emph{generalized base affine space} of $G$. Identifying $\mfrak{g}$ with $\mfrak{g}^*$ and $G \times_P [\mfrak{p},\mfrak{p}]^\perp$ with $T^* Y/A$, where $A$ is the torus $P/(P,P)$, $\mu_P$ coincides with the morphism $T^* Y/A \rarr \mfrak{g}^*$ induced by the momentum mapping $T^* Y \rarr \mfrak{g}^*$ for the $G$-action. Moreover, one finds that the natural algebra homomorphism
\begin{align*}
   \widetilde{\Phi}_X \colon  U(\mfrak{g}) \rarr  \Gamma(X, (\pi_*\mcal{D}_Y)^A)
\end{align*}
gives us a quantization of $\mu_P$, where $X=G/P$ is the \emph{generalized flag variety} for $G$ and $\pi\colon Y \rarr X$ is the canonical projection. A result in \cite{Borho-Brylinski1982} says that the associated variety of the two-sided ideal $\ker \smash{\widetilde{\Phi}_X}$ of $U(\mfrak{g})$ coincides with $\smash{\widebar{\mcal{S}^*_{\mfrak{p}}}}$.

Let $U$ be the maximal unipotent subgroup contained in the opposite parabolic subgroup of $P$. We denote by $U_e$ the big cell $UP \subset X$. Since the restriction $\Gamma(X, \pi_*(\mcal{D}_Y)^A)\rarr \Gamma(U_e, \pi_*(\mcal{D}_Y)^A)$ is injective, we have $\ker \smash{\widetilde{\Phi}_{U_e}} = \ker \smash{\widetilde{\Phi}_X}$, where
\begin{align*}
\widetilde{\Phi}_{U_e} \colon U(\mfrak{g}) \rarr \Gamma(X, \pi_*(\mcal{D}_Y)^A)\rarr \Gamma(U_e, \pi_*(\mcal{D}_Y)^A)
\end{align*}
is the composition. In fact, the algebra homomorphism $\smash{\widetilde{\Phi}_{U_e}}$ is a quantization of the restriction
\begin{align*}
 \widetilde{\mu}_P \colon U \times  [\mfrak{p},\mfrak{p}]^\perp  \rarr  \mfrak{g}
\end{align*}
of the morphism $\mu_P$.
\smallskip

Now let us consider the sheaf  $\mcal{D}_{Y,\kappa}^{ch}$ of the \emph{chiral differential operators (cdo)} \cite{BD,MSV} on $Y$, which is a natural chiralization of the sheaf $\mcal{D}_Y$ of differential operators on $Y$. It is a sheaf of vertex algebras on $Y$, whose sheaf of associated graded Poisson vertex algebras is isomorphic to $(\pi_Y)_*\mcal{O}_{J_\infty(\smash{T^*Y})}$.
Here $J_\infty(T^*Y)$ is the arc space of $T^*Y$ and $\pi_Y \colon J_\infty(T^*Y) \rarr T^*Y \rarr Y$ is the composition of the natural projections. A cdo on $X$ may not exist for a general smooth algebraic variety $X$, but for the generalized base affine space $Y$ the argument in \cite{Gorbounov-Malikov-Schechtman2001, Gorbounov-Malikov-Schechtman2004} shows that one can construct a cdo $\mcal{D}_{Y,\kappa}^{ch}$ on $Y$ from the cdo $\mcal{D}_{G,\kappa}$ on $G$ at level $\kappa$ by
performing the BRST reduction provided that 
\begin{align}
   (\kappa-\kappa_c^\mfrak{p})_{|[\mfrak{l},\mfrak{l}]} = 0, \label{eq;level-condition}
\end{align}
where $\mfrak{l}$ is the Levi subalgebra of $\mfrak{p}$
and $\kappa_c^\mfrak{p}$ is the $\mfrak{p}$-invariant symmetric bilinear form on $\mfrak{p}$ defined by $\kappa_c^\mfrak{p}(a,b) = - \tr_{\mfrak{g}/\mfrak{p}}(\ad(a)\ad(b))$, see Section \ref{subsection:Chiral differential operators}.

The $A$-action on $Y$ induces the $J_\infty(A)$-action on $\smash{\mcal{D}_{Y,\kappa}^{ch}}$.
Then one finds that the generalized Grothendieck's simultaneous resolution $T^*Y/A \rarr \mfrak{g}^*$ is chiralized to the vertex algebra homomorphism
\begin{align}
  \mcal{V}^{\kappa}(\mfrak{g}) \rarr \Gamma(X,(\pi_*\mcal{D}_{Y,\kappa}^{ch})^{A[[t]]}). \label{eq:chiral_generalized_Grothendieck's_simultaneous_resolution}
\end{align}
By composing \eqref{eq:chiral_generalized_Grothendieck's_simultaneous_resolution}
with the restriction mapping $\Gamma(X,(\pi_*\mcal{D}_{Y,\kappa}^{ch})^{A[[t]]}) \rarr \Gamma(U_e,(\pi_*\mcal{D}_{Y,\kappa}^{ch})^{A[[t]]})$, we obtain the vertex algebra homomorphism
\begin{align}\label{eq:gen.Wakimoto.intro}
  \widetilde{w}_{\kappa,\mfrak{g}}^\mfrak{p}\colon \mcal{V}^\kappa(\mfrak{g}) \rarr \mcal{M}_\mfrak{u} \otimes_\C \mcal{V}^{\kappa-\kappa_c^\mfrak{p}}(\mfrak{z}(\mfrak{l})),
\end{align}
which gives us a quantization of the morphism
\begin{align*}
   J_\infty(\widetilde{\mu}_P) \colon J_\infty(U \times [\mfrak{p},\mfrak{p}]^\perp) \rarr  J_\infty(\mfrak{g})
\end{align*}
under the assumption \eqref{eq;level-condition}.

Let us note that if $P$ is the Borel subgroup, then the condition $(\kappa-\kappa_c^\mfrak{p})_{|[\mfrak{l},\mfrak{l}]}=0$ is empty and \eqref{eq:gen.Wakimoto.intro} is the Wakimoto free field realization (\cite{Feigin-Frenkel1990,Frenkel2005}) of $\mcal{V}^\kappa(\mfrak{g})$, see Remark \ref{remark:Wakimoto}.

Our first main result is the following.
\medskip

\maintheorem{\label{mainth1} (Theorem \ref{thm:associated variety}) Let us suppose that $(\kappa-\kappa_c^\mfrak{p})_{|[\mfrak{l},\mfrak{l}]} = 0$. Then for the vertex algebra $\mcal{N}^\kappa_\mfrak{p}(\mfrak{g})=\mcal{V}^\kappa(\mfrak{g})/ \ker \widetilde{w}_{\kappa,\mfrak{g}}^\mfrak{p}$ holds
\begin{align*}
   X_{\smash{\mcal{N}^\kappa_\mfrak{p}\!(\mfrak{g})}} \simeq  \widebar{\mcal{S}^*_{\mfrak{p}}}.
\end{align*}
In particular, we have
$X_{\mcal{V}_{\kappa}(\mfrak{g})}\subset \widebar{\mcal{S}^*_{\mfrak{p}}}$.}

Next, we study the cases where the above condition is satisfied in more details. We exclude the case  $\mfrak{g}=\mfrak{sl}_2$, when $\smash{\widebar{\mcal{S}^*_{\mfrak{p}}}} =\mfrak{g}^*$, and focus on the cases with the abelian nilradical of $\mfrak{p}$, so that  we can obtain the explicit form of the morphism $w_{\kappa,\mfrak{g}}^\mfrak{p}$, see Theorem \ref{thm:FF homomorphism 1-graded}. The classification of simple Lie algebras and their opposite standard parabolic subalgebras with commutative nilradical is given in Table \ref{tab:parabolic commutative nilradical}.

Our second main result is the following.
\medskip

\maintheorem{\label{mainth2} (Theorem \ref{MainTheorem}) Let us suppose that the nilradical of $\mfrak{p}$ is commutative and that $(\kappa-\kappa_c^\mfrak{p})_{|[\mfrak{l},\mfrak{l}]} = 0$. The associated variety $\smash{X_{\mcal{V}_{\kappa}(\mfrak{g})}}$ is  isomorphic either to $\widebar{\mcal{S}_{\mfrak{p}}^*}$ or to $\widebar{\mcal{O}^*}$, where $\mcal{O}$ is a nilpotent orbit contained in $\smash{\widebar{\mcal{S}_{\mfrak{p}}}}$ and is explicitly determined in Table \ref{tab:associated varieties}. Moreover, we have that $\mcal{V}_{\kappa}(\mfrak{g})\simeq \mcal{N}^\kappa_\mfrak{p}(\mfrak{g})$ if $X_{\mcal{V}_{\kappa}(\mfrak{g})} \simeq \widebar{\mcal{S}_{\mfrak{p}}^*}$.}

Main\,Theorem \ref{mainth2} reproves \cite[Theorems 1.1, 1.2, 1.3]{Arakawa-Moreau2017}, \cite[Theorem 1.1\,(3)]{Arakawa-Moreau2018} and also \cite[Theorems 6.1, 7.1]{Arakawa-Moreau2018b}. The results for $\mfrak{g}=\mfrak{sl}_3, \mfrak{sp}_{2n}, \mfrak{e}_6, \mfrak{e}_7$ are new.

The proof of Main\,Theorem \ref{mainth2} does not involve the explicit form of the singular vectors. The main ingredient of the proof is Main\,Theorem \ref{mainth1} together with Theorem \ref{thm:DS reduction}, which allows us to make an inductive argument.

Finally, our third main result is the following.
\medskip

\maintheorem{(Theorem \ref{thm:DS reduction}\,(iii)) Let us suppose that the nilradical of $\mfrak{p}$ is commutative and that $(\kappa-\kappa_c^\mfrak{p})_{|[\mfrak{l},\mfrak{l}]} = 0$. Then the level $\kappa$ is collapsing for the minimal nilpotent element $f_\theta$.}

We denote by $\C$, $\R$, $\Z$, $\N_0$ and $\N$ the set of complex numbers, real numbers, integers, non-negative integers and positive integers, respectively. All algebras and modules are considered over the field of complex numbers.


\section{Affine Kac--Moody algebras, Weyl and Clifford algebras}

In this section we define smooth and induced modules for affine Kac--Moody algebras and introduce
a formalism for infinite-dimensional Weyl algebras and Clifford algebras.


\subsection{Affine Kac--Moody algebras}

Let $\mfrak{g}$ be a reductive finite-dimensional complex Lie algebra and let $\mfrak{h}$ be a Cartan subalgebra of $\mfrak{g}$. We denote by $\Delta$ the root system of $\mfrak{g}$ with respect to $\mfrak{h}$, by $\Delta_+$ a positive root system in $\Delta$ and by $\Pi \subset \Delta_+$ the set of simple roots. The standard Borel subalgebra $\mfrak{b}$ of $\mfrak{g}$ is defined through $\mfrak{b} = \mfrak{h} \oplus \mfrak{n}$ with the nilradical $\mfrak{n}$ and the opposite nilradical $\widebar{\mfrak{n}}$ given by
\begin{align*}
  \mfrak{n} = \bigoplus_{\alpha \in \Delta_+} \mfrak{g}_\alpha \qquad \text{and} \qquad \widebar{\mfrak{n}} = \bigoplus_{\alpha \in \Delta_+} \mfrak{g}_{-\alpha},
\end{align*}
where $\mfrak{g}_\alpha$ is the root subspace corresponding to a root $\alpha \in \Delta$. Besides, we have the corresponding triangular decomposition
\begin{align*}
  \mfrak{g} = \widebar{\mfrak{n}} \oplus \mfrak{h} \oplus \mfrak{n}
\end{align*}
of the Lie algebra $\mfrak{g}$. For a Lie subalgebra $\mfrak{a}$ of $\mfrak{g}$, we shall use the notation $\Delta_+^\mfrak{a} = \{\alpha \in \Delta_+;\, \mfrak{g}_\alpha \subset \mfrak{a}\}$.

Let us consider a $\mfrak{g}$-invariant symmetric bilinear form $\kappa$ on $\mfrak{g}$. The affine Kac--Moody algebra $\widehat{\mfrak{g}}_\kappa$ associated to $\mfrak{g}$ of level $\kappa$ is the $1$-dimensional central extension $\widehat{\mfrak{g}}_\kappa = \mfrak{g}(\!(t)\!) \oplus \C c$ of the formal loop algebra $\mfrak{g}(\!(t)\!)= \mfrak{g} \otimes_\C \C(\!(t)\!)$ with the commutation relations
\begin{align}
  [a \otimes f(t), b \otimes g(t)] = [a,b] \otimes f(t)g(t) - \kappa(a,b)\Res_{t=0} (f(t)dg(t))c, \label{eq:commutation relation}
\end{align}
where $a, b \in \mfrak{g}$, $f(t), g(t) \in \C(\!(t)\!)$ and $c$ is the central element of $\widehat{\mfrak{g}}_\kappa$. The Lie algebras $\widehat{\mfrak{g}}_\kappa$ and $\widehat{\mfrak{g}}_{\kappa'}$ are isomorphic, for $\mfrak{g}$-invariant symmetric bilinear forms $\kappa$ and $\kappa'$, if $\kappa'=k \kappa$ for some $k \in \C^\times$.
By introducing the notation $a_n = a \otimes t^n$ for $a \in \mfrak{g}$ and $n \in \Z$, the commutation relations \eqref{eq:commutation relation} can be simplified into the form
\begin{align}
  [a_m,b_n]=[a,b]_{m+n}+m \kappa(a,b) \delta_{m,-n} c \label{eq:commutation relation modes}
\end{align}
for $m,n \in \Z$ and $a,b \in \mfrak{g}$. Let us note that $\widehat{\mfrak{g}}_\kappa$ has the structure of a $\Z$-graded complete topological Lie algebra with the gradation defined by $\deg c =0$ and $\deg a_n = -n$ for $a\in \mfrak{g}$, $n\in \Z$.

Since $\mfrak{h}$ is a Cartan subalgebra of $\mfrak{g}$, we may define a Cartan subalgebra $\smash{\widehat{\mfrak{h}}}$ of $\widehat{\mfrak{g}}_\kappa$ by
\begin{align*}
  \widehat{\mfrak{h}} = \mfrak{h} \otimes_\C \C 1 \oplus \C c.
\end{align*}
The \emph{standard Borel subalgebra} $\smash{\widehat{\mfrak{b}}_{\rm st}}$ of $\widehat{\mfrak{g}}_\kappa$ is given by $\smash{\widehat{\mfrak{b}}_{\rm st}} = \smash{\widehat{\mfrak{h}}} \oplus \widehat{\mfrak{n}}_{\rm st}$, where the nilradical $\widehat{\mfrak{n}}_{\rm st}$ and the opposite nilradical $\widehat{\widebar{\mfrak{n}}}_{\rm st}$ are defined through
\begin{align*}
  \widehat{\mfrak{n}}_{\rm st} = \mfrak{n} \otimes_\C \C 1 \oplus \mfrak{g} \otimes_\C t\C[[t]] \qquad \text{and} \qquad \widehat{\widebar{\mfrak{n}}}_{\rm st} = \widebar{\mfrak{n}} \otimes_\C \C 1 \oplus \mfrak{g} \otimes_\C t^{-1}\C[t^{-1}].
\end{align*}
We have also the corresponding triangular decomposition
\begin{align*}
  \widehat{\mfrak{g}}_\kappa = \widehat{\widebar{\mfrak{n}}}_{\rm st} \oplus \widehat{\mfrak{h}} \oplus \widehat{\mfrak{n}}_{\rm st}.
\end{align*}
Let us note that while any two Borel subalgebras of $\mfrak{g}$ are conjugate by an automorphism of $\mfrak{g}$, any two Borel subalgebras of $\widehat{\mfrak{g}}_\kappa$ may not be conjugate by an
automorphism of $\widehat{\mfrak{g}}_\kappa$, see \cite{Futorny1997}. Besides, we introduce the maximal standard parabolic subalgebra $\widehat{\mfrak{g}}_{\rm st}$ of $\widehat{\mfrak{g}}_\kappa$ by
\begin{align*}
  \widehat{\mfrak{g}}_{\rm st} = \mfrak{g} \otimes_\C \C[[t]] \oplus \C c,
\end{align*}
whereas $\mfrak{g} \otimes_\C \C 1 \oplus \C c$ is the Levi subalgebra of $\widehat{\mfrak{g}}_{\rm st}$.
\medskip

\definition{\label{def:smooth module} Let $M$ be a $\widehat{\mfrak{g}}_\kappa$-module. We say that $M$ is a \emph{smooth} $\widehat{\mfrak{g}}_\kappa$-module if for each vector $v \in M$ there exists a positive integer $N_v \in \N$ such that
\begin{align*}
  (\mfrak{g} \otimes_\C t^{N_v}\C[[t]]) v = 0,
\end{align*}
or in other words that the Lie subalgebra $\mfrak{g} \otimes_\C t^{N_v}\C[[t]]$ annihilates $v$.
The category of smooth $\widehat{\mfrak{g}}_\kappa$-modules on which the central element $c$ acts as the identity we will denote by $\mcal{E}(\widehat{\mfrak{g}}_\kappa)$.}

Let us recall that by a graded $\widehat{\mfrak{g}}_\kappa$-module $M$ we mean a $\C$-graded vector space $M$ having the structure of a $\widehat{\mfrak{g}}_\kappa$-module compatible with the gradation of $\widehat{\mfrak{g}}_\kappa$. Let us also note that by shifting a given gradation on $M$ by a complex number we obtain a new gradation on $M$.
\medskip

\definition{\label{def:positive energy module} Let $M$ be a graded $\widehat{\mfrak{g}}_\kappa$-module. We say that $M$ is a \emph{positive energy} $\widehat{\mfrak{g}}_\kappa$-module if $M=\bigoplus_{n=0}^\infty M_{\lambda+n}$ and $M_\lambda \neq 0$, where $\lambda \in \C$. The category of positive energy $\widehat{\mfrak{g}}_\kappa$-modules on which the central element $c$ acts as the identity we will denote by $\mcal{E}_+(\widehat{\mfrak{g}}_\kappa)$.}

If $M$ is a positive energy $\widehat{\mfrak{g}}_\kappa$-module, then it follows immediately that $M$ is also a smooth $\widehat{\mfrak{g}}_\kappa$-module. Therefore, the category $\mcal{E}_+(\widehat{\mfrak{g}}_\kappa)$ is a full subcategory of $\mcal{E}(\widehat{\mfrak{g}}_\kappa)$.
\medskip

Now, we construct a class of $\widehat{\mfrak{g}}_\kappa$-modules, the so-called induced modules, which belong to the category $\mcal{E}_+(\widehat{\mfrak{g}}_\kappa)$. Let $E$ be a $\mfrak{g}$-module. Then the induced $\widehat{\mfrak{g}}_\kappa$-module
\begin{align*}
\mathbb{M}_{\kappa,\mfrak{g}}(E) = U(\widehat{\mfrak{g}}_\kappa)\otimes_{U(\widehat{\mfrak{g}}_{\rm st})}\! E,
\end{align*}
where $E$ is considered as the $\widehat{\mfrak{g}}_{\rm st}$-module on which $\mfrak{g} \otimes_\C t\C[[t]]$ acts trivially and $c$ acts as the identity, has the unique maximal
graded $\widehat{\mfrak{g}}_\kappa$-submodule $\mathbb{K}_{\kappa,\mfrak{g}}(E)$ having zero intersection with the $\mfrak{g}$-submodule $E$ of $\mathbb{M}_{\kappa,\mfrak{g}}(E)$. Therefore, we may set
\begin{align*}
  \mathbb{L}_{\kappa,\mfrak{g}}(E) = \mathbb{M}_{\kappa,\mfrak{g}}(E)/\mathbb{K}_{\kappa,\mfrak{g}}(E)
\end{align*}
for a $\mfrak{g}$-module $E$. Moreover, it is easy to see that if $E$ is a simple $\mfrak{g}$-module, then $\mathbb{L}_{\kappa,\mfrak{g}}(E)$ is also a simple $\widehat{\mfrak{g}}_\kappa$-module. The $\widehat{\mfrak{g}}_\kappa$-module $\mathbb{M}_{\kappa,\mfrak{g}}(E)$ is called the \emph{generalized Verma module} induced from $E$ for the maximal standard parabolic subalgebra $\widehat{\mfrak{g}}_{\rm st}$. Hence, we obtain the induction functor
\begin{align*}
  \mathbb{M}_{\kappa,\mfrak{g}} \colon \mcal{M}(\mfrak{g}) \rarr \mcal{E}_+(\widehat{\mfrak{g}}_\kappa)
\end{align*}
and the functor
\begin{align*}
  \mathbb{L}_{\kappa,\mfrak{g}} \colon \mcal{M}(\mfrak{g}) \rarr \mcal{E}_+(\widehat{\mfrak{g}}_\kappa),
\end{align*}
where $\mcal{M}(\mfrak{g})$ stands for the category of $\mfrak{g}$-modules. Let us recall that when $E$ is a simple finite-dimensional $\mfrak{g}$-module, then $\mathbb{M}_{\kappa,\mfrak{g}}(E)$ is usually called the \emph{Weyl module}.


\subsection{Weyl and Clifford algebras}

Let us consider the commutative algebra $\mcal{K}=\C(\!(t)\!)$ and let $\Omega_\mcal{K} = \C(\!(t)\!)dt$ be the corresponding module of K\"ahler differentials. Then for a finite-dimensional complex vector space $V$ we introduce the infinite-dimensional vector spaces $\mcal{K}(V) = V \otimes_\C \mcal{K}$ and $\Omega_\mcal{K}(V) = V \otimes_\C \Omega_\mcal{K}$. Let us note that $\mcal{K}(V)$ and $\Omega_\mcal{K}(V)$ are complete topological vector spaces with respect to the topology in which the basis of open neighbourhoods of $0$ is formed by $V \otimes_\C t^N\C[[t]]$ and $V \otimes_\C t^N\C[[t]]dt$ for $N \in \N_0$, respectively. Moreover, there is a pairing $(\cdot\,,\cdot) \colon \Omega_\mcal{K}(V^*) \otimes_\C \mcal{K}(V) \rarr \C$ given by
\begin{align}
  (\alpha \otimes f(t)dt, v \otimes g(t)) = \alpha(v) \Res_{t=0}(g(t)f(t)dt) \label{eq:non-degenerate form}
\end{align}
for $v \in V$, $\alpha \in V^*$ and $f(t),g(t) \in \mcal{K}$, which enables us to identify the restricted dual space to $\mcal{K}(V)$ with the vector space $\Omega_\mcal{K}(V^*)$, and vice versa. The pairing gives us non-degenerate bilinear forms $\omega_+$ and $\omega_-$ on $\Omega_\mcal{K}(V^*) \oplus \mcal{K}(V)$ defined by
\begin{align*}
   \omega_\pm(\alpha \otimes f(t)dt,v \otimes g(t)) = \pm \omega_\pm(v \otimes g(t), \alpha \otimes f(t)dt) = (\alpha \otimes f(t)dt, v \otimes g(t))
\end{align*}
and
\begin{align*}
  \omega_\pm(v \otimes f(t), w \otimes g(t)) = \omega_\pm(\alpha \otimes f(t)dt, \beta \otimes g(t)dt) = 0
\end{align*}
for $v,w \in V$, $\alpha,\beta \in V^*$, $f(t),g(t) \in \mcal{K}$. We have that $\omega_+$ is symmetric and $\omega_-$ is skew-symmetric. The Weyl algebra $\eus{A}_{\mcal{K}(V)}$ and the Clifford algebra $\eus{C}_{\mcal{K}(V)}$ are given by
\begin{align*}
  \eus{A}_{\mcal{K}(V)} = T(\Omega_\mcal{K}(V^*) \oplus \mcal{K}(V))/I_{\mcal{K}(V)} \qquad \text{and} \qquad \eus{C}_{\mcal{K}(V)} = T(\Omega_\mcal{K}(V^*) \oplus \mcal{K}(V))/J_{\mcal{K}(V)},
\end{align*}
where $I_{\mcal{K}(V)}$ and $J_{\mcal{K}(V)}$ are two-sided ideals of the tensor algebra $T(\Omega_\mcal{K}(V^*) \oplus \mcal{K}(V))$ generated by
\begin{align*}
  a \otimes b - b \otimes a + \omega_-(a,b)\cdot 1 \qquad \text{and} \qquad a \otimes b + b \otimes a - \omega_+(a,b)\cdot 1
\end{align*}
for $a, b  \in \Omega_\mcal{K}(V^*) \oplus \mcal{K}(V)$, respectively.

We define also a class of $\eus{A}_{\mcal{K}(V)}$-modules and $\eus{C}_{\mcal{K}(V)}$-modules called induced modules. Let us consider the vector subspaces $\mcal{L}_+$, $\mcal{L}_-$ and $\mcal{L}_0$ of $\Omega_\mcal{K}(V^*) \oplus \mcal{K}(V)$ defined by
\begin{gather*}
  \mcal{L}_- = V^*\! \otimes_\C t^{-2}\C[t^{-1}]dt \oplus V\! \otimes_\C t^{-1}\C[t^{-1}], \qquad  \mcal{L}_+ = V^*\! \otimes_\C \C[[t]]dt \oplus V\! \otimes_\C t\C[[t]], \\
  \mcal{L}_0 = V^*\! \otimes_\C \C t^{-1}dt \oplus V\! \otimes_\C \C 1.
\end{gather*}
Then we have the direct sum decomposition
\begin{align*}
  \Omega_\mcal{K}(V^*) \oplus \mcal{K}(V) = \mcal{L}_- \oplus \mcal{L}_0 \oplus \mcal{L}_+
\end{align*}
of $\Omega_\mcal{K}(V^*) \oplus \mcal{K}(V)$. It induces the triangular decomposition
\begin{align*}
  \eus{A}_{\mcal{K}(V)} \simeq \eus{A}_{\mcal{K}(V),-} \otimes_\C \eus{A}_{\mcal{K}(V),0} \otimes_\C \eus{A}_{\mcal{K}(V),+}
\end{align*}
of the Weyl algebra $\eus{A}_{\mcal{K}(V)}$, where
\begin{align*}
  \eus{A}_{\mcal{K}(V),-} \simeq S(\mcal{L}_-), \qquad \eus{A}_{\mcal{K}(V),0} \simeq \eus{A}_V, \qquad \eus{A}_{\mcal{K}(V),+} \simeq S(\mcal{L}_+)
\end{align*}
and $\eus{A}_V$ is the Weyl algebra of the vector space $V$, and also the triangular decomposition
\begin{align*}
  \eus{C}_{\mcal{K}(V)} \simeq \eus{C}_{\mcal{K}(V),-} \otimes_\C \eus{C}_{\mcal{K}(V),0} \otimes_\C \eus{C}_{\mcal{K}(V),+}
\end{align*}
of the Clifford algebra $\eus{C}_{\mcal{K}(V)}$, where
\begin{align*}
  \eus{C}_{\mcal{K}(V),-} \simeq \Lambda(\mcal{L}_-), \qquad \eus{C}_{\mcal{K}(V),0} \simeq \eus{C}_V, \qquad \eus{C}_{\mcal{K}(V),+} \simeq \Lambda(\mcal{L}_+)
\end{align*}
and $\eus{C}_V$ is the Clifford algebra associated with the vector space $V^* \oplus V$. Let us recall that $S(E)$ and $\Lambda(E)$ stand for the symmetric algebra and the exterior algebra, respectively, of a vector space $E$. Moreover, the algebras $\eus{A}_{\mcal{K}(V)}$ and $\eus{C}_{\mcal{K}(V)}$ are $\Z$-graded with the gradation determined by
\begin{align*}
  \deg(v \otimes t^n) = -n, \qquad \deg 1 = 0, \qquad \deg(\alpha \otimes t^{-n-1}dt) = n
\end{align*}
for $v \in V$, $\alpha \in V^*$ and $n \in \Z$.
\medskip

\definition{Let $M$ be an $\eus{A}_{\mcal{K}(V)}$-module. Then we say that $M$ is a \emph{smooth} $\eus{A}_{\mcal{K}(V)}$-module if for each vector $v \in M$ there exists a positive integer $N_v \in \N$ such that
\begin{align*}
  (V^*\! \otimes_\C t^{N_v}\C[[t]]dt \oplus V\! \otimes_\C t^{N_v}\C[[t]]) v = 0.
\end{align*}
We shall denote by $\mcal{E}(\eus{A}_{\mcal{K}(V)})$ the category of smooth  $\eus{A}_{\mcal{K}(V)}$-modules.}

Analogously as for $\widehat{\mfrak{g}}_\kappa$-modules, we may introduce graded $\eus{A}_{\mcal{K}(V)}$-modules and positive energy $\eus{A}_{\mcal{K}(V)}$-modules. We shall denote by $\mcal{E}_+(\eus{A}_{\mcal{K}(V)})$ the category of positive energy $\eus{A}_{\mcal{K}(V)}$-modules.
\medskip

Let $E$ be an $\eus{A}_V$-module. Then the induced $\eus{A}_{\mcal{K}(V)}$-module
\begin{align*}
  \mathbb{M}_{\mcal{K}(V)}^\eus{A}(E) = \eus{A}_{\mcal{K}(V)} \otimes_{\eus{A}_{\mcal{K}(V),0} \otimes_\C \eus{A}_{\mcal{K}(V),+}}\! E,
\end{align*}
where $E$ is considered as the $\eus{A}_{\mcal{K}(V),0} \otimes_\C \eus{A}_{\mcal{K}(V),+}$-module on which the Weyl algebra $\eus{A}_{\mcal{K}(V),0}$ acts via the canonical isomorphism $\eus{A}_{\mcal{K}(V),0} \simeq \eus{A}_V$ and $\eus{A}_{\mcal{K}(V),+}$ acts trivially, has a unique maximal $\eus{A}_{\mcal{K}(V)}$-submodule $\smash{\mathbb{K}_{\mcal{K}(V)}^\eus{A}(E)}$ having zero intersection with the $\eus{A}_V$-submodule $E$ of $\smash{\mathbb{M}_{\mcal{K}(V)}^\eus{A}(E)}$. Hence, we may set
\begin{align*}
  \mathbb{L}_{\mcal{K}(V)}^\eus{A}(E) = \mathbb{M}_{\mcal{K}(V)}^\eus{A}(E)/\mathbb{K}_{\mcal{K}(V)}^\eus{A}(E)
\end{align*}
for an $\eus{A}_V$-module $E$. Moreover, it is easy to see that if $E$ is a simple $\eus{A}_V$-module, then $\smash{\mathbb{L}_{\mcal{K}(V)}^\eus{A}(E)}$ is also a simple $\eus{A}_{\mcal{K}(V)}$-module. Therefore, we have the induction functor
\begin{align}
  \mathbb{M}_{\mcal{K}(V)}^\eus{A} \colon \mcal{M}(\eus{A}_V) \rarr \mcal{E}_+(\eus{A}_{\mcal{K}(V)})
\end{align}
and the functor
\begin{align}
  \mathbb{L}_{\mcal{K}(V)}^\eus{A} \colon \mcal{M}(\eus{A}_V) \rarr \mcal{E}_+(\eus{A}_{\mcal{K}(V)}),
\end{align}
where $\mcal{M}(\eus{A}_V)$ is the category of $\eus{A}_V$-modules. The induced $\eus{A}_{\mcal{K}(V)}$-module $\smash{\mathbb{M}_{\mcal{K}(V)}^\eus{A}(S(V^*))}$ will be of a particular importance since it has a natural structure of a vertex algebra as we will see in the next section.
\medskip

Let $\mfrak{g}$ be a reductive Lie algebra and let $\mfrak{b}=\mfrak{h} \oplus \mfrak{n}$ be a Borel subalgebra of $\mfrak{g}$ with the Cartan subalgebra $\mfrak{h}$, the nilradical $\mfrak{n}$ and the opposite nilradical $\widebar{\mfrak{n}}$. Let us consider a parabolic subalgebra $\mfrak{p}$ of $\mfrak{g}$. We assume that $\mfrak{p}$ contains the opposite Borel subalgebra $\mfrak{h} \oplus \widebar{\mfrak{n}}$ of $\mfrak{g}$. Let $\mfrak{p} = \mfrak{l} \oplus \widebar{\mfrak{u}}$ be the Levi decomposition of $\mfrak{p}$ with the Levi subalgebra $\mfrak{l}$, the nilradical $\widebar{\mfrak{u}}$ and the opposite nilradical $\mfrak{u}$. Let $\{e_\alpha;\, \alpha \in \Delta_+^\mfrak{u}\}$ be a root basis of the nilpotent Lie subalgebra $\mfrak{u}$. We denote by $\{x_\alpha;\, \alpha \in \Delta_+^\mfrak{u}\}$ the linear coordinate functions on $\mfrak{u}$ with respect to the given basis of $\mfrak{u}$. It follows immediately that the set $\{e_\alpha \otimes t^n;\, \alpha \in \Delta_+^\mfrak{u},\, n \in \Z\}$ forms a topological basis of $\mcal{K}(\mfrak{u})$ and the set $\{x_\alpha \otimes t^{-n-1}dt;\, \alpha \in \Delta_+^\mfrak{u},\, n \in \Z\}$ forms the dual topological basis of $\Omega_\mcal{K}(\mfrak{u}^*)$ with respect to the pairing \eqref{eq:non-degenerate form}, i.e.\ we have
\begin{align*}
  (x_\alpha \otimes t^{-n-1}dt,e_\beta \otimes t^m)=x_\alpha(e_\beta)\Res_{t=0} t^{m-n-1}dt= \delta_{\alpha,\beta} \delta_{n,m}
\end{align*}
for $\alpha,\beta \in \Delta_+^\mfrak{u}$ and $m,n \in \Z$. If we denote $x_{\alpha,n}=x_\alpha \otimes t^{-n-1}dt$ and $\partial_{x_{\alpha,n}}=e_\alpha \otimes t^n$ for $\alpha \in \Delta_+^\mfrak{u}$ and $n \in \Z$, then the two-sided ideal $I_{\mcal{K}(\mfrak{u})}$ is generated by elements
\begin{align*}
\bigg(\sum_{n \in \Z} a_nx_{\alpha,n}\bigg)\! \otimes \!\bigg(\sum_{m\in \Z} b_m\partial_{x_{\beta,m}}\! \bigg) - \bigg(\sum_{m\in \Z} b_m\partial_{x_{\beta,m}}\! \bigg) \!\otimes\! \bigg(\sum_{n \in \Z} a_nx_{\alpha,n}\bigg)+\delta_{\alpha,\beta}\bigg(\sum_{n\in \Z} a_nb_n\bigg)\!\cdot 1
\end{align*}
and it coincides with the canonical commutation relations
\begin{align*}
[x_{\alpha,n},\partial_{x_{\beta,m}}]=-\delta_{\alpha,\beta}\delta_{n,m}
\end{align*}
for $\alpha,\beta \in \Delta_+^\mfrak{u}$ and $m,n \in \Z$. Therefore, we obtain that the Weyl algebra $\eus{A}_{\mcal{K}(\mfrak{u})}$ is topologically generated by the set $\{x_{\alpha,n}, \partial_{x_{\alpha,n}};\, \alpha \in \Delta_+^\mfrak{u},\, n \in \Z\}$ with the canonical commutation relations.
\medskip

Let us note that all the notions introduced above can also be defined for the Clifford algebra $\eus{C}_{\mcal{K}(V)}$ with appropriate modifications.


\section{Vertex algebras and their modules}

In this section we recall some basic facts on vertex algebras, for more details see \cite{Borcherds1986}, \cite{Kac1998}, \cite{Dong-Li-Mason1998}, \cite{Frenkel-Ben-Zvi2004}, \cite{Frenkel2007-book}.


\subsection{Vertex algebras}

Let $R$ be an algebra over $\C$. Then an $R$-valued formal power series (or formal distribution) in the variables $z_1,\dots,z_n$ is a series
\begin{align*}
  a(z_1,\dots,z_n)=\sum_{m_1,\dots,m_n \in \Z} a_{m_1,\dots,m_n} z_1^{m_1} \dotsm z_n^{m_n},
\end{align*}
where $a_{m_1,\dots,m_n} \in R$. The complex vector space of all $R$-valued formal power series is denoted by $R[[z_1^{\pm 1},\dots,z_n^{\pm 1}]]$. For a formal power series $a(z) = \sum_{n \in \Z} a_n z^n$, the residue is given by
\begin{align*}
  \Res_{z=0} a(z) = a_{-1}.
\end{align*}
A particulary important example of a $\C$-valued formal power series in two variables $z, w$ is the formal delta function $\delta(z-w)$ defined by
\begin{align*}
  \delta(z-w) = \sum_{n \in \Z} z^n w^{-n-1}.
\end{align*}
Let $V$ be a complex vector space, so $\End V$ is an associative algebra over $\C$. We say that a formal power series $a(z) \in \End V[[z^{\pm 1}]]$ is a field, if $a(z)v \in V(\!(z)\!)$ for all $v \in V$.
We shall write the field $a(z)$ as
\begin{align}
  a(z) = \sum_{n \in \Z} a_{(n)} z^{-n-1}.
\end{align}
The vector space of all fields on $V$ in the variable $z$ we will denote by $\eus{F}(V)(z)$.
\medskip

\definition{A vertex algebra consists of the following data:
\begin{enumerate}[topsep=3pt,itemsep=0pt]
  \item[1)] a complex vector space $\mcal{V}$ (the space of states);
  \item[2)] a non-zero vector $\vac \in \mcal{V}$ (the vacuum vector);
  \item[3)] an endomorphism $T \colon \mcal{V} \rarr \mcal{V}$ (the translation operator);
  \item[4)] a linear mapping $Y(\,\cdot\,,z) \colon \mcal{V} \rarr \End \mcal{V}[[z^{\pm 1}]]$ sending
  \begin{align*}
    a \in \mcal{V} \mapsto Y(a,z) = \sum_{n \in \Z} a_{(n)} z^{-n-1} \in \eus{F}(\mcal{V})(z)
  \end{align*}
  (the state-field correspondence)
\end{enumerate}
satisfying the subsequent axioms:
\begin{enumerate}[topsep=3pt,itemsep=0pt]
  \item[1)] $Y(\vac,z)= \id_\mcal{V}$, $Y(a,z)\vac_{|z=0}=a$ (the vacuum axiom);
  \item[2)] $T\vac=0$, $[T,Y(a,z)] = \partial_z Y(a,z)$ (the translation axiom);
  \item[3)] for all $a,b \in \mcal{V}$, there is a non-negative integer $N_{a,b} \in \N_0$ such that
  \begin{align*}
    (z-w)^{N_{a,b}}[Y(a,z),Y(b,w)]=0
  \end{align*}
  (the locality axiom).
\end{enumerate}

A vertex algebra $\mcal{V}$ is called $\Z$-graded if $\mcal{V}$ is a $\Z$-graded vector space, $\vac$ is a vector of degree $0$, $T$ is an endomorphism of degree $1$, and for $a \in \mcal{V}_m$ the field $Y(a,z)$ has conformal dimension $m$, i.e.\ we have
\begin{align*}
  \deg a_{(n)} = -n+m-1
\end{align*}
for all $n \in \Z$.}

Let us note that according to the translation axiom, the action of $T$ on the space of states $\mcal{V}$ is completely determined by $Y$, since we have $T(a) = a_{(-2)}\vac$. Moreover, we have $a = a_{(-1)}\vac$. The field $Y(a,z)$ for $a \in \mcal{V}$ is usually called the \emph{vertex operator}.
\medskip

We say that a $\Z$-graded vertex algebra $\mcal{V}$ is conformal, of central charge $c_\mcal{V} \in \C$, if we are given a non-zero conformal vector $\omega \in \mcal{V}_2$ such that the Fourier coefficients of the corresponding vertex operator
\begin{align}
  Y(\omega,z) = \sum_{n\in \Z} L_n z^{-n-2}
\end{align}
satisfy the defining relations of the Virasoro algebra with central charge $c_\mcal{V} \in \C$, and in addition we have $L_{-1} = T$ and $L_0{}_{|\mcal{V}_m} = m\, \id_{\mcal{V}_m}$ for $m \in \Z$. We call $L_0$ the energy operator, it is a Hamiltonian operator.


\subsection{Associated variety and singular support of vertex algebras}

In \cite{Li2005}, Li has shown that every vertex algebra has a canonical decreasing filtration. For a vertex algebra $\mcal{V}$, let $F^p\mcal{V}$ for $p \in \N_0$ be a vector subspace of $\mcal{V}$ spanned by the elements
\begin{align*}
  a_{1,(-n_1-1)}a_{2,(-n_2-1)} \cdots a_{r,(-n_r-1)}\vac
\end{align*}
for $r \in \N_0$, $a_1,a_2,\dots,a_r \in \mcal{V}$, $n_1,n_2,\dots,n_r \in \N_0$ and $n_1+n_2+\dots +n_r \geq p$. It follows easily that
\begin{align*}
  \mcal{V} = F^0\mcal{V} \supset F^1 \mcal{V} \supset F^2\mcal{V} \supset \dots
\end{align*}
and $T(F^p\mcal{V}) \subset F^{p+1}\mcal{V}$ for $p\in \N_0$. Further, we define
\begin{align*}
  \gr^F\!\mcal{V} = \bigoplus_{p \in \N_0} F^p\mcal{V}/F^{p+1}\mcal{V}
\end{align*}
together with the canonical quotient mapping
\begin{align*}
  \sigma_p \colon F^p\mcal{V} \rarr F^p\mcal{V}/F^{p+1}\mcal{V}
\end{align*}
for $p \in \N_0$. Then the vector space $\gr^F\!\mcal{V}$ has the natural structure of a Poisson vertex algebra, see \cite{Li2005}. The vector subspace
\begin{align}
  R_\mcal{V} = F^0\mcal{V}/F^1\mcal{V} = \mcal{V}/C_2(\mcal{V}) \subset \gr^F\!\mcal{V},
\end{align}
where
\begin{align*}
  C_2(\mcal{V}) = \vspan_\C\{a_{(-2)}b;\, a,b \in \mcal{V}\} = F^1\mcal{V},
\end{align*}
is called the Zhu's $C_2$-algebra of $\mcal{V}$. The Poisson vertex algebra structure of $\gr^F\!\mcal{V}$ restricts to the Poisson algebra structure of $R_\mcal{V}$ given by
\begin{align*}
  \sigma_{\rm Zhu}(a) \cdot \sigma_{\rm Zhu}(b) = \sigma_{\rm Zhu}(a_{(-1)}b) \qquad \text{and} \qquad \{\sigma_{\rm Zhu}(a), \sigma_{\rm Zhu}(b)\} = \sigma_{\rm Zhu}(a_{(0)}b)
\end{align*}
for $a, b \in \mcal{V}$, where $\sigma_{\rm Zhu}$ denotes the canonical projection from $\mcal{V}$ to $R_\mcal{V}$. Let us note that we have
\begin{align*}
  C_2(\mcal{V}) = \vspan_\C\{a_{(-n)}b;\, a,b \in \mcal{V},\, n \geq 2\}
\end{align*}
by the translation axiom. Since $R_\mcal{V}$ is a commutative algebra, we introduce its spectrum $\smash{\widetilde{X}}_\mcal{V}$ and maximal spectrum $X_\mcal{V}$ through
\begin{align}
  \widetilde{X}_\mcal{V} = \Spec R_\mcal{V} \qquad \text{and} \qquad X_\mcal{V} = \Specm R_\mcal{V}
\end{align}
and we shall call them the \emph{associated scheme} of $\mcal{V}$ and the \emph{associated variety} of $\mcal{V}$, respectively, see \cite{Arakawa2012}. A vertex algebra $\mcal{V}$ is called \emph{finitely strongly generated} if $R_\mcal{V}$ is finitely generated as an algebra.
\medskip

For a scheme $X$ of finite type over $\C$, we may introduce the scheme $J_\infty(X)$ called the \emph{arc space} of $X$ (or the \emph{infinite jet scheme} of $X$). In the case when $X = \Spec R$ for a finitely generated commutative algebra $R$ over $\C$, i.e.\ we have
\begin{align*}
  R \simeq \C[x^i;\, i=1,2,\dots,n]/\langle f_r;\, r=1,2,\dots,m \rangle
\end{align*}
with $f_r \in \C[x^i;\, i=1,2,\dots,n]$ for $r=1,2,\dots,m$, the infinite jet scheme $J_\infty(X)$ is given by
\begin{align*}
  J_\infty(X) = \Spec J_\infty(R),
\end{align*}
where
\begin{align*}
  J_\infty(R) \simeq \C[x^i_{(-j)};\, j \in \N_0,\, i=1,2,\dots,n]/\langle T^j\!f_r;\, j \in \N_0,\, r=1,2,\dots,m\rangle
\end{align*}
and $T$ is a derivation of the algebra $\C[x^i_{(-j)};\, j \in \N_0,\, i=1,2,\dots,n]$ defined by
\begin{align*}
  Tx^i_{(-j)} = jx^i_{(-j-1)}
\end{align*}
for $i=1,2,\dots,n$ and $j \in \N_0$. We identify $x^i$ with $\smash{x^i_{(-1)}}$ for $i=1,2,\dots,n$. In fact, the coordinate ring $J_\infty(R) = \C[J_\infty(X)]$ of its infinite jet scheme $J_\infty(X)$ is a differential algebra, hence it is a commutative vertex algebra, since there is a one-to-one correspondence between commutative vertex algebras and differential algebras (\cite{Borcherds1986}). Besides, the differential algebra $(J_\infty(R),T)$ has the following universal property. We have the canonical algebra homomorphism $j \colon R \rarr J_\infty(R)$ such that if we are given a differential algebra $(A,\partial)$ and an algebra homomorphism $\varphi \colon R \rarr A$, then there is a unique homomorphism $\widetilde{\varphi} \colon J_\infty(R) \rarr A$ of differential algebras such that the diagram
\begin{align*}
\bfig
\Vtriangle(0,0)|alr|<500,500>[R`J_\infty(R)`A;j`\varphi`\widetilde{\varphi}]
\efig
\end{align*}
is commutative.

It was shown in \cite[Lemma 4.2]{Li2005} that $\gr^F\!\mcal{V}$ is generated by $R_\mcal{V}$ as a differential algebra. Hence, we have a surjective homomorphism $J_\infty(R_\mcal{V}) \rarr \gr^F\!\mcal{V}$ of differential algebras. In fact, it is a homomorphism of Poisson vertex algebras, see \cite[Lemma 4.2]{Li2005}, \cite[Proposition 2.5.1]{Arakawa2012}.
\medskip

\definition{\cite{Arakawa2012} Let $\mcal{V}$ be a vertex algebra. Then the \emph{singular support} of $\mcal{V}$ is
defined as
\begin{align*}
  SS(\mcal{V}) = \Spec \gr^F\!\mcal{V} \subset J_\infty(\widetilde{X}_\mcal{V}),
\end{align*}
where $\widetilde{X}_\mcal{V} = \Spec R_\mcal{V}$.}

If $\mcal{V}$ is an $\N_0$-graded vertex algebra with the gradation given by a Hamiltonian operator $H$, i.e.\ we have
\begin{align*}
  \mcal{V} = \bigoplus_{\Delta \in \N_0} \mcal{V}_\Delta, \qquad \mcal{V}_\Delta = \{v\in\mcal{V};\, Hv =\Delta v\},
\end{align*}
then there is another increasing filtration on $\mcal{V}$ introduced in \cite{Li2005}, called the \emph{conformal weight filtration}, defined as follows. Let $\{a^i;\, i \in I\}$ be a set of homogenous strong generators of $\mcal{V}$. Let $G_p\mcal{V}$ for $p \in \N_0$ be a vector subspace of $\mcal{V}$ spanned by the elements
\begin{align*}
  a_{(-n_1-1)}^{i_1}a_{(-n_2-1)}^{i_2} \cdots a_{(-n_r-1)}^{i_r}\vac
\end{align*}
for $r \in \N_0$, $i_1,i_2,\dots,i_r \in I$, $n_1,n_2,\dots,n_r \in \N_0$ and $\Delta_{a_1}+\Delta_{a_2}+ \dots + \Delta_{a_r} \leq p$. It easily follows that
\begin{align*}
  0 = G_{-1}\mcal{V} \subset G_0\mcal{V} \subset G_1\mcal{V} \dots
\end{align*}
and $T(G_p\mcal{V}) \subset G_p\mcal{V}$ for $p\in \N_0$. Further, the vector space
\begin{align*}
  \gr^G\!\mcal{V} = \bigoplus_{p \in \N_0} G_p\mcal{V}/G_{p-1}\mcal{V}
\end{align*}
has the natural structure of a Poisson vertex algebra, see \cite{Li2005}. Moreover, we have that
\begin{align*}
  \gr^F\! \mcal{V} \simeq \gr^G\! \mcal{V}
\end{align*}
as Poisson vertex algebras, see \cite[Proposition 2.6.1]{Arakawa2012}.


\subsection{Positive energy modules}

In \cite{Zhu1996}, Zhu introduced a functorial construction which assigns to a $\Z$-graded vertex algebra an associative algebra known as the Zhu's algebra. For a $\Z$-graded vertex algebra $\mcal{V}$, we define a bilinear mapping on $\mcal{V}$ by
\begin{align}
  a * b = \Res_{z=0}\!\bigg(Y(a,z){(1+z)^{\deg a} \over z}\,b\bigg) = \sum_{i=0}^{\deg a} \binom{\deg a}{i} a_{(i-1)}b \label{eq:Zhu algebra mult}
\end{align}
for homogeneous elements $a, b \in \mcal{V}$ and extend linearly. The Zhu's algebra $A(\mcal{V})$ is defined as
\begin{align}
  A(\mcal{V}) = \mcal{V}/O(\mcal{V}),
\end{align}
where $O(\mcal{V})$ is the vector subspace of $\mcal{V}$ spanned by
\begin{align*}
  \Res_{z=0}\!\bigg(Y(a,z){(1+z)^{\deg a} \over z^2}\,b\bigg) = \sum_{i=0}^{\deg a} \binom{\deg a}{i} a_{(i-2)}b
\end{align*}
for homogeneous elements $a,b \in \mcal{V}$. The bilinear mapping \eqref{eq:Zhu algebra mult} gives rise to an associative structure on the quotient $A(\mcal{V})$. We denote by $\pi_{\rm Zhu}$ the canonical projection from $\mcal{V}$ to $A(\mcal{V})$.

In addition, if $\mcal{V}$ is an $\N_0$-graded vertex algebra, we may define an increasing, separated and exhaustive filtration $\{F_p\mcal{V}\}_{p \in \Z}$ on $\mcal{V}$ by
\begin{align*}
  F_p \mcal{V} = \bigoplus_{m=0}^p \mcal{V}_m.
\end{align*}
Let $F_pA(\mcal{V})$ be the image of $F_p\mcal{V}$ in $A(\mcal{V})=\mcal{V}/O(\mcal{V})$. Then the associated $\Z$-graded algebra
\begin{align*}
  \gr A(\mcal{V}) = \bigoplus_{p \in \Z} F_pA(\mcal{V})/F_{p-1}A(\mcal{V})
\end{align*}
has the natural structure of a Poisson algebra given by
\begin{align*}
  \tau_{\deg a}(a) \cdot \tau_{\deg b}(b) = \tau_{\deg a + \deg b}(a_{(-1)}b) \qquad \text{and} \qquad \{\tau_{\deg a}(a), \tau_{\deg b}(b)\} = \tau_{\deg a + \deg b-1}(a_{(0)}b)
\end{align*}
for homogeneous elements $a, b \in \mcal{V}$, where
\begin{align*}
  \tau_p \colon F_p\mcal{V} \rarr \gr_p A(\mcal{V}) = F_pA(\mcal{V})/F_{p-1}A(\mcal{V})
\end{align*}
for $p \in \Z$ is the canonical projection. Besides, there is a natural surjective homomorphism
\begin{align}
  \eta_\mcal{V} \colon R_\mcal{V} \rarr \gr A(\mcal{V})
\end{align}
of Poisson algebras given by
\begin{align*}
  \eta_\mcal{V}(\sigma_{\rm Zhu}(a)) = \tau_{\deg a}(a)
\end{align*}
for a homogeneous element $a \in \mcal{V}$, see \cite{Arakawa-Lam-Yamada2014}. Let us note that the homomorphism $\eta_\mcal{V}$ is not an isomorphism in general.
\medskip

Let us suppose that we are given a $\Z$-graded vertex algebra $\mcal{V}$. Then a $\mcal{V}$-module $M$ is called graded if $M$ is a $\C$-graded vector space and for $a \in \mcal{V}_m$ the field $Y_M(a,z)$ has conformal dimension $m$, i.e.\ we have
\begin{align*}
  \deg a_{(n)}^M = -n+m-1
\end{align*}
for all $n \in \Z$. Let us note that by shifting a given gradation on $M$ by a complex number we obtain a new gradation on $M$. For a $\mcal{V}$-module $M$, we define
\begin{align*}
  o_M(a) =a^M_{(\deg a -1)}
\end{align*}
for a homogeneous element $a \in \mcal{V}$. It easily follows that for a homogeneous element $a \in \mcal{V}$ the operator $o_M(a)$ preserves the homogeneous components of any graded $\mcal{V}$-module $M$.
\medskip

\definition{Let $\mcal{V}$ be a $\Z$-graded vertex algebra. We say that a graded $\mcal{V}$-module $M$ is a \emph{positive energy $\mcal{V}$-module} if $M=\bigoplus_{n=0}^\infty M_{\lambda+n}$ and $M_\lambda \neq 0$, where $\lambda \in \C$. Moreover, we denote by $M_{\rm top}$ the top degree component $M_\lambda$ of $M$. The category of positive energy $\mcal{V}$-modules we will denote by $\mcal{E}_+(\mcal{V})$.}

As the following theorem proved in \cite{Zhu1996} shows, the Zhu's algebra $A(\mcal{V})$ plays a prominent role in the representation theory of vertex algebras.
\medskip

\theorem{\label{thm:Zhu correspondence}
Let $\mcal{V}$ be a $\Z$-graded vertex algebra and let $M$ be a positive energy $\mcal{V}$-module. Then the top degree component $M_{\rm top}$ is an $A(\mcal{V})$-module, where the action of $\pi_{\rm Zhu}(a) \in A(\mcal{V})$ for $a \in \mcal{V}$ is given by $o_M(a)$. In addition, the correspondence $M \mapsto M_{\rm top}$ gives a bijection between the set of simple positive energy $\mcal{V}$-modules and that of simple $A(\mcal{V})$-modules.}

To a $\Z$-graded vertex algebra $\mcal{V}$ we may associate a complete topological Lie algebra $U(\mcal{V})$, first introduced by Borcherds \cite{Borcherds1986}, by
\begin{align*}
  U(\mcal{V}) = (\mcal{V} \otimes_\C \C(\!(t)\!)) / \im \partial,
\end{align*}
where
\begin{align*}
  \partial = T \otimes \id + \id \otimes \partial_t.
\end{align*}
If we denote by $a_{[n]}$ for $a \in \mcal{V}$ and $n \in \Z$ the projection of $a \otimes t^n \in \mcal{V} \otimes_\C \C(\!(t)\!)$ onto $U(\mcal{V})$, then the Lie bracket on $U(\mcal{V})$ is determined by
\begin{align*}
  [a_{[m]}, b_{[n]}] = \sum_{k=0}^\infty \binom{m}{k} (a_{(k)}b)_{[m+n-k]}
\end{align*}
for $a, b \in \mcal{V}$ and $m,n \in \Z$.

For a homogeneous element $a \in \mcal{V}$, we set $\deg a_{[n]} = -n+\deg a -1$. Then the degree assignment to the elements of $U(\mcal{V})$ gives us a triangular decomposition
\begin{align}
  U(\mcal{V}) = U(\mcal{V})_- \oplus U(\mcal{V})_0 \oplus U(\mcal{V})_+ \label{eq:triangular docomposition U(V)}
\end{align}
of the Lie algebra $U(\mcal{V})$ together with a canonical surjective homomorphism
\begin{align*}
  U(\mcal{V})_0  \rarr A(\mcal{V})
\end{align*}
of Lie algebras defined by
\begin{align*}
  a_{[\deg a -1 ]} \mapsto \pi_{\rm Zhu}(a)
\end{align*}
for a homogeneous element $a \in \mcal{V}$, see \cite{DeSole-Kac2006}.

Further, let us consider a $\mcal{V}$-module $M$. Then it has also a natural structure of a $U(\mcal{V})$-module since we have a canonical homomorphism
\begin{align*}
  U(\mcal{V}) \rarr \End \mcal{V}
\end{align*}
of Lie algebras defined through
\begin{align*}
  a_{[n]} \mapsto a_{(n)}
\end{align*}
for $a \in \mcal{V}$ and $n \in \Z$. We denote by $\Omega_\mcal{V}(M)$ the vector subspace of $M$ consisting of lowest weight vectors, i.e.\ we have
\begin{align}
  \Omega_\mcal{V}(M) = \{v\in M;\, U(\mcal{V})_-v=0\}.
\end{align}
It follows immediately using the triangular decomposition of $U(\mcal{V})$ that $\Omega_\mcal{V}(M)$ is a $U(\mcal{V})_0$-module. Moreover, by \cite{Zhu1996} we have that $\Omega_\mcal{V}(M)$ is an $A(\mcal{V})$-module, where the action of $\pi_{\rm Zhu}(a) \in A(\mcal{V})$ for $a \in \mcal{V}$ is given by $o_M(a)$. It is clear that
\begin{align}
  \Omega_\mcal{V} \colon \mcal{E}(\mcal{V}) \rarr \mcal{M}(A(\mcal{V}))
\end{align}
is a functor, where $\mcal{E}(\mcal{V})$ stands for the category of $\mcal{V}$-modules and $\mcal{M}(A(\mcal{V}))$ for the category of $A(\mcal{V})$-modules. Let us note that if $M$ is a positive energy $\mcal{V}$-module, then we have $\Omega_\mcal{V}(M) \supset M_{\rm top}$ and $\Omega_\mcal{V}(M) = M_{\rm top}$ provided $M$ is a simple $\mcal{V}$-module.

Therefore, we may consider a functor $\Omega_\mcal{V} \colon \mcal{E}_+(\mcal{V}) \rarr \mcal{M}(A(\mcal{V}))$. On the other hand, there exists also an induction functor
\begin{align}
  \mathbb{M}_\mcal{V} \colon \mcal{M}(A(\mcal{V})) \rarr \mcal{E}_+(\mcal{V})
\end{align}
which is a left adjoint functor to $\Omega_\mcal{V}$ and has the following universal property. For a $\mcal{V}$-module $M$ and a morphism $\varphi \colon  E \rarr \Omega_\mcal{V}(M)$ of $A(\mcal{V})$-modules, there is a unique morphism $\widetilde{\varphi} \colon \mathbb{M}_\mcal{V}(E) \rarr M$ of $\mcal{V}$-modules which extends $\varphi$, see \cite{Zhu1996}, \cite{Dong-Li-Mason1998}. Moreover, for an $A(\mcal{V})$-module $E$ we have $\mathbb{M}_\mcal{V}(E)_{\rm top} \simeq E$ as modules over $A(\mcal{V})$. Besides, since the $\mcal{V}$-module $\mathbb{M}_\mcal{V}(E)$ has a unique maximal $\mcal{V}$-submodule $\mathbb{K}_\mcal{V}(E)$ having zero intersection with the $A(\mcal{V})$-submodule $E$ of $\mathbb{M}_\mcal{V}(E)$, we may set
\begin{align}
  \mathbb{L}_\mcal{V}(E) = \mathbb{M}_\mcal{V}(E)/\mathbb{K}_\mcal{V}(E)
\end{align}
for an $A(\mcal{V})$-module $E$.


\subsection{Affine vertex algebras}
\label{subsec:affine vertex}

Let $\mfrak{g}$ be a finite-dimensional reductive Lie algebra and let $\kappa$ be a $\mfrak{g}$-invariant symmetric bilinear form on $\mfrak{g}$. The induced $\widehat{\mfrak{g}}_\kappa$-module $\mathbb{M}_{\kappa,\mfrak{g}}(\C)$, where $\C$ is the trivial $1$-dimensional $\mfrak{g}$-module, is of a special importance in the theory of vertex algebras since it is equipped with the natural structure of an $\N_0$-graded vertex algebra, called the \emph{universal affine vertex algebra}, see \cite{Kac1998}, which we will denote by $\mcal{V}^\kappa(\mfrak{g})$. For an element $a \in \mfrak{g}$, we denote by $a(z) \in \widehat{\mfrak{g}}_\kappa[[z^{\pm 1}]]$ the formal distribution defined by
\begin{align}
a(z) = \sum_{n \in \Z} a_n z^{-n-1}.  \label{eq:field a}
\end{align}
The state-field correspondence $Y \colon \mcal{V}^\kappa(\mfrak{g}) \rarr \End \mcal{V}^\kappa(\mfrak{g})[[z^{\pm 1}]]$ is given by
\begin{align*}
  Y(a_{1,-n_1-1} \dots a_{k,-n_k-1}\vac,z) = {1 \over n_1! \dots n_k!}\, \normOrd{\partial_z^{n_1} a_1(z) \dots \partial_z^{n_k}a_k(z)}
\end{align*}
for $k \in \N$, $n_1,n_2,\dots,n_k \in \N_0$ and $a_1,a_2,\dots,a_k \in \mfrak{g}$, where $\vac$ is the vacuum vector (a highest weight vector of $\mathbb{M}_{\kappa,\mfrak{g}}(\C)$, i.e.\ we set $\vac = 1 \otimes 1$). The translation operator $T \colon \mcal{V}^\kappa(\mfrak{g}) \rarr \mcal{V}^\kappa(\mfrak{g})$ is determined by $T\vac = 0$ and $[T,a_n] = -na_{n-1}$ for $a \in \mfrak{g}$ and $n \in \Z$. Further, we define the $\N_0$-grading on $\mcal{V}^\kappa(\mfrak{g})$ by the formula $\deg \vac = 0$ and $\deg a_n = -n$ for $n \in \Z$, $a \in \mfrak{g}$. Moreover, for the commutator of fields we have
\begin{align}
  [a(z),b(w)] = [a,b](w)\delta(z-w) + \kappa(a,b)\partial_w \delta(z-w) \label{eq:comm relations vertex algebra}
\end{align}
for $a, b \in \mfrak{g}$.

For the Zhu's algebra $A(\mcal{V}^\kappa(\mfrak{g}))$ we have a canonical isomorphism
\begin{align*}
  A(\mcal{V}^\kappa(\mfrak{g})) \simeq U(\mfrak{g})
\end{align*}
of associative algebras determined by
\begin{align}
  a_{1,-n_1-1}a_{2,-n_2-1}\dots a_{k,-n_k-1}\vac \mapsto
  (-1)^{n_1+n_2+\dots +n_k} a_k\dots a_2a_1
\end{align}
for $k \in \N$, $n_1,n_2,\dots,n_k \in \N_0$ and $a_1,a_2,\dots,a_k \in \mfrak{g}$, see \cite{Frenkel-Zhu1992}. On the other hand, for the Zhu's $C_2$-algebra $R_{\mcal{V}^\kappa(\mfrak{g})}$ we have an isomorphism
\begin{align*}
  R_{\mcal{V}^\kappa(\mfrak{g})} \simeq S(\mfrak{g})
\end{align*}
of Poisson algebras given by
\begin{align}
   a_{1,-n_1-1}a_{2,-n_2-1}\dots a_{k,-n_k-1}\vac \mapsto \delta_{n_1,0}\delta_{n_2,0}\dots \delta_{n_k,0} a_1a_2\dots a_k
\end{align}
for $k \in \N$, $n_1,n_2,\dots,n_k \in \N_0$ and $a_1,a_2,\dots,a_k \in \mfrak{g}$, where the Poisson algebra structure on $S(\mfrak{g})$ is induced from the Kirillov--Kostant--Souriau Poisson structure on $\mfrak{g}^*$. Hence, for the associated variety $X_{\mcal{V}^\kappa(\mfrak{g})}$ of $\mcal{V}^\kappa(\mfrak{g})$ we get
\begin{align*}
  X_{\mcal{V}^\kappa(\mfrak{g})} \simeq \mfrak{g}^*.
\end{align*}
Besides, the mapping $\eta_{\mcal{V}^\kappa(\mfrak{g})}$ is an isomorphism of Poisson algebras.
\medskip

Let us consider a graded vertex algebra ideal $\mcal{J}_\kappa(\mfrak{g})$ of $\mcal{V}^\kappa(\mfrak{g})$ not containing the vacuum vector $\vac$. Then for the Zhu's algebra $A(\mcal{V}^\kappa(\mfrak{g})/\mcal{J}_\kappa(\mfrak{g}))$ we have
\begin{align*}
  A(\mcal{V}^\kappa(\mfrak{g})/\mcal{J}_\kappa(\mfrak{g})) \simeq U(\mfrak{g})/J_\kappa(\mfrak{g}),
\end{align*}
where the two-sided ideal $J_\kappa(\mfrak{g})$ of $U(\mfrak{g})$ is the image of $\mcal{J}_\kappa(\mfrak{g})$ in $A(\mcal{V}^\kappa(\mfrak{g})) \simeq U(\mfrak{g})$, see \cite[Proposition 1.4.2]{Frenkel-Zhu1992}. In many cases, the vertex algebra ideal $\mcal{J}_\kappa(\mfrak{g})$ is generated by singular vectors. If we identify the Lie algebra $\mfrak{g}$ with the Lie subalgebra $\mfrak{g} \otimes_\C \C 1$ of $\widehat{\mfrak{g}}_\kappa$ and endow $S(\mfrak{g})$ with the adjoint action of $\mfrak{g}$, then the projection $\sigma_{\rm Zhu} \colon \mcal{V}^\kappa(\mfrak{g}) \rarr R_{\mcal{V}^\kappa(\mfrak{g})} \simeq S(\mfrak{g})$ is a homomorphism of $\mfrak{g}$-modules. Hence, the image of a singular vector in $\mcal{V}^\kappa(\mfrak{g})$ under $\sigma_{\rm Zhu}$ is a singular vector in $S(\mfrak{g})$. Further, let us introduce the symmetrization mapping
\begin{align*}
  \beta_{-1} \colon S(\mfrak{g}) \rarr \mcal{V}^\kappa(\mfrak{g})
\end{align*}
by
\begin{align*}
  \beta_{-1}(a_1a_2\dots a_k) = {1 \over k!} \sum_{\sigma \in \mfrak{S}_k} a_{\sigma(1),-1}a_{\sigma(2),-1} \dots a_{\sigma(k),-1}\vac
\end{align*}
for $k \in \N$ and $a_1,a_2,\dots,a_k \in \mfrak{g}$. Then $\beta_{-1}$ is an injective homomorphism of $\mfrak{g}$-modules such that
\begin{align*}
  \sigma_{\rm Zhu} \circ \beta_{-1} = \id_{S(\mfrak{g})}.
\end{align*}
If we assume that $\mfrak{g}$ is simple, then by \cite[Theorem 1.2]{Arakawa-Jiang-Moreau2021} we have that $\sigma_{\rm Zhu}(v) \in S(\mfrak{g})$ is a non-zero singular vector whenever $v \in \mcal{V}^\kappa(\mfrak{g})$ is a non-zero singular vector provided $\kappa$ is a non-critical level. On the other hand, if $w \in S(\mfrak{g})$ is a singular vector, then $\beta_{-1}(w) \in \mcal{V}^\kappa(\mfrak{g})$ is a singular vector if and only if $f_{\theta,1}\beta_{-1}(w)=0$. We call $\sigma_{\rm Zhu}(v)$ the symbol of a singular vector $v \in \mcal{V}^\kappa(\mfrak{g})$.

In particular, the universal affine vertex algebra $\mcal{V}^\kappa(\mfrak{g})$ has a unique maximal graded vertex algebra ideal $\mcal{I}_\kappa(\mfrak{g})$ not containing the vacuum vector. The $\N_0$-graded quotient vertex algebra $\mcal{V}^\kappa(\mfrak{g})/\mcal{I}_\kappa(\mfrak{g})$ is
called the \emph{simple affine vertex algebra} and we will denote it by $\mcal{V}_\kappa(\mfrak{g})$.
\medskip

Let us note that for the universal affine vertex algebra $\mcal{V}^\kappa(\mfrak{g})$ the functors $\mathbb{M}_{\mcal{V}^\kappa(\mfrak{g})}$ and $\mathbb{L}_{\mcal{V}^\kappa(\mfrak{g})}$ coincide with the functors $\mathbb{M}_{\kappa,\mfrak{g}}$ and $\mathbb{L}_{\kappa,\mfrak{g}}$, respectively. Therefore, according to Theorem \ref{thm:Zhu correspondence} the assignment $E \mapsto \mathbb{L}_{\kappa,\mfrak{g}}(E)$ gives a one-to-one correspondence between the isomorphism classes of simple $\mfrak{g}$-modules and simple positive energy $\mcal{V}^\kappa(\mfrak{g})$-modules and also provides a one-to-one correspondence between isomorphism classes of simple modules over $U(\mfrak{g})/I_\kappa(\mfrak{g})$ and simple positive energy $\mcal{V}_\kappa(\mfrak{g})$-modules.


\subsection{Weyl vertex algebras}
\label{subsec:Weyl vertex}

Let $V$ be a finite-dimensional vector space. The induced $\eus{A}_{\mcal{K}(V)}$-module $\smash{\mathbb{M}_{\mcal{K}(V)}^\eus{A}(S(V^*))}$ carries the natural structure of an $\N_0$-graded vertex algebra, called the \emph{Weyl vertex algebra} (or \emph{$\beta\gamma$-system}), which we denote by $\mcal{M}_V$. Let $\{x_1,x_2,\dots,x_m\}$, where $m=\dim V$, be linear coordinate functions on $V$ and let $a_i(z), a_i^*(z) \in \eus{A}_{\mcal{K}(V)}[[z^{\pm 1}]]$ for $i=1,2,\dots,m$ be the formal distributions defined by
\begin{align}
  a_i(z)= \sum_{n\in \Z} a_{i,n}z^{-n-1} \qquad \text{and} \qquad a_i^*(z) = \sum_{n\in \Z} a^*_{i,n}z^{-n}, \label{eq:field a,a*}
\end{align}
where $a_{i,n}=\partial_{x_{i,n}}$ and $a^*_{i,n}=x_{i,-n}$ for $n \in \Z$ and $i=1,2,\dots,m$.  The state-field correspondence $Y \colon \mcal{M}_V \rarr \End \mcal{M}_V[[z^{\pm 1}]]$ is given by
\begin{multline*}
  Y(a_{i_1,-n_1-1}\dots a_{i_r,-n_r-1}a^*_{j_1,-m_1}\dots a^*_{j_s,-m_s}\vac,z) = \\ {1 \over n_1! \dots n_r!} {1 \over m_1!\dots m_s!} \, \normOrd{\partial_z^{n_1}a_{i_1}\!(z) \dots \partial_z^{n_r}a_{i_r}\!(z) \partial_z^{m_1}a_{j_1}^*\!(z) \dots \partial_z^{m_s}a_{j_s}^*\!(z)}
\end{multline*}
for $r,s \in \N$ and $n_1,\dots,n_r,m_1,\dots,m_s \in \N_0$, where $\vac$ is the vacuum vector (a highest weight vector of $\smash{\mathbb{M}_{\mcal{K}(V)}^\eus{A}(S(V^*))}$, i.e.\ we set $\vac = 1 \otimes 1$). The translation operator $T \colon \mcal{M}_V \rarr \mcal{M}_V$ is defined by $T\vac = 0$, $[T,a_{i,n}] = -na_{i,n-1}$ and $[T,a^*_{i,n}] = -(n-1)a^*_{i,n-1}$ for $n \in \Z$ and $i=1,2,\dots,m$. We define the $\N_0$-grading on $\mcal{M}_V$ by the formula $\deg \vac = 0$ and $\deg \partial_{x_{i,n}} = -n$, $\deg x_{i,n} = n$ for $n \in \Z$, $i=1,2,\dots,m$. Moreover, for the commutator of fields we have
\begin{align}
  [a_i(z), a_j(w)] = 0, \qquad [a_i(z), a_j^*(w)] = \delta_{i,j} \delta(z-w), \qquad [a_i^*(z), a_j^*(w)] = 0
\end{align}
for $i,j = 1,2,\dots,m$.

For the Zhu's algebra $A(\mcal{M}_V)$, we have a canonical isomorphism
\begin{align*}
  A(\mcal{M}_V) \simeq \eus{A}_V
\end{align*}
of associative algebras determined by
\begin{multline*}
  a_{i_1,-n_1-1}\dots a_{i_r,-n_r-1}a^*_{j_1,-m_1}\dots a^*_{j_s,-m_s}\vac \mapsto \\
  \delta_{m_1,0} \dots \delta_{m_s,0}(-1)^{n_1+\dots + n_r}x_{j_1}\dots x_{j_s}\partial_{x_{i_1}}\dots\partial_{x_{i_r}}
\end{multline*}
for $r,s \in \N$, $n_1,\dots,n_r,m_1,\dots,m_s \in \N_0$. On the other hand, for the Zhu's $C_2$-algebra $R_{\mcal{M}_V}$ we have an isomorphism
\begin{align*}
  R_{\mcal{M}_V} \simeq \C[V \oplus V^*]
\end{align*}
of Poisson algebras given by
\begin{multline*}
  a_{i_1,-n_1-1}\dots a_{i_r,-n_r-1}a^*_{j_1,-m_1}\dots a^*_{j_s,-m_s}\vac \mapsto \\
    \delta_{n_1,0} \dots \delta_{n_r,0}\delta_{m_1,0} \dots \delta_{m_s,0} x_{j_1}\dots x_{j_s}\xi_{i_1}\dots \xi_{i_r}
\end{multline*}
for $r,s \in \N$, $n_1,\dots,n_r,m_1,\dots,m_s \in \N_0$, where $\xi_1,\xi_2,\dots,\xi_m$ are linear coordinate functions on $V^*$ dual to $x_1,x_2,\dots,x_m$. The Poisson bracket on $\C[V \oplus V^*]$ is defined by
\begin{align*}
  \{f,g\} = \sum_{i=1}^m \!\bigg({\partial f \over \partial \xi_i} {\partial g \over \partial x_i} - {\partial g \over \partial \xi_i} {\partial f \over \partial x_i}\bigg)
\end{align*}
for $f,g \in \C[V\oplus V^*]$. Hence, for the associated variety $X_{\mcal{M}_V}$ of $\mcal{M}_V$ we get
\begin{align*}
  X_{\mcal{M}_V} \simeq V\oplus V^*.
\end{align*}
Besides, the mapping $\eta_{\mcal{M}_V}$ is an isomorphism of Poisson algebras.
\medskip

Let us note that for the Weyl vertex algebra $\mcal{M}_V$ the functors $\mathbb{M}_{\mcal{M}_V}$ and $\mathbb{L}_{\mcal{M}_V}$ coincide with the functors $\smash{\mathbb{M}_{\mcal{K}(V)}^\eus{A}}$ and $\smash{\mathbb{L}_{\mcal{K}(V)}^\eus{A}}$, respectively.


\subsection{Clifford vertex superalgebras}
\label{subsec:Clifford vertex}

Let $V$ be a finite-dimensional vector space. Then the induced $\eus{C}_{\mcal{K}(V)}$-module $\smash{\mathbb{M}_{\mcal{K}(V)}^\eus{C}(\Lambda(V^*))}$ carries the natural structure of an $\N_0$-graded vertex superalgebra, called the \emph{Clifford vertex superalgebra}, which we will denote by $\smash{\bigwedge^{\infty/2+\bullet}_V}$. Let $\{\theta_1,\theta_2,\dots,\theta_m\}$, where $m=\dim V$, be linear coordinate functions on $V$ and let $\psi_i(z), \psi_i^*(z) \in \eus{C}_{\mcal{K}(V)}[[z^{\pm 1}]]$ for $i=1,2,\dots,m$ be the formal distributions defined by
\begin{align}
  \psi_i(z)= \sum_{n\in \Z} \psi_{i,n}z^{-n-1} \qquad \text{and} \qquad \psi_i^*(z) = \sum_{n\in \Z} \psi^*_{i,n}z^{-n}, \label{eq:field psi,psi*}
\end{align}
where $\psi_{i,n}=\partial_{\theta_{i,n}}$ and $\psi^*_{i,n}=\theta_{i,-n}$ for $n \in \Z$ and $i=1,2,\dots,m$.  The state-field correspondence $Y \colon \smash{\bigwedge^{\infty/2+\bullet}_V} \rarr \End \smash{\bigwedge^{\infty/2+\bullet}_V}[[z^{\pm 1}]]$ is given by
\begin{multline*}
  Y(\psi_{i_1,-n_1-1}\dots \psi_{i_r,-n_r-1}\psi^*_{j_1,-m_1}\dots \psi^*_{j_s,-m_s}\vac,z) = \\ {1 \over n_1! \dots n_r!} {1 \over m_1!\dots m_s!} \, \normOrd{\partial_z^{n_1}\psi_{i_1}\!(z) \dots \partial_z^{n_r}\psi_{i_r}\!(z) \partial_z^{m_1}\psi_{j_1}^*\!(z) \dots \partial_z^{m_s}\psi_{j_s}^*\!(z)}
\end{multline*}
for $r,s \in \N$ and $n_1,\dots,n_r,m_1,\dots,m_s \in \N_0$, where $\vac$ is the vacuum vector (a highest weight vector of $\smash{\mathbb{M}_{\mcal{K}(V)}^\eus{C}(\Lambda(V^*))}$, i.e.\ we set $\vac = 1 \otimes 1$). The translation operator $T \colon \smash{\bigwedge^{\infty/2+\bullet}_V} \rarr \smash{\bigwedge^{\infty/2+\bullet}_V}$ is defined by $T\vac = 0$, $[T,\psi_{i,n}] = -n\psi_{i,n-1}$ and $[T,\psi^*_{i,n}] = -(n-1)\psi^*_{i,n-1}$ for $n \in \Z$ and $i=1,2,\dots,m$. We define the $\N_0$-grading on $\smash{\bigwedge^{\infty/2+\bullet}_V}$ by the formula $\deg \vac = 0$ and $\deg \partial_{\theta_{i,n}} = -n$, $\deg \theta_{i,n} = n$ for $n \in \Z$, $i=1,2,\dots,m$. We introduce an additional $\Z$-gradation, called the \emph{charge gradation}, by setting $\cdeg \vac = 0$ and $\cdeg \partial_{\theta_{i,n}} = -1$, $\cdeg \theta_{i,n} = 1$ for $n \in \Z$, $i=1,2,\dots,m$. Moreover, for the anti-commutator of fields we have
\begin{align}
  [\psi_i(z), \psi_j(w)]_+ = 0, \qquad [\psi_i(z), \psi_j^*(w)]_+ = \delta_{i,j} \delta(z-w), \qquad [\psi_i^*(z), \psi_j^*(w)]_+ = 0
\end{align}
for $i,j = 1,2,\dots,m$. For the Clifford vertex algebra $\smash{\bigwedge\nolimits^{\infty/2+\bullet}_V}$, we have a canonical isomorphism
\begin{align*}
A({\textstyle \bigwedge\nolimits^{\infty/2+\bullet}_V}) \simeq \eus{C}_V \qquad \text{and} \qquad
\widetilde{X}_{\smash{\bigwedge^{\infty/2+\bullet}_V}} \simeq  T^*\Pi V,
\end{align*}
where $\eus{C}_V$ is the Clifford algebra associated with $V \oplus V^*$ and $\Pi V$
is the vector space $V$ considered as a purely odd affine space.


\subsection{Vertex algebras of neutral free bosons and superfermions}
\label{subsec:neutral free}

Let $\mfrak{a}$ be a finite-dimensional commutative Lie algebra and let $\omega$ be a bilinear form on $\mfrak{a}$. Then we may introduce either the Kac--Moody affinization of $\mfrak{a}$ if $\omega$ is symmetric or the Clifford affinization of $\mfrak{a}$ if $\omega$ is skew-symmetric. Let us consider a $1$-dimensional central extension $\widehat{\mfrak{a}}_\omega = \mfrak{a}(\!(t)\!) \oplus \C c$ of the formal loop algebra $\mfrak{a}(\!(t)\!) = \mfrak{a} \otimes_\C \C(\!(t)\!)$ with the commutation relations
\begin{align*}
  [a\otimes f(t), b\otimes g(t)] =
  \begin{cases}
    -\omega(a,b) \Res_{t=0}(f(t)dg(t)) c & \text{if $\omega$ is symmetric}, \\
    \omega(a,b) \Res_{t=0}(f(t)g(t)dt) c & \text{if $\omega$ is skew-symmetric}
  \end{cases}
\end{align*}
for $a,b \in \mfrak{a}$ and $f(t), g(t) \in \C(\!(t)\!)$, where $c$ is the central element of $\widehat{a}_\omega$. By introducing the notation $a_n = a \otimes t^n$ for $a \in \mfrak{a}$ and $n \in \Z$, the previous commutation relations can be simplified into the form
\begin{align*}
  [a_m,b_n] =
  \begin{cases}
    m\omega(a,b) \delta_{m,-n} c & \text{if $\omega$ is symmetric}, \\
    \omega(a,b) \delta_{m,-n-1} c & \text{if $\omega$ is skew-symmetric}
  \end{cases}
\end{align*}
for $m,n \in \Z$ and $a,b \in \mfrak{a}$. Let us note that $\widehat{\mfrak{a}}_\omega$ has the structure of a $\Z$-graded topological Lie algebra with the gradation defined by $\deg c =0$ and $\deg a_n = -n$ for $a\in \mfrak{a}$, $n\in \Z$.

Since $\mfrak{a}[[t]] \oplus \C c$ is a Lie subalgebra of $\widehat{\mfrak{a}}_\omega$, we define the induced $\widehat{\mfrak{a}}_\omega$-module
\begin{align*}
  \mcal{F}^\omega(\mfrak{a}) = U(\widehat{\mfrak{a}}_\omega) \otimes_{U(\mfrak{a}[[t]] \oplus \C c)}\! \C,
\end{align*}
where $\C$ is considered as a module over $\mfrak{a}[[t]] \oplus \C c$ on which $\mfrak{a}[[t]]$ acts trivially and $c$ acts as the identity. The induced $\widehat{\mfrak{a}}_\omega$-module $\mcal{F}^\omega(\mfrak{a})$ is equipped with the natural structure of an $\N_0$-graded vertex algebra, called the \emph{vertex algebra of neutral free bosons} if $\omega$ is symmetric and the \emph{vertex algebra of neutral free superfermions} if $\omega$ is skew-symmetric. For an element $a \in \mfrak{a}$, we denote by $a(z) \in \widehat{\mfrak{a}}_\omega[[z^{\pm 1}]]$ the formal distribution defined by
\begin{align}
a(z) = \sum_{n \in \Z} a_n z^{-n-1}.
\end{align}
The state-field correspondence $Y \colon \mcal{F}^\omega(\mfrak{a}) \rarr \End \mcal{F}^\omega(\mfrak{a})[[z^{\pm 1}]]$ is given by
\begin{align*}
  Y(a_{1,-n_1-1} \dots a_{k,-n_k-1}\vac,z) = {1 \over n_1! \dots n_k!}\, \normOrd{\partial_z^{n_1} a_1(z) \dots \partial_z^{n_k}a_k(z)}
\end{align*}
for $k \in \N$, $n_1,n_2,\dots,n_k \in \N_0$ and $a_1,a_2,\dots,a_k \in \mfrak{a}$, where $\vac \in \mcal{F}^\omega(\mfrak{a})$ is the vacuum vector (a highest weight vector of $\mcal{F}^\omega(\mfrak{a})$, i.e.\ we set $\vac = 1 \otimes 1$). The translation operator $T \colon \mcal{F}^\omega(\mfrak{a}) \rarr \mcal{F}^\omega(\mfrak{a})$ is determined by $T\vac = 0$ and $[T,a_n] = -na_{n-1}$ for $a \in \mfrak{a}$ and $n \in \Z$. Further, we define the $\N_0$-grading on $\mcal{F}^\omega(\mfrak{a})$ by the formula $\deg \vac = 0$ and $\deg a_n = -n$ for $n \in \Z$, $a \in \mfrak{a}$. Moreover, for the commutator of fields we have
\begin{align}
  [a(z),b(w)] =
  \begin{cases}
    \omega(a,b)\partial_w \delta(z-w) & \text{if $\omega$ is symmetric}, \\
    \omega(a,b)\delta(z-w) & \text{if $\omega$ is skew-symmetric}
  \end{cases}
\end{align}
for $a, b \in \mfrak{a}$.


\section{Associated varieties of affine vertex algebras}

We show that the Dixmier sheet closures appear as associated varieties of quotients of universal affine vertex algebras and determine associated varieties of simple affine vertex algebras for certain non-admissible levels.


\subsection{Dixmier sheets and Richardson orbits}

Let $G$ be a connected reductive algebraic group with its Lie algebra $\mfrak{g}$. We denote by $\mcal{N}(\mfrak{g})$ the nilpotent cone of $\mfrak{g}$, i.e.\ the set of nilpotent elements of $\mfrak{g}$. It is an irreducible closed algebraic subvariety of $\mfrak{g}$ and a finite union of $G$-orbits. There is a unique nilpotent orbit of $\mfrak{g}$, denoted by $\mcal{O}_{\rm reg}$ and called the \emph{regular nilpotent orbit} of $\mfrak{g}$, which is a dense open subset of $\mcal{N}(\mfrak{g})$. When $\mfrak{g}$ is a simple Lie algebra, there exists a unique nilpotent orbit of $\mfrak{g}$ that is a dense open subset of $\mcal{N}(\mfrak{g}) \setminus \mcal{O}_{\rm reg}$, denoted by $\mcal{O}_{\rm subreg}$ and called the \emph{subregular nilpotent orbit} of $\mfrak{g}$. Besides, there is a unique non-zero nilpotent orbit of $\mfrak{g}$ of minimal dimension, denoted by $\mcal{O}_{\rm min}$ and called the \emph{minimal nilpotent orbit} of $\mfrak{g}$, such that it is contained in the closure of all non-zero nilpotent orbits of $\mfrak{g}$. By $\mcal{O}_{\rm zero}$ we denote the \emph{zero nilpotent orbit} of $\mfrak{g}$.
For the dimension of these distinguished nilpotent orbits of $\mfrak{g}$ see Figure \ref{fig:hasse diagram general}.

\begin{figure}[ht]
\centering
{\begin{tikzpicture}
[yscale=1.1,xscale=1.7,vector/.style={circle,draw=white,fill=black,ultra thick, inner sep=0.8mm},vector2/.style={circle,draw=white,fill=white,ultra thick, inner sep=1mm}]
\begin{scope}
  \node (A) at (0,0)  {$\mcal{O}_{{\rm zero}}$};
  \node (B) at (0,1)  {$\mcal{O}_{{\rm min}}$};
  \node (C1) at (0,1.7)  {};
  \node (C) at (0,2)  {$\,\dots$};
  \node (C2) at (0,2.3) {};
  \node (D) at (0,3)  {$\mcal{O}_{{\rm subreg}}$};
  \node (E) at (0,4)  {$\mcal{O}_{{\rm reg}}$};
  \node at (2,0) {$0$};
  \node at (2,1) {$2h^\vee - 2$};
  \node at (2,2) {$\dots$};
  \node at (2,3) {$\dim \mfrak{g} - \rank \mfrak{g} - 2$};
  \node at (2,4) {$\dim \mfrak{g} - \rank \mfrak{g}$};
  \draw [thin, -] (B) -- (C1);
  \draw [thin, -] (D) -- (C2);
  \draw [thin, -] (A) -- (B);
  \draw [thin, -] (D) -- (E);
  \node at (0,4.7) {nilpotent orbit};
  \node at (2,4.7) {dimension};
\end{scope}
\end{tikzpicture}}
\caption{Hasse diagram of nilpotent orbits of a simple Lie algebra $\mfrak{g}$}
\label{fig:hasse diagram general}
\vspace{-2mm}
\end{figure}

For $x \in \mfrak{g}$, we denote by $\mfrak{g}^x$ the centralizer of $x$ in $\mfrak{g}$. Besides, for a subset $X$ of $\mfrak{g}$ we define the set
\begin{align*}
  X^{\rm reg} = \{x \in X;\, \dim \mfrak{g}^x = \min\nolimits_{y \in X} \dim \mfrak{g}^y \}
\end{align*}
and call it the \emph{set of regular elements} in $X$. The \emph{Jordan class} of an element $x \in \mfrak{g}$ can be defined by
\begin{align*}
  J_G(x) =G.(\mfrak{z}(\mfrak{g}^{x_s})^{\rm reg} +x_n),
\end{align*}
where $x=x_s+x_n$ is the Jordan decomposition of $x$ with the semisimple part $x_s$ and the nilpotent part $x_n$. It is a $G$-invariant, irreducible and locally closed subset of $\mfrak{g}$, hence it is a subvariety of $\mfrak{g}$. A \emph{sheet} of $\mfrak{g}$ is an irreducible component of the subset
\begin{align*}
  \mfrak{g}^{(m)} = \{x \in \mfrak{g};\, \dim \mfrak{g}^x =m\}
\end{align*}
for some $m \in \N_0$. It is a finite disjoint union of Jordan classes and contains a unique nilpotent orbit. A sheet of $\mfrak{g}$ which contains a semisimple element is called a \emph{Dixmier sheet} of $\mfrak{g}$ and the unique nilpotent orbit in the sheet is a \emph{Richardson orbit} of $\mfrak{g}$. If $\mcal{S}$ is a Dixmier sheet of $\mfrak{g}$, then there exists a parabolic subalgebra $\mfrak{p}$ of $\mfrak{g}$ such that
\begin{align*}
  \widebar{\mcal{S}} =  G.[\mfrak{p},\mfrak{p}]^\perp, \qquad \mcal{S} = (G.[\mfrak{p},\mfrak{p}]^\perp)^{\rm reg}
\end{align*}
and the corresponding Richardson orbit is given by
\begin{align*}
  \widebar{\mcal{O}} = G.\mfrak{p}^\perp, \qquad \mcal{O} = (G.\mfrak{p}^\perp)^{\rm reg},
\end{align*}
where $[\mfrak{p},\mfrak{p}]^\perp$ and $\mfrak{p}^\perp$ denote the orthogonal complements of $[\mfrak{p},\mfrak{p}]$ and $\mfrak{p}$, respectively, according to a fixed non-degenerate $\mfrak{g}$-invariant symmetric bilinear form $\kappa_\omega$ on $\mfrak{g}$. The Dixmier sheet and the Richardson orbit given by a parabolic subalgebra $\mfrak{p}$ of $\mfrak{g}$ we will denote by $\mcal{S}_\mfrak{p}$ and $\mcal{O}_\mfrak{p}$, respectively. Let us note that we have
\begin{align}
  \mcal{S}_\mfrak{p} \cap \mcal{N}(\mfrak{g}) = \mcal{O}_\mfrak{p} \qquad \text{and} \qquad \widebar{\mcal{S}_{\mfrak{p}}} \cap \mcal{N}(\mfrak{g}) = \widebar{\mcal{O}_{\mfrak{p}}} \label{eq:intersection}
\end{align}
by \cite[Corrolary 39.2.6, Proposition 39.3.3]{Tauvel-Yu2005-book}.
\medskip

Let us consider the PBW filtration on the universal enveloping algebra $U(\mfrak{g})$ of $\mfrak{g}$ and the associated graded algebra $\gr U(\mfrak{g}) \simeq S(\mfrak{g}) \simeq \C[\mfrak{g}^*]$. The \emph{associated variety} $\mcal{V}(I)$ of a left ideal $I$ of $U(\mfrak{g})$ is defined as the zero locus in $\mfrak{g}^*$ of the associated graded ideal $\gr I$ of $S(\mfrak{g})$. Moreover, if $I$ is a two-sided ideal of $\mfrak{g}$, then $I$ and $\gr I$ are invariant under the adjoint action of $G$. Consequently, the associated variety $\mcal{V}(I)$ is a union of $G$-orbits of $\mfrak{g}^*$. Obviously, we have
\begin{align}
  \mcal{V}(I) = \Specm (S(\mfrak{g})/\gr I) = \Specm(S(\mfrak{g})/\sqrt{\gr I}),
\end{align}
where $\sqrt{\gr I}$ stands for the radical of $\gr I$. As $\kappa_\omega$ is a non-degenerate $\mfrak{g}$-invariant symmetric bilinear form on $\mfrak{g}$, it provides a one-to-one correspondence between adjoint orbits of $\mfrak{g}$ and coadjoint orbits of $\mfrak{g}^*$. For an adjoint orbit $\mcal{O}$ of $\mfrak{g}$ we denote by $\mcal{O}^*$ the corresponding coadjoint orbit of $\mfrak{g}^*$, and similarly for a sheet $\mcal{S}$ of $\mfrak{g}$ we denote by $\mcal{S}^*$ the corresponding subset of $\mfrak{g}^*$.


\subsection{Momentum mapping of homogeneous spaces}
\label{subsec:momentum}

Let $G$ be a connected reductive algebraic group, $H$ be a maximal torus of $G$ and $P$ be a connected parabolic subgroup of $G$ containing $H$ with the Levi subgroup $L$, the unipotent radical $U$ and the opposite unipotent radical $\widebar{U}$. Let us recall that $P=LU$. We denote by $\mfrak{g}$, $\mfrak{p}$, $\mfrak{l}$, $\mfrak{u}$, $\widebar{\mfrak{u}}$ and $\mfrak{h}$ the Lie algebras of $G$, $P$, $L$, $U$, $\widebar{U}$ and $H$, respectively. The homogenous space
\begin{align*}
  X = G/P
\end{align*}
is a smooth algebraic variety, the \emph{generalized flag variety} of $G$. Besides, we have also the canonical $G$-equivariant projection
\begin{align*}
  p \colon G \rarr G/P.
\end{align*}
By a result of \cite{Kashiwara1989}, for any $\lambda \in \mfrak{z}(\mfrak{l})^* \subset \mfrak{h}^*$ there is a $G$-equivariant sheaf of rings of twisted differential operators $\mcal{D}_X^\lambda$ on $X$. The normalization is such that $\smash{\mcal{D}_X^{-\rho_\mfrak{p}}}$ represents the usual $G$-equivariant sheaf of rings of differential operators on $X$, where
\begin{align*}
  \rho_\mfrak{p}(a) = -{1\over 2}\tr_{\mfrak{g}/\mfrak{p}} \ad(a)
\end{align*}
for $a \in \mfrak{l}$. In fact, we have that $\rho_\mfrak{p}{}_{|[\mfrak{l},\mfrak{l}]}=0$, which means that $\rho_\mfrak{p} \in \mfrak{z}(\mfrak{l})^*$. As $\mcal{D}_X^\lambda$ is $G$-equivariant, we have a Lie algebra homomorphism $\mfrak{g} \rarr \Gamma(X,\mcal{D}_X^\lambda)$, which extends to a surjective homomorphism
\begin{align}
  \Phi^\lambda_X \colon U(\mfrak{g}) \rarr \Gamma(X,\mcal{D}_X^\lambda) \label{eq:Phi homomorphism}
\end{align}
of associative algebras. Hence, for any $\mcal{D}_X^\lambda$-module $\mcal{M}$ the vector space $\Gamma(X,\mcal{M})$ of global sections of $\mcal{M}$ has the natural structure of a $\mfrak{g}$-module. Further, let us consider the homogenous space
\begin{align*}
  Y = G/(P,P),
\end{align*}
where $(P,P)$ is the commutator subgroup of $P$, the \emph{generalized base affine space} of $G$. We have the canonical $G$-equivariant projection
\begin{align*}
  q \colon G \rarr G/(P,P).
\end{align*}
In addition, we have that $Y$ is a principal $P/(P,P)$-bundle over $X$. Let us note also that 
\begin{align*}
  A=P/(P,P) 
\end{align*}
is a torus. As $\mcal{D}_Y$ is a $G$-equivariant sheaf of rings of differential operators on $Y$, we have a Lie algebra homomorphism $\mfrak{g} \rarr \Gamma(Y,\mcal{D}_Y)$, which extends to a homomorphism
\begin{align}
  \Phi_Y \colon U(\mfrak{g}) \rarr \Gamma(Y,\mcal{D}_Y)
\end{align}
of associative algebras. It has been shown in \cite{Borho-Brylinski1982} that
\begin{align}
   \ker \Phi_Y = \bigcap_{\lambda \in \smash{\mfrak{z}(\mfrak{l})^*}} \ker \Phi^\lambda_X = \bigcap_{\lambda \in \smash{\mfrak{z}(\mfrak{l})^*}} \Ann_{U(\mfrak{g})}\! M^\mfrak{g}_\mfrak{p}(\lambda) \label{eq:ker Phi_Y}
\end{align}
and that the associated variety of the two-sided ideal $\ker \Phi_Y$ of $U(\mfrak{g})$ coincides with the closure of the Dixmier sheet attached to $\mfrak{p}$, i.e.\ we have
\begin{align}
  \mcal{V}(\ker \Phi_Y) = \widebar{\mcal{S}^*_{\mfrak{p}}}. \label{eq:ker Phi_Y variety}
\end{align}
As $\ker \Phi_Y$ is a two-sided ideal, the associated variety $\mcal{V}(\ker \Phi_Y)$ is a union of $G$-orbits. In addition, since $Y$ is a homogeneous $G$-space, there is the momentum mapping
\begin{align*}
  \mu_Y \colon T^*Y \rarr \mfrak{g}^*
\end{align*}
for the symplectic form $\omega_Y = d\alpha_Y$, where $\alpha_Y$ is the canonical $1$-form on $T^*Y$. Moreover, we have that
\begin{align}
  \langle \mu_Y, a \rangle = -\alpha_Y(\Phi_Y(a))
\end{align}
for $a \in \mfrak{g}$. Besides, the momentum mapping $\mu_Y$ induces a $G$-equivariant morphism
\begin{align*}
  \widebar{\mu}_Y \colon T^*Y/A \rarr \mfrak{g}^*.
\end{align*}
If we identify the vector bundle $T^*Y/A$ over $X$ with the associated vector bundle $G \times_P [\mfrak{p},\mfrak{p}]^\perp$ and $\mfrak{g}^*$ with $\mfrak{g}$ through $\kappa_\omega$, then $\widebar{\mu}_Y$ coincides with the \emph{generalized Grothendieck's simultaneous resolution}
\begin{align}
  \mu_P \colon G \times_P [\mfrak{p},\mfrak{p}]^\perp \rarr \mfrak{g} \label{eq:Grothendieck resolusion}
\end{align}
defined by
\begin{align*}
  \mu_P(g,a) = \Ad(g)(a)
\end{align*}
for $g \in G$ and $a \in [\mfrak{p},\mfrak{p}]^\perp$. Hence, it follows immediately that the image of $\widebar{\mu}_Y$ is $\smash{\widebar{\mcal{S}_{\mfrak{p}}^*}}$.


\subsection{Dynkin index}

Let us assume that $\mfrak{g}$ is a simple Lie algebra. For a finite-dimensional $\mfrak{g}$-module $\sigma_E \colon \mfrak{g} \rarr \End E$, we introduce the $\mfrak{g}$-invariant symmetric bilinear form $\kappa_E$ on $\mfrak{g}$ by
\begin{align*}
  \kappa_E(a,b) = \tr_E(\sigma_E(a)\sigma_E(b))
\end{align*}
for $a,b \in \mfrak{g}$. As $\mfrak{g}$ is a simple Lie algebra, any $\mfrak{g}$-invariant symmetric bilinear form on $\mfrak{g}$ is of the form $k\kappa_0$ for some $k \in \C$, where $\kappa_0$ is the normalized $\mfrak{g}$-invariant symmetric bilinear form on $\mfrak{g}$ such that
\begin{align*}
  \kappa_\mfrak{g} = 2h^\vee \kappa_0.
\end{align*}
The \emph{Dynkin index} of a finite-dimensional $\mfrak{g}$-module $E$, denoted by $\ind_D(\mfrak{g},E)$, is defined by
\begin{align}
  \kappa_E = \ind_D(\mfrak{g},E) \kappa_0.
  \label{eq:Dynkin-index}
\end{align}
A theorem of Dynkin \cite{Dynkin1952} states that if $L^\mfrak{g}_\mfrak{b}(\lambda)$ is the simple finite-dimensional $\mfrak{g}$-module with highest weight $\lambda \in \mfrak{h}^*$ with respect to a Borel subalgebra $\mfrak{b}$ of $\mfrak{g}$, then
\begin{align}
  \ind_D(\mfrak{g},L^\mfrak{g}_\mfrak{b}(\lambda)) =  {\dim L^\mfrak{g}_\mfrak{b}(\lambda) \over \dim \mfrak{g}}\,(\lambda,\lambda+2\rho), \label{eq:Dynkin index}
\end{align}
where $(\cdot\,,\cdot)$ is the bilinear form on $\mfrak{h}^*$ induced by $\kappa_0$. Although it is not obvious from the definition, the Dynkin index is a non-negative integer. Moreover, we have
\begin{align*}
  \ind_D(\mfrak{g},E_1 \oplus E_2) = \ind_D(\mfrak{g},E_1) + \ind_D(\mfrak{g},E_2)
\end{align*}
for finite-dimensional $\mfrak{g}$-modules $E_1$ and $E_2$. For the adjoint representation of $\mfrak{g}$, we get
\begin{align*}
  \ind_D(\mfrak{g},\mfrak{g}) = (\theta,\theta+2\rho) = 2(1+(\theta,\rho))= 2h^\vee,
\end{align*}
where $\theta$ is the maximal root of $\mfrak{g}$. For example, let $\mfrak{g}=\mfrak{sl}_n$ with $n\geq 2$. Then the Dynkin index of the simple finite-dimensional $\mfrak{g}$-module $L^\mfrak{g}_\mfrak{b}(a\omega_1)$ with $a \in \N_0$ is given by
\begin{align*}
  \ind_D(\mfrak{g},L^\mfrak{g}_\mfrak{b}(a\omega_1))= \binom{n+a}{n+1},
\end{align*}
which implies that the normalized $\mfrak{g}$-invariant symmetric bilinear form $\kappa_0$ on $\mfrak{g}$ is the trace form of the basic representation of $\mfrak{g}$.
\medskip

Let $\mfrak{p} = \mfrak{l} \oplus \mfrak{u}$ be a parabolic subalgebra of $\mfrak{g}$ with the Levi subalgebra $\mfrak{l}$, the nilradical $\mfrak{u}$ and the opposite nilradical $\widebar{\mfrak{u}}$. Since $\mfrak{l}$ is a reductive Lie algebra, we have the direct sum decomposition
\begin{align}
  \mfrak{l} = \bigoplus_{i=0}^r \mfrak{l}_i \label{eq:othogonal decomposition}
\end{align}
of $\mfrak{l}$ into the direct sum of simple Lie subalgebras $\mfrak{l}_i$ for $i=1,2,\dots,r$ and an abelian Lie subalgebra $\mfrak{l}_0$ such that these direct summands are mutually orthogonal with respect to the Cartan--Killing form $\kappa_\mfrak{g}$. We denote by $\smash{\kappa_0^{\mfrak{l}_i}}$ the normalized $\mfrak{l}_i$-invariant symmetric bilinear form on $\mfrak{l}_i$ and by $h_i^\vee$ the dual Coxeter number of $\mfrak{l}_i$ for $i=1,2,\dots,r$. Besides, we denote by $\smash{\kappa_0^{\mfrak{l}_0}}$ the restriction $\kappa_0{}_{|\mfrak{l}_0}$.
Let us recall that on $\mfrak{p}$ we have a $\mfrak{p}$-invariant symmetric bilinear form defined by
\begin{align}
  \kappa_c^\mfrak{p}(a,b) = - \tr_{\mfrak{g}/\mfrak{p}}(\ad(a)\ad(b)) \label{eq:kappa_c parabolic}
\end{align}
for $a, b \in \mfrak{p}$. We say that a $\mfrak{g}$-invariant symmetric bilinear form $\kappa$ on $\mfrak{g}$ is \emph{critical} if $(\kappa-\kappa_c^\mfrak{b})_{|\mfrak{h}}=0$ for some Borel subalgebra $\mfrak{b}=\mfrak{h} \oplus \mfrak{n}$ of $\mfrak{g}$. The critical $\mfrak{g}$-invariant symmetric bilinear form on $\mfrak{g}$ we will denote by $\kappa_c$. In fact, we have that $\kappa_c = -{1\over 2}\kappa_\mfrak{g}=-h^\vee \kappa_0$.
\medskip

\proposition{\label{prop:critical levels}
Let $\mfrak{p}$ be a parabolic subalgebra of $\mfrak{g}$. Then for $k \in \C$ we have
\begin{align*}
\begin{aligned}
  (k\kappa_0 - \kappa_c^\mfrak{p})_{|\mfrak{l}_i}  &= {(k + h^\vee)\ind_D(\mfrak{l}_i,\mfrak{u}) + kh^\vee_i \over h^\vee}\,\kappa_0^{\mfrak{l}_i}, \\
  (k\kappa_0 - \kappa_c^\mfrak{p})_{|\mfrak{l}_0} &= (k+h^\vee) \kappa_0^{\mfrak{l}_0}
\end{aligned}
\end{align*}
for $i=1,2,\dots,r$. In particular, we have $(k\kappa_0 - \kappa_c^\mfrak{p})_{|\mfrak{l}_i}=0$ if and only if
\begin{align*}
  k = - {\ind_D(\mfrak{l}_i,\mfrak{u})h^\vee \over \ind_D(\mfrak{l}_i,\mfrak{u}) + h^\vee_i}
\end{align*}
for $i=1,2,\dots,r$ and $(k\kappa_0-\kappa_c^\mfrak{p})_{|\mfrak{l}_0}=0$ if $k=-h^\vee$.}

\proof{For $a,b \in \mfrak{l}$, we have
\begin{align*}
  \kappa_\mfrak{g}(a,b)&= \tr_\mfrak{g}(\ad(a)\ad(b)) = \tr_{\widebar{\mfrak{u}}}(\ad(a)\ad(b)) + \tr_\mfrak{l}(\ad(a)\ad(b)) + \tr_\mfrak{u}(\ad(a)\ad(b)) \\
  &=2\tr_\mfrak{u}(\ad(a)\ad(b)) + \tr_\mfrak{l}(\ad(a)\ad(b))
\end{align*}
and
\begin{align*}
  \kappa_c^\mfrak{p}(a,b) = -\tr_{\mfrak{g}/\mfrak{p}}(\ad(a)\ad(b)) = -\tr_\mfrak{u}(\ad(a)\ad(b)).
\end{align*}
Hence, for $k \in \C$ we may write
\begin{align*}
  (k\kappa_0-\kappa_c^\mfrak{p})(a,b) &= {k \over 2h^\vee}\,\kappa_\mfrak{g}(a,b) - \kappa_c^\mfrak{p}(a,b) \\
  &= {k \over 2h^\vee}\, (2\tr_\mfrak{u}(\ad(a)\ad(b))+\tr_\mfrak{l}(\ad(a)\ad(b))) + \tr_\mfrak{u}(\ad(a)\ad(b)) \\
  & = {k \over 2h^\vee}\,(2\ind_D(\mfrak{l}_i,\mfrak{u})\kappa_0^{\mfrak{l}_i}(a,b) +\ind_D(\mfrak{l}_i,\mfrak{l})\kappa_0^{\mfrak{l}_i}(a,b)) + \ind_D(\mfrak{l}_i,\mfrak{u})\kappa_0^{\mfrak{l}_i}(a,b)\\
  & = \bigg({k \over h^\vee}\,(\ind_D(\mfrak{l}_i,\mfrak{u})+h^\vee_i)+ \ind_D(\mfrak{l}_i,\mfrak{u})\!\bigg)\kappa_0^{\mfrak{l}_i}(a,b)
\end{align*}
for $a,b \in \mfrak{l}_i$ and $i=1,2,\dots,r$. Similarly, for $k \in \C$ we get
\begin{align*}
  (k\kappa_0-\kappa_c^\mfrak{p})(a,b) = {k + h^\vee \over h^\vee} \tr_\mfrak{u}(\ad(a)\ad(b)) = {k + h^\vee \over 2h^\vee}\, \kappa_\mfrak{g}(a,b)
   = (k+h^\vee)\kappa_0(a,b)
\end{align*}
for $a,b \in \mfrak{l}_0$.}

\vspace{-2mm}


\subsection{Chiral differential operators}
\label{subsection:Chiral differential operators}

Let $G$ be a connected reductive algebraic group with its Lie algebra $\mfrak{g}$. We shall use the notation introduced in Section \ref{subsec:momentum}. Let $\kappa$ be a $\mfrak{g}$-invariant symmetric bilinear form on $\mfrak{g}$. Then there exists a sheaf of chiral differential operators $\mcal{D}_{G,\kappa}^{ch}$ on $G$ together with two injective homomorphisms
\begin{align}
  j_{\kappa,\mfrak{g}}^L \colon \mcal{V}^\kappa(\mfrak{g}) \rarr \Gamma(G,\mcal{D}_{G,\kappa}^{ch}) \qquad \text{and} \qquad j_{\kappa,\mfrak{g}}^R \colon \mcal{V}^{\kappa^*}\!(\mfrak{g}) \rarr \Gamma(G,\mcal{D}_{G,\kappa}^{ch})
\end{align}
of vertex algebras, where $\kappa^*$ is the \emph{dual} $\mfrak{g}$-invariant symmetric bilinear form on $\mfrak{g}$ defined by
\begin{align}
  \kappa + \kappa^* = - \kappa_\mfrak{g}, \label{eq:dual-level}
\end{align}
such that
\begin{align*}
\begin{aligned}
  \mcal{V}^\kappa(\mfrak{g}) &= \Com(\mcal{V}^{\kappa^*}\!(\mfrak{g}), \Gamma(G,\mcal{D}_{G,\kappa}^{ch})) = (\Gamma(G,\mcal{D}_{G,\kappa}^{ch}))^{j_{\kappa,\mfrak{g}}^R(\mfrak{g}[[t]])}, \\
  \mcal{V}^{\kappa^*}\!(\mfrak{g}) &= \Com(\mcal{V}^\kappa(\mfrak{g}), \Gamma(G,\mcal{D}_{G,\kappa}^{ch})) = (\Gamma(G,\mcal{D}_{G,\kappa}^{ch}))^{j_{\kappa,\mfrak{g}}^L(\mfrak{g}[[t]])}.
\end{aligned}
\end{align*}
Moreover, we have that
\begin{align*}
  \Gamma(G,\mcal{D}_{G,\kappa}^{ch}) \simeq U(\widehat{\mfrak{g}}_\kappa)\otimes_{U(\mfrak{g}[[t]] \oplus \C c)}\mcal{O}_{J_\infty(G)}(J_\infty(G))
\end{align*}
as $\mcal{V}^\kappa(\mfrak{g})$-modules, where $\mfrak{g}[[t]]$ acts naturally on $\mcal{O}_{J_\infty(G)}(J_\infty(G))$ and $c$ acts as the identity.
\medskip

Now let us recall that for a Lie algebra $\mfrak{a}$ and an $\mfrak{a}$-invarint symmetric bilinear form $\kappa$ on $\mfrak{a}$, the BRST cohomology $H^{\infty/2+\bullet}(\widehat{\mfrak{a}}_\kappa,M)$ with coefficient in a $\mcal{V}^\kappa(\mfrak{a})$-module $M$ is well-defined if and only if $\kappa=-\kappa_\mfrak{a}$, where $\kappa_\mfrak{a}$ is the Cartan--Killing form of $\mfrak{a}$.
\medskip

\lemma{\label{Lemma:killing form} We have that $(\kappa-\kappa_c^\mfrak{p})_{|[\mfrak{l},\mfrak{l}]}=0$ if and only if
$\kappa^*{}_{|[\mfrak{p},\mfrak{p}]}=-\kappa_{[\mfrak{p},\mfrak{p}]}$.}

\proof{As we have
\begin{align*}
  \kappa_\mfrak{g}(a,b) =\tr_\mfrak{g}(\ad(a)\ad(b)) = \tr_\mfrak{p}(\ad(a)\ad(b)) + \tr_{\mfrak{g}/\mfrak{p}}(\ad(a)\ad(b))= \kappa_\mfrak{p}(a,b) -\kappa_c^\mfrak{p}(a,b)
\end{align*}
for $a,b \in \mfrak{p}$, by using \eqref{eq:dual-level} we may write
\begin{align*}
  \kappa^*{}_{|[\mfrak{p},\mfrak{p}]} = -(\kappa+\kappa_\mfrak{g})_{|[\mfrak{p},\mfrak{p}]} = - \kappa_\mfrak{p}{}_{|[\mfrak{p},\mfrak{p}]}- (\kappa - \kappa_c^\mfrak{p})_{|[\mfrak{p},\mfrak{p}]} = -\kappa_{[\mfrak{p},\mfrak{p}]}- (\kappa - \kappa_c^\mfrak{p})_{|[\mfrak{p},\mfrak{p}]},
\end{align*}
which gives us that $\kappa^*{}_{|[\mfrak{p},\mfrak{p}]} = -\kappa_{[\mfrak{p},\mfrak{p}]}$ if and only if $(\kappa - \kappa_c^\mfrak{p})_{|[\mfrak{p},\mfrak{p}]}=0$. Further, it follows immediately that $(\kappa - \kappa_c^\mfrak{p})_{|[\mfrak{p},\mfrak{p}]} = 0$ is equivalent to the condition $(\kappa - \kappa_c^\mfrak{p})_{|[\mfrak{l},\mfrak{l}]}=0$, which implies the required assertion.}

Let us assume that $(\kappa-\kappa_c^\mfrak{p})_{|[\mfrak{l},\mfrak{l}]}=0$.
Then by using Lemma \ref{Lemma:killing form} and following \cite{Gorbounov-Malikov-Schechtman2001} we can define a sheaf of chiral differential operators $\smash{\mcal{D}_{Y,\kappa}^{ch}}$ on the generalized base affine space $Y=G/(P,P)$ as the sheaf associated to a presheaf
\begin{align*}
  V \mapsto H^{{\infty \over 2} + 0}(\widehat{[\mfrak{p},\mfrak{p}]}_{\kappa^*}, \mcal{D}_{G,\kappa}^{ch}(q^{-1}(V)))
\end{align*}
for an open subset $V$ of $Y$, where $q \colon G \rarr G/(P,P)$ is the canonical $G$-equivariant projection. Of the two mappings $\smash{j_{\kappa,\mfrak{g}}^L}$ and $\smash{j_{\kappa,\mfrak{g}}^R}$, the former survives intact because for the BRST cohomology we use the action $\smash{j_{\kappa,\mfrak{g}}^R}$, the latter is maimed. Hence, we get homomorphisms
\begin{align}
  \pi_{\kappa,\mfrak{g}}^\mfrak{p} \colon \mcal{V}^\kappa(\mfrak{g}) \rarr \Gamma(Y,\mcal{D}_{Y,\kappa}^{ch}) \qquad \text{and} \qquad \sigma_{\kappa,\mfrak{g}}^\mfrak{p} \colon \mcal{V}^{\kappa^*- \kappa_c^\mfrak{p}}(\mfrak{z}(\mfrak{l})) \rarr \Gamma(Y,\mcal{D}_{Y,\kappa}^{ch}) \label{eq:chiral moment map X}
\end{align}
of vertex algebras. Let us note that $\mfrak{z}(\mfrak{l})$ is the Lie algebra of the torus $A=P/(P,P)$.

Let us define the sheaf $\smash{\widetilde{\mcal{D}}_{X,\kappa}^{ch}}$ of vertex algebras on $X$
by
\begin{align*}
  \widetilde{\mcal{D}}_{X,\kappa}^{ch} = (\pi_*\mcal{D}_{Y,\kappa}^{ch})^{\sigma_{\kappa,\mfrak{g}}^\mfrak{p}(\mfrak{z}(\mfrak{l})[[t]])},
\end{align*}
where $\pi \colon Y \rarr X$ is the canonical $G$-equivariant projection. The principal $A$-fibration $\pi \colon Y \rarr X$ induces the principal $J_\infty(A)$-fibration $J_\infty(\pi) \colon J_\infty(Y) \rarr J_\infty(X)$. The mapping $\pi_{\kappa,\mfrak{g}}^\mfrak{p}$  descends to the vertex algebra homomorphism
\begin{align}
   \widetilde{\pi}_{\kappa,\mfrak{g}}^\mfrak{p} \colon \mcal{V}^\kappa(\mfrak{g}) \rarr \Gamma(X,\widetilde{\mcal{D}}_{X,\kappa}^{ch}). \label{eq:chiral moment map Y}
\end{align}
The sheaf of chiral differential operators $\mcal{D}_{Y,\kappa}^{ch}$ on $Y$ is a quantization of $(\pi_{Y})_* \mcal{O}_{J_\infty(\smash{T^*Y})}$, where $\pi_Y \colon J_\infty(T^*Y) \rarr T^*Y \rarr Y$ is the composition of the natural projections. Moreover, the mapping $\smash{\pi_{\kappa,\mfrak{g}}^\mfrak{p}}$ in \eqref{eq:chiral moment map X} is a quantization of the morphism $J_\infty(j \circ \mu_Y) \colon J_\infty(T^*Y) \rarr J_\infty(\mfrak{g}^*)$, where the morphism $j \colon \mfrak{g}^* \rarr \mfrak{g}^*$ is defined through $j(\alpha) = -\alpha$ for $\alpha \in \mfrak{g}^*$. Similarly, the sheaf $\smash{\widetilde{\mcal{D}}_{X,\kappa}^{ch}}$ is a quantization of $(\pi_X)_*\mcal{O}_{J_\infty(\smash{T^*Y/A})}$, where $\pi_X \colon J_\infty(T^*Y/A) \rarr T^*Y/A \rarr X$ is the composition of the natural projections. The  mapping $\widetilde{\pi}_{\kappa,\mfrak{g}}^\mfrak{p}$ in \eqref{eq:chiral moment map Y} is a quantization of the morphism
\begin{align}\label{eq:our-chiai-moment}
   J_\infty(j \circ \widebar{\mu}_Y) \colon J_\infty(T^*Y/A) \simeq J_\infty(G) \times_{J_\infty(P)}  J_\infty([\mfrak{p},\mfrak{p}]^\perp) \rarr J_\infty(\mfrak{g}^*)
\end{align}
induced by the generalized Grothendieck's simultaneous resolution
\begin{align*}
  \widebar{\mu}_Y \colon T^*Y/A \simeq G\times_P  [\mfrak{p},\mfrak{p}]^{\perp} \rarr \mfrak{g}^*.
\end{align*}
Let $U_e=\widebar{U}P \subset X$ be the big cell. By using the group isomorphism $P \simeq A \times (P,P)$ and the fact that $(\pi \circ q)^{-1}(U_e) = p^{-1}(U_e) = \widebar{U}P \subset G$, we obtain
\begin{align*}
   \mcal{D}^{ch}_{G,\kappa}(q^{-1}(\pi^{-1}(U_e))) \simeq  \mcal{D}^{ch}_{\widebar{U}}(\widebar{U}) \otimes_\C \mcal{D}_{A, \kappa-\kappa^{\mfrak{p}}_{c}}^{ch}(A) \otimes_\C
   \mcal{D}_{(P,P), \kappa-\kappa^{\mfrak{p}}_{c}}^{ch}((P,P)).
\end{align*}
Therefore, it follows immediately that
\begin{align*}
   \mcal{D}^{ch}_{Y,\kappa}(\pi^{-1}(U_e)) &\simeq \mcal{D}^{ch}_{\widebar{U}}(\widebar{U}) \otimes_\C \mcal{D}_{A, \kappa-\kappa^{\mfrak{p}}_{c}}^{ch}(A) \otimes_\C
   H^{{\infty \over 2} + 0}(\widehat{[\mfrak{p},\mfrak{p}]}_{\kappa^*}, \mcal{D}_{(P,P), \kappa-\kappa^{\mfrak{p}}_{c}}^{ch}((P,P))) \\
   & \simeq \mcal{D}^{ch}_{\widebar{U}}(\widebar{U}) \otimes_\C \mcal{D}_{A,\kappa-\kappa^{\mfrak{p}}_{c}}^{ch}(A).
\end{align*}
Moreover, we have $\smash{(\mcal{D}_{A,\kappa-\kappa^{\mfrak{p}}_{c}}^{ch} (A))^{\sigma_{\kappa,\mfrak{g}}^\mfrak{p}(\mfrak{z}(\mfrak{l})[[t]])}} \simeq
\mcal{V}^{\kappa-\kappa^{\mfrak{p}}_c}(\mfrak{z}(\mfrak{l}))$ and $\smash{\mcal{D}^{ch}_{\widebar{U}}(\widebar{U})} \simeq \mcal{M}_{\widebar{\mfrak{u}}}$, where $\mcal{M}_{\widebar{\mfrak{u}}}$ is the Weyl vertex algebra associated to $\widebar{\mfrak{u}}$, which gives us
\begin{align*}
   \widetilde{\mcal{D}}_{X,\kappa}^{ch}(U_e) \simeq \mcal{M}_{\widebar{\mfrak{u}}} \otimes_\C \mcal{V}^{\kappa-\kappa_c^\mfrak{p}}(\mfrak{z}(\mfrak{l})).
\end{align*}
The vertex algebra homomorphism
\begin{align}
  \widetilde{w}_{\kappa,\mfrak{g}}^\mfrak{p} \colon \mcal{V}^{\kappa}(\mfrak{g}) \rarr \Gamma(X,\widetilde{\mcal{D}}_{X,\kappa}^{ch}) \rarr \widetilde{\mcal{D}}_{X,\kappa}^{\ch}(U_e) \rarr \mcal{M}_{\widebar{\mfrak{u}}} \otimes_\C \mcal{V}^{\kappa-\kappa_c^\mfrak{p}}(\mfrak{z}(\mfrak{l})) \label{eq:gen.Wakimoto}
\end{align}
is a quantization of the restriction of the morphism \eqref{eq:our-chiai-moment} to the open set $J_\infty(\smash{\widetilde{U}_e}) \subset J_\infty(T^*Y/A)$, where $\smash{\widetilde{U}_e}$ is the preimage of $U_e$ by the projection $T^*Y/A \rarr X$. The vertex algebra homomorphism $\widetilde{w}_{\kappa,\mfrak{g}}^\mfrak{p}$ also induces the homomorphism
\begin{align*}
    \widetilde{\pi}_\mfrak{g}^\mfrak{p} \colon U(\mfrak{g}) \rarr \Gamma(X,(\pi_*\mcal{D}_Y)^A) \rarr (\pi_* \mcal{D}_{Y})^{A}(U_e) \rarr \eus{A}_{\widebar{\mfrak{u}}} \otimes_\C U(\mfrak{z}(\mfrak{l}))
\end{align*}
between their Zhu's algebras, where $\eus{A}_{\widebar{\mfrak{u}}}$ is the Weyl algebra associated to $\widebar{\mfrak{u}}$.
\medskip

Let $\{f_\alpha;\, \alpha \in \Delta_+^\mfrak{u}\}$ be a root basis of the nilpotent Lie algebra $\widebar{\mfrak{u}}$. We denote by $\{x_\alpha;\, \alpha \in \Delta_+^\mfrak{u}\}$ the linear coordinate functions on $\widebar{\mfrak{u}}$ with respect to the given basis of $\widebar{\mfrak{u}}$. The Weyl algebra $\eus{A}_{\widebar{\mfrak{u}}}$ of the vector space $\widebar{\mfrak{u}}$ is generated by $\{x_\alpha,\partial_{x_\alpha};\, \alpha \in \Delta_+^\mfrak{u}\}$ together with the canonical commutation relations. In addition, generating fields of the Weyl vertex algebra $\mcal{M}_{\widebar{\mfrak{u}}}$ are given by
\begin{align*}
  a_\alpha(z) = \sum_{n \in \Z} \partial_{x_{\alpha,n}} z^{-n-1} \qquad \text{and} \qquad a_\alpha^*(z) = \sum_{n \in \Z} x_{\alpha,-n} z^{-n}
\end{align*}
for $\alpha \in \Delta_+^\mfrak{u}$ and satisfy the following commutation relations
\begin{align*}
  [a_\alpha(z), a_\beta(w)] = 0, \qquad [a_\alpha(z), a_\beta^*(w)] = \delta_{\alpha,\beta} \delta(z-w), \qquad [a_\alpha^*(z), a_\beta^*(w)] = 0
\end{align*}
for $\alpha,\beta \in \Delta_+^\mfrak{u}$. Besides, for an element $h \in \mfrak{h}$, we denote by $h_{\mfrak{z}(\mfrak{l})}$ the $\mfrak{z}(\mfrak{l})$-part of $h$ with respect to the direct sum decomposition $\mfrak{l} = \mfrak{z}(\mfrak{l}) \oplus [\mfrak{l},\mfrak{l}]$.

We have
\begin{align*}
    & \widetilde{\pi}_\mfrak{g}^\mfrak{p}(f_\gamma) = -\partial_{x_\gamma} \otimes 1 +\sum_{\alpha\in \Delta^\mfrak{u}_+,\, \alpha>\gamma} q_\alpha^\gamma(x_\beta) \partial_{x_\alpha} \otimes 1 \qquad  (\gamma \in \Delta^\mfrak{u}_+), \\
    & \widetilde{\pi}_\mfrak{g}^\mfrak{p}(f_\gamma) = \sum_{\alpha \in \Delta^\mfrak{u}_+} q_\alpha^\gamma(x_\beta) \partial_{x_\alpha} \otimes 1 \qquad (\gamma \in  \Delta^\mfrak{l}_+),\\
    & \widetilde{\pi}_\mfrak{g}^\mfrak{p}(e_\gamma) = \sum_{\alpha \in \Delta^\mfrak{u}_+} p_\alpha^\gamma(x_\beta) \partial_{x_\alpha} \otimes 1 \qquad (\gamma \in  \Delta^\mfrak{l}_+),\\
    & \widetilde{\pi}_\mfrak{g}^\mfrak{p}(h) = \sum_{\alpha\in \Delta^\mfrak{u}_+} c_\alpha^h\, x_\alpha \partial_{x_\alpha} \otimes 1 + 1 \otimes h_{\mfrak{z}(\mfrak{l})} \qquad (h\in \mfrak{h}),
\end{align*}
where $p_\alpha^\gamma$ and $q_\alpha^\gamma$ are polynomials on the vector space $\widebar{\mfrak{u}}$ of weights $\gamma+\alpha$ and $-\gamma+\alpha$, respectively, and $c^h_\alpha \in \C$. Here we have set the weights of $x_\alpha$ and $\partial_{x_\alpha}$ equal to $\alpha$ and $-\alpha$, respectively.
\medskip

\proposition{\label{prop:cdo formula}Let $\kappa$ be a $\mfrak{g}$-invariant symmetric bilinear form on $\mfrak{g}$ such that $(\kappa-\kappa_c^\mfrak{p})_{|[\mfrak{l},\mfrak{l}]}=0$. Then we have
\begin{align*}
    & \widetilde{w}_{\kappa,\mfrak{g}}^\mfrak{p}(f_\gamma(z))=
    -a_\gamma(z) \otimes 1 +\sum_{\alpha\in \Delta^\mfrak{u}_+,\, \alpha>\gamma} \normOrd{q_\alpha^\gamma(a_\beta^*(z))a_\alpha(z)} \otimes 1
    \qquad (\gamma\in \Delta^{\mfrak{u}}_+),\\
    & \widetilde{w}_{\kappa,\mfrak{g}}^\mfrak{p}(f_\gamma(z))= \sum_{\alpha \in \Delta^\mfrak{u}_+} \normOrd{q_\alpha^\gamma(a_\beta^*(z))a_\alpha(z)} \otimes 1 \qquad  (\gamma \in  \Delta^{\mfrak{l}}_+),\\
    & \widetilde{w}_{\kappa,\mfrak{g}}^\mfrak{p}(e_\gamma(z))= \sum_{\alpha \in \Delta^\mfrak{u}_+} \normOrd{p_\alpha^\gamma(a_\beta^*(z))a_\alpha(z)} \otimes 1 \qquad  (\gamma \in  \Delta^{\mfrak{l}}_+),\\
    & \widetilde{w}_{\kappa,\mfrak{g}}^\mfrak{p}(h(z))= \sum_{\alpha\in \Delta^\mfrak{u}_+} c_\alpha^h\, \normOrd{a_\alpha^*(z)a_\alpha(z)} \otimes 1 + 1\otimes h_{\mfrak{z}(\mfrak{l})}(z) \qquad (h\in \mfrak{h}).
\end{align*}
Moreover, the homomorphism $\smash{\widetilde{w}_{\kappa,\mfrak{g}}^\mfrak{p}} \colon \mcal{V}^\kappa(\mfrak{g}) \rarr \mcal{M}_{\widebar{\mfrak{u}}} \otimes_\C \mcal{V}^{\kappa-\kappa_c^\mfrak{p}}(\mfrak{z}(\mfrak{l}))$ of $\N_0$-graded vertex algebras is the unique extension of the vertex algebra homomorphism $\mcal{V}^\kappa(\widebar{\mfrak{p}}) \rarr \mcal{M}_{\widebar{\mfrak{u}}} \otimes_\C \mcal{V}^{\kappa-\kappa_c^\mfrak{p}} (\mfrak{z}(\mfrak{l}))$ defined by the above formulas, where $\widebar{\mfrak{p}} = \mfrak{l} \oplus \widebar{\mfrak{u}}$ is the opposite parabolic subalgebra of $\mfrak{g}$.}

\proof{The first statement follows from the weight consideration. Indeed, the conformal weights of $\widetilde{w}_{\kappa,\mfrak{g}}^\mfrak{p}(a(z))$ for $a\in \mfrak{g}$, $a_\alpha(z)$, $a^*_\alpha(z)$ for $\alpha \in \Delta_+^\mfrak{u}$ and $h(z)$ for $h \in \mfrak{z}(\mfrak{l})$ are $1$, $1$, $0$ and $1$, respectively. Next, let us consider the $\Z$-grading
\begin{align}
   \mfrak{g} = \bigoplus_{j\in \Z} \mfrak{g}_j \label{eq:defining-grading-of-p}
\end{align}
of $\mfrak{g}$ such that $\mfrak{u}=\smash{\bigoplus_{j>0}}\,\mfrak{g}_j$, $\widebar{\mfrak{u}}= \smash{\bigoplus_{j<0}}\,\mfrak{g}_j$ and $\mfrak{l}=\mfrak{g}_0$. There exists the grading element $\xi \in \mfrak{z}(\mfrak{l})$ such that $\mfrak{g}_j=\{a\in \mfrak{g};\, [\xi,a]=ja\}$. We extend this grading to $\smash{\widetilde{\mcal{D}}_{X,\kappa}^{ch}(U_e)}$ by setting the degree of $h(z)$ to be $0$ for $h \in \mfrak{z}(\mfrak{l})$ and the degree of $a_\alpha(z)$,  $a_\alpha^*(z)$ to be $-\alpha(\xi)$, $\alpha(\xi)$, respectively, for $\alpha \in \Delta_+^\mfrak{u}$. Then $\smash{\widetilde{w}_{\kappa,\mfrak{g}}^\mfrak{p}(a(z))}$ has a non-positive degree for all $a \in \widebar{\mfrak{p}}$. Now let us consider the difference of $\smash{\widetilde{w}_{\kappa,\mfrak{g}}^\mfrak{p}(a(z))}$ and the right-hand side of the formula in the preposition, which vanishes in the Zhu's $C_2$-algebra of $\smash{\widetilde{\mcal{D}}_{X,\kappa}^{ch}(U_e)}$. It is a sum of monomials
of the form $\normOrd{a^*_{\beta_1}(z)\dots a^*_{\beta_r}(z) \partial_z a^*_{\gamma}(z)}$
since it has the conformal weight one. But such an element has a positive degree and thus must be zero.

Let us prove the second statement. Let $\smash{\widehat{w}_{\kappa,\mfrak{g}}^\mfrak{p}} \colon \mcal{V}^\kappa(\mfrak{p}) \rarr \mcal{M}_{\widebar{\mfrak{u}}} \otimes_\C \mcal{V}^{\kappa-\kappa_c^\mfrak{p}}(\mfrak{z}(\mfrak{l}))$ be another extension. We show by induction on $j \geq 0$ that $\smash{\widehat{w}_{\kappa,\mfrak{g}}^\mfrak{p}(a(z))} = \smash{\widetilde{w}_{\kappa,\mfrak{g}}^\mfrak{p}(a(z))}$ for $a \in \mfrak{g}_j$. By the induction hypothesis, we have that $\smash{\widehat{w}_{\kappa,\mfrak{g}}^\mfrak{p}(a(z))} - \smash{\widetilde{w}_{\kappa,\mfrak{g}}^\mfrak{p}(a(z))}$ commutes with $\smash{\widetilde{w}_{\kappa,\mfrak{g}}^\mfrak{p}(x(z))}$ for all $x\in \widebar{\mfrak{u}}$. Hence, it belongs to $\smash{(\mcal{M}_{\widebar{\mfrak{u}}})^{\widebar{\mfrak{u}}[[t]]}} \otimes_\C \mcal{V}^{\kappa-\kappa_c^\mfrak{p}}(\mfrak{z}(\mfrak{l}))$. However, we have that $\smash{(\mcal{M}_{\widebar{\mfrak{u}}})^{\widebar{\mfrak{u}}[[t]]}}$ is spanned by the monomials $\smash{j^R_{\kappa,\widebar{\mfrak{u}}}(f_{\alpha_1,-n_1}) j^R_{\kappa,\widebar{\mfrak{u}}}(f_{\alpha_2,-n_2})\dots j^R_{\kappa,\widebar{\mfrak{u}}}(f_{\alpha_r,-n_r})}$, where $j^R_{\kappa,\widebar{\mfrak{u}}}$ is the right action of $\mcal{V}^{\kappa^*}\!(\widebar{\mfrak{u}})$ on $\smash{\mcal{D}^{ch}_{\widebar{U}}(\widebar{U})} \simeq \mcal{M}_{\widebar{\mfrak{u}}}$. It follows that $\smash{(\mcal{M}_{\widebar{\mfrak{u}}})^{\widebar{\mfrak{u}}[[t]]} \otimes_\C \mcal{V}^{\kappa-\kappa_c^\mfrak{p}}(\mfrak{z}(\mfrak{l}))}$ is spanned by the elements of non-positive degrees with respect to the grading \eqref{eq:defining-grading-of-p}. Therefore $\smash{\widehat{w}_{\kappa,\mfrak{g}}^\mfrak{p}(a(z))} - \smash{\widetilde{w}_{\kappa,\mfrak{g}}^\mfrak{p}(a(z))}$ must be zero.}

\remark{\label{remark:Wakimoto} In the case where $P$ is the Borel subgroup the condition $(\kappa-\kappa_c^\mfrak{p})_{|[\mfrak{l},\mfrak{l}]}=0$ is empty and it follows from Proposition \ref{prop:cdo formula} that \eqref{eq:gen.Wakimoto} is the Wakimoto free field realization (\cite{Feigin-Frenkel1990,Frenkel2005}) of $\mcal{V}^\kappa(\mfrak{g})$.}


\subsection{Feigin--Frenkel homomorphism for parabolic subalgebras}
\label{subsec:FF homomorphism}

Let $\mfrak{g}$ be a reductive Lie algebra and let $\mfrak{b}=\mfrak{h} \oplus \mfrak{n}$ be a Borel subalgebra of $\mfrak{g}$ with the Cartan subalgebra $\mfrak{h}$, the nilradical $\mfrak{n}$ and the opposite nilradical $\widebar{\mfrak{n}}$. Let us consider a parabolic subalgebra $\mfrak{p}$ of $\mfrak{g}$. We will assume that $\mfrak{p}$ contains the opposite Borel subalgebra $\mfrak{h} \oplus \widebar{\mfrak{n}}$ of $\mfrak{g}$. Let
\begin{align*}
  \mfrak{p} = \mfrak{l} \oplus \widebar{\mfrak{u}}
\end{align*}
be the Levi decomposition of $\mfrak{p}$ with the Levi subalgebra $\mfrak{l}$, the nilradical $\widebar{\mfrak{u}}$ and the opposite nilradical $\mfrak{u}$. Let us also note that $\mfrak{p}$ is the opposite parabolic subalgebra to the standard parabolic subalgebra $\mfrak{l} \oplus \mfrak{u}$ of $\mfrak{g}$, it will be called an opposite standard parabolic subalgebra. By a standard parabolic subalgebra of $\mfrak{g}$ we mean a parabolic subalgebra of $\mfrak{g}$, which contains the Borel subalgebra $\mfrak{b}$.
\medskip

Let $\{e_\alpha;\, \alpha \in \Delta_+^\mfrak{u}\}$ be a root basis of the nilpotent Lie subalgebra $\mfrak{u}$. We denote by $\{x_\alpha;\, \alpha \in \Delta_+^\mfrak{u}\}$ the linear coordinate functions on $\mfrak{u}$ with respect to the given basis of $\mfrak{u}$. The Weyl algebra $\eus{A}_\mfrak{u}$ of the vector space $\mfrak{u}$ is generated by $\{x_\alpha,\partial_{x_\alpha};\, \alpha \in \Delta_+^\mfrak{u}\}$ together with the canonical commutation relations. Further, let us consider the Weyl vertex algebra $\mcal{M}_\mfrak{u}$. Then generating fields of $\mcal{M}_\mfrak{u}$ are given by
\begin{align*}
  a_\alpha(z) = \sum_{n \in \Z} \partial_{x_{\alpha,n}} z^{-n-1} \qquad \text{and} \qquad a_\alpha^*(z) = \sum_{n \in \Z} x_{\alpha,-n} z^{-n}
\end{align*}
for $\alpha \in \Delta_+^\mfrak{u}$ and satisfy the following commutation relations
\begin{align*}
  [a_\alpha(z), a_\beta(w)] = 0, \qquad [a_\alpha(z), a_\beta^*(w)] = \delta_{\alpha,\beta} \delta(z-w), \qquad [a_\alpha^*(z), a_\beta^*(w)] = 0
\end{align*}
for $\alpha,\beta \in \Delta_+^\mfrak{u}$. Let us denote by $\eus{F}^{\mfrak{g},\mfrak{p}}_{{\rm loc}}(z)$ the vector space of all differential polynomials in $a^*_\alpha(z)$ for $\alpha \in \Delta_+^\mfrak{u}$, and by $\eus{C}^\mfrak{g}_{{\rm loc}}(z)$ the vector space of all fields of the form $a(z)$ for $a \in \mfrak{g}$. As the linear mapping $a \mapsto a(z)$ from $\mfrak{g}$ to $\eus{C}^\mfrak{g}_{{\rm loc}}(z)$ is an isomorphism of vector spaces, we endow $\eus{C}^\mfrak{g}_{{\rm loc}}(z)$ with a Lie algebra structure via this isomorphism. Hence, we obtain that the vector spaces $\mfrak{g} \otimes_\C \eus{F}^{\mfrak{g},\mfrak{p}}_{{\rm loc}}(z)$ and $\eus{C}^\mfrak{g}_{{\rm loc}}(z) \otimes_\C \eus{F}^{\mfrak{g},\mfrak{p}}_{{\rm loc}}(z)$ have the natural structure of Lie algebras, and moreover $\mfrak{g} \otimes_\C \eus{F}^{\mfrak{g},\mfrak{p}}_{{\rm loc}}(z)$ and $\eus{C}^\mfrak{g}_{{\rm loc}}(z) \otimes_\C \eus{F}^{\mfrak{g},\mfrak{p}}_{{\rm loc}}(z)$ are $\mfrak{g} \otimes_\C \eus{F}^{\mfrak{g},\mfrak{p}}_{{\rm loc}}(z)$-modules through the adjoint action.
\medskip

Let $\kappa$ be a $\mfrak{g}$-invariant symmetric bilinear form on $\mfrak{g}$. Then by \cite[Theorem 6.3.1]{Frenkel2007-book} there is a homomorphism
\begin{align*}
  w_{\kappa,\mfrak{g}}^\mfrak{p} \colon \mcal{V}^\kappa(\mfrak{g}) \rarr \mcal{M}_\mfrak{u} \otimes_\C \mcal{V}^{\kappa-\kappa_c^\mfrak{p}}(\mfrak{p})
\end{align*}
of $\N_0$-graded vertex algebras, where the $\mfrak{p}$-invariant symmetric bilinear form $\kappa_c^\mfrak{p}$ on $\mfrak{p}$ is defined by the formula \eqref{eq:kappa_c parabolic}.
For our purposes, we need also an explicit form of this homomorphism which is given in the following theorem.
\medskip

\theorem{\label{thm:FF homomorphism 1-graded}
Let $\kappa$ be a $\mfrak{g}$-invariant symmetric bilinear form on $\mfrak{g}$. Let us assume that $\mfrak{u}$ is a commutative Lie algebra. Then there exists a homomorphism
\begin{align*}
  w_{\kappa,\mfrak{g}}^\mfrak{p} \colon \mcal{V}^\kappa(\mfrak{g}) \rarr \mcal{M}_\mfrak{u} \otimes_\C \mcal{V}^{\kappa - \kappa_c^\mfrak{p}}(\mfrak{p})
\end{align*}
of $\N_0$-graded vertex algebras such that
\begin{align*}
  w_{\kappa,\mfrak{g}}^\mfrak{p}(a(z)) = \begin{cases}
    - {\displaystyle \sum_{\alpha \in \Delta_+^\mfrak{u}}} [a]_\alpha a_\alpha(z) & \text{if $a \in \mfrak{u}$}, \\[5mm]
    {\displaystyle \sum_{\alpha \in \Delta_+^\mfrak{u}}}\, \normOrd{\,[\ad(u(z))(a)]_\alpha a_\alpha(z)} + a(z) & \text{if $a \in \mfrak{l}$}, \\[5mm]
    - {1 \over 2} {\displaystyle \sum_{\alpha \in \Delta_+^\mfrak{u}}} \normOrd{\,[(\ad(u(z)))^2(a)]_\alpha a_\alpha(z)} - \ad(u(z))(a(z)) - \kappa(\partial_zu(z),a) + a(z) & \text{if $a \in \widebar{\mfrak{u}}$},
  \end{cases}
\end{align*}
for $a \in \mfrak{g}$, where $[a]_\alpha$ denotes the $\alpha$-th coordinate of $a \in \mfrak{u}$ with respect to the basis $\{e_\alpha;\, \alpha \in \Delta_+^\mfrak{u}\}$ of $\mfrak{u}$ and the element $u(z) \in \mfrak{u} \otimes_\C \smash{\eus{F}^{\mfrak{g},\mfrak{p}}_{\rm loc}(z)}$ is given by
\begin{align}
  u(z)=\sum_{\alpha \in \Delta_+^\mfrak{u}}  a_\alpha^*(z) e_\alpha. \label{eq:u(z) def}
\end{align}
The $\mfrak{p}$-invariant symmetric bilinear form $\kappa_c^\mfrak{p}$ on $\mfrak{p}$ is given through $\kappa_c^\mfrak{p}(a,b) = - \tr_\mfrak{u}(\ad(a)\ad(b))$ for $a,b \in \mfrak{l}$ and $\kappa_c^\mfrak{p}(a,b)=0$ provided $a \in \widebar{\mfrak{u}}$ or $b \in \widebar{\mfrak{u}}$.}

\proof{To simplify the notation we shall write $\Phi(a(z))$ instead of $w_{\kappa,\mfrak{g}}^\mfrak{p}(a(z))$ for $a\in \mfrak{g}$. We need to show that
\begin{align*}
  [\Phi(a(z)),\Phi(b(w))] = \Phi([a,b](w))\delta(z-w) + \kappa(a,b)\partial_w \delta(z-w)
\end{align*}
for $a, b \in \mfrak{g}$. The proof is a straightforward but tedious calculation based on the Wick theorem. We will use the following identities
\begin{gather*}
  [a_\alpha(z), [\ad(u(w))(a)]_\beta] = [e_\alpha,a]_\beta \delta(z-w), \\[3mm]
  [a_\alpha(z), [(\ad(u(w)))^2(a)]_\beta] = 2 [u(w),[e_\alpha,a]]_\beta \delta(z-w)
\end{gather*}
and
\begin{align*}
  [a_\alpha(z), \ad(u(w))(a(w))] = [e_\alpha,a](w) \delta(z-w), \qquad [a_\alpha(z), \kappa(\partial_wu(w),a)] = \kappa(e_\alpha,a)\partial_w \delta(z-w)
\end{align*}
for $a \in \mfrak{g}$ and $\alpha,\beta \in \Delta_+^\mfrak{u}$. For $a \in \mfrak{u}$, we may write
\begin{align*}
  [\Phi(a(z)),\Phi(b(w))] = 0 = \Phi([a,b](w))\delta(z-w)
\end{align*}
if $b \in \mfrak{u}$,
\begin{align*}
  [\Phi(a(z)),\Phi(b(w))] &= - \sum_{\alpha,\beta \in \Delta_+^\mfrak{u}} [a]_\alpha [e_\alpha,b]_\beta a_\beta(w) \delta(z-w) = - \sum_{\beta \in \Delta_+^\mfrak{u}} [a,b]_\beta a_\beta(w) \delta(z-w) \\
  & = \Phi([a,b](w)) \delta(z-w)
\end{align*}
if $b \in \mfrak{l}$, and
\begin{align*}
  [\Phi(a(z)),\Phi(b(w))] &= \sum_{\alpha,\beta \in \Delta_+^\mfrak{u}} [a]_\alpha \normOrd{\, [u(w),[e_\alpha,b]]_\beta a_\beta(w)}\delta(z-w) \\
  &\quad + \sum_{\alpha \in \Delta_+^\mfrak{u}} [a]_\alpha [e_\alpha,b](w)\delta(z-w) +  \sum_{\alpha  \in \Delta_+^\mfrak{u}} [a]_\alpha \kappa(e_\alpha,b) \partial_w \delta(z-w) \\
  &=  \sum_{\beta \in \Delta_+^\mfrak{u}} \!\normOrd{\,[u(w),[a,b]]_\beta a_\beta(w)}\delta(z-w) + [a,b](w)\delta(z-w) \\
  &\quad + \kappa(a,b) \partial_w \delta(z-w) = \Phi([a,b](w))\delta(z-w) + \kappa(a,b) \partial_w \delta(z-w)
\end{align*}
if $b \in \widebar{\mfrak{u}}$. Further, for $a \in \mfrak{l}$ we have
\begin{align*}
  [\Phi(a(z)),\Phi(b(w))] &= \sum_{\alpha,\beta \in \Delta_+^\mfrak{u}} [e_\alpha,b]_\beta \normOrd{\,[u(w), a]_\alpha a_\beta(w)}\delta(z-w) \\
  &\quad -\sum_{\alpha,\beta \in \Delta_+^\mfrak{u}} [e_\alpha,a]_\beta \normOrd{\,[u(w), b]_\alpha a_\beta(w)}\delta(z-w) -\sum_{\alpha,\beta \in \Delta_+^\mfrak{u}} [e_\alpha,b]_\beta [e_\beta,a]_\alpha \partial_w \delta(z-w) \\
  &\quad + [a,b](w)\delta(z-w) + (\kappa-\kappa_c^\mfrak{p})(a,b) \partial_w \delta(z-w) \\
  &= \sum_{\beta \in \Delta_+^\mfrak{u}} \normOrd{\,([[u(w),a],b]_\beta + [a,[u(w), b]]_\beta) a_\beta(w)}\delta(z-w) + [a,b](w)\delta(z-w) \\
  &\quad  + \kappa(a,b) \partial_w \delta(z-w) \\
  &= \sum_{\beta \in \Delta_+^\mfrak{u}} \normOrd{\,[u(w), [a,b]]_\beta a_\beta(w)} \delta(z-w) + [a,b](w)\delta(z-w) + \kappa(a,b) \partial_w \delta(z-w) \\
  &= \Phi([a,b](w))\delta(z-w) + \kappa(a,b)\partial_w \delta(z-w)
\end{align*}
if $b \in \mfrak{l}$, where we used the identity
\begin{align*}
  \sum_{\alpha,\beta \in \Delta_+^\mfrak{u}} [e_\alpha,b]_\beta [e_\beta,a]_\alpha = \sum_{\alpha \in \Delta_+^\mfrak{u}} [[e_\alpha,b],a]_\alpha = \tr_\mfrak{u}(\ad(a)\ad(b)) = -\kappa_c^\mfrak{p}(a,b),
\end{align*}
and
\begin{align*}
  [\Phi(a(z)),\Phi(b(w))] &= -\sum_{\alpha,\beta \in \Delta_+^\mfrak{u}} \normOrd{\, [u(w),[e_\alpha,b]]_\beta [u(w),a]_\alpha a_\beta(w)} \delta(z-w) \\
  & \quad +{1\over 2} \sum_{\alpha,\beta \in \Delta_+^\mfrak{u}} [e_\beta,a]_\alpha \normOrd{\,[u(w),[u(w),b]]_\beta a_\alpha(w)} \delta(z-w)\\
  & \quad + \sum_{\alpha,\beta \in \Delta_+^\mfrak{u}} [e_\beta,a]_\alpha [u(w),[e_\alpha,b]]_\beta \partial_w\delta(z-w)\\
  & \quad -\sum_{\alpha \in \Delta_+^\mfrak{u}} [u(w),a]_\alpha [e_\alpha,b](w) \delta(z-w) - \sum_{\alpha \in \Delta_+^\mfrak{u}} [u(z),a]_\alpha \kappa(e_\alpha,b)\partial_w \delta(z-w) \\
  &\quad -\sum_{\alpha \in \Delta_+^\mfrak{u}} a_\alpha^*(w)[a,[e_\alpha,b]](w)\delta(z-w) \\
  & \quad -\sum_{\alpha \in \Delta_+^\mfrak{u}} a_\alpha^*(w)(\kappa-\kappa_c^\mfrak{p})(a,[e_\alpha,b])\partial_w \delta(z-w)  +[a,b](w)\delta(z-w)
\end{align*}
if $b \in \widebar{\mfrak{u}}$. The coefficient in front of $\partial_w \delta(z-w)$ is equal to
\begin{align*}
  \sum_{\alpha,\beta \in \Delta_+^\mfrak{u}} [e_\beta,a]_\alpha [u(w),[e_\alpha,b]]_\beta - \sum_{\alpha \in \Delta_+^\mfrak{u}} \kappa(e_\alpha,b)[u(z),a]_\alpha -\sum_{\alpha \in \Delta_+^\mfrak{u}} (\kappa-\kappa_c^\mfrak{p})(a,[e_\alpha,b])a_\alpha^*(w).
\end{align*}
As we have
\begin{gather*}
\begin{aligned}
  \sum_{\alpha,\beta \in \Delta_+^\mfrak{u}} [e_\beta,a]_\alpha [u(w),[e_\alpha,b]]_\beta &= \sum_{\alpha \in \Delta_+^\mfrak{u}} [[u(w),[e_\alpha,b]],a]_\alpha = \sum_{\alpha \in \Delta_+^\mfrak{u}} [[e_\alpha,[u(w),b]],a]_\alpha \\
  &= -\kappa_c^\mfrak{p}(a,[u(w),b]),
\end{aligned} \\[2mm]
  \sum_{\alpha \in \Delta_+^\mfrak{u}} \kappa(e_\alpha,b)[u(z),a]_\alpha  = \kappa([u(z),a],b) = - \kappa(a,[u(z),b]),
\end{gather*}
we obtain that the coefficient is $\kappa(a,[u(z),b]) - \kappa(a,[u(w),b])= \kappa(u(w),[a,b]) - \kappa(u(z),[a,b])$, which gives us
\begin{align*}
  (\kappa(u(w),[a,b]) - \kappa(u(z),[a,b]))\partial_w \delta(z-w) = -\kappa(\partial_w u(w),[a,b]) \delta(z-w).
\end{align*}
On the other hand, the coefficient in front of $\delta(z-w)$ is equal to
\begin{gather*}
  -\sum_{\alpha,\beta \in \Delta_+^\mfrak{u}} \normOrd{\, [u(w),[e_\alpha,b]]_\beta [u(w),a]_\alpha a_\beta(w)} +{1\over 2} \sum_{\alpha,\beta \in \Delta_+^\mfrak{u}} [e_\beta,a]_\alpha \normOrd{\,[u(w),[u(w),b]]_\beta a_\alpha(w)} \\
  -\sum_{\alpha \in \Delta_+^\mfrak{u}} [u(w),a]_\alpha [e_\alpha,b](w) -\sum_{\alpha \in \Delta_+^\mfrak{u}} a_\alpha^*(w)[a,[e_\alpha,b]](w) + [a,b](w).
\end{gather*}
By using
\begin{align*}
  \sum_{\alpha \in \Delta_+^\mfrak{u}} [u(w),a]_\alpha [e_\alpha,b](w) = \sum_{\alpha \in \Delta_+^\mfrak{u}} a_\alpha^*(w)[[e_\alpha,a],b](w)
\end{align*}
and
\begin{gather*}
\begin{aligned}
    \sum_{\alpha \in \Delta_+^\mfrak{u}} [u(w),[e_\alpha,b]]_\beta [u(w),a]_\alpha & = [u(w),[[u(w),a],b]]_\beta = [[u(w),a], [u(w),b]]_\beta
\end{aligned}  \\[2mm]
\begin{aligned}
    \sum_{\beta \in \Delta_+^\mfrak{u}} [e_\beta,a]_\alpha [u(w),[u(w),b]]_\beta &= [[u(w),[u(w),b]], a]_\alpha \\[-5mm]
    &= [[u(w),a],[u(w),b]]_\alpha + [u(w),[[u(w),b],a]]_\alpha \\
    &= [[u(w),a],[u(w),b]]_\alpha + [u(w),[[u(w),a],b]]_\alpha + [u(w),[u(w),[b,a]]]_\alpha \\
    &= 2 [[u(w),a],[u(w),b]]_\alpha -  [u(w),[u(w),[a,b]]]_\alpha,
\end{aligned}
\end{gather*}
the coefficient can be simplified into the form
\begin{align*}
  - {1 \over 2} \sum_{\alpha \in \Delta_+^\mfrak{u}} \normOrd{\,[u(w),[u(w),[a,b]]]_\alpha a_\alpha(w)}  -\ad(u(w))([a,b](w)) + [a,b](w).
\end{align*}
Putting all together, we obtain that $[\Phi(a(z)),\Phi(b(w))] = \Phi([a,b](w)) \delta(z-w)$. Finally, for $a, b \in \widebar{\mfrak{u}}$ we may write
\begin{align*}
  [\Phi(a(z)),\Phi(b(w))] &= {1 \over 2}  \sum_{\alpha,\beta \in \Delta_+^\mfrak{u}} \normOrd{\,[u(w),[e_\alpha,b]]_\beta [u(w),[u(w),a]]_\alpha a_\beta(w)} \delta(z-w) \\
  & \quad -{1 \over 2}  \sum_{\alpha,\beta \in \Delta_+^\mfrak{u}} \normOrd{\,[u(w),[e_\alpha,a]]_\beta [u(w),[u(w),b]]_\alpha a_\beta(w)} \delta(z-w) \\
  & \quad - \sum_{\alpha,\beta \in \Delta_+^\mfrak{u}} [u(w),[e_\alpha,b]]_\beta  [u(z),[e_\beta,a]]_\alpha \partial_w\delta(z-w) \\
  & \quad + {1 \over 2} \sum_{\alpha \in \Delta_+^\mfrak{u}} [u(w),[u(w),a]]_\alpha [e_\alpha,b](w) \delta(z-w)\\
  & \quad - {1 \over 2} \sum_{\alpha \in \Delta_+^\mfrak{u}} [u(w),[u(w),b]]_\alpha [e_\alpha,a](w) \delta(z-w)\\
  & \quad + \sum_{\alpha,\beta \in \Delta_+^\mfrak{u}} a_\alpha^*(w)a_\beta^*(w) [[e_\alpha,a], [e_\beta,b]](w)\delta(z-w) \\
  & \quad - \sum_{\alpha \in \Delta_+^\mfrak{u}} a_\alpha^*(w) [[e_\alpha,a],b](w)\delta(z-w) - \sum_{\alpha \in \Delta_+^\mfrak{u}} a_\alpha^*(w) [a, [e_\alpha,b]](w)\delta(z-w) \\
  & \quad + {1 \over 2} \sum_{\alpha \in \Delta_+^\mfrak{u}} [u(z),[u(z),a]]_\alpha \kappa(e_\alpha,b) \partial_w \delta(z-w) \\
  & \quad + {1 \over 2} \sum_{\alpha \in \Delta_+^\mfrak{u}} [u(w),[u(w),b]]_\alpha \kappa(e_\alpha,a) \partial_w \delta(z-w) \\
  & \quad + \sum_{\alpha,\beta \in \Delta_+^\mfrak{u}} a_\alpha^*(z)a_\beta^*(w)(\kappa-\kappa_c^\mfrak{p})([e_\alpha,a], [e_\beta,b]) \partial_w\delta(z-w).
\end{align*}
We have
\begin{align*}
  \sum_{\alpha \in \Delta_+^\mfrak{u}} a_\alpha^*(w) [[e_\alpha,a],b](w) + \sum_{\alpha \in \Delta_+^\mfrak{u}} a_\alpha^*(w) [a, [e_\alpha,b]](w) = \sum_{\alpha \in \Delta_+^\mfrak{u}} a_\alpha^*(w) [e_\alpha,[a,b]](w) = 0.
\end{align*}
Further, since we have
\begin{align*}
    \sum_{\alpha \in \Delta_+^\mfrak{u}} [u(w),[u(w),a]]_\alpha [e_\alpha,b](w) &= \sum_{\beta,\gamma \in \Delta_+^\mfrak{u}} a_\beta^*(w)a_\gamma^*(w)[[e_\beta, [e_\gamma,a]],b](w) \\
    &=  \sum_{\beta,\gamma \in \Delta_+^\mfrak{u}} a_\beta^*(w)a_\gamma^*(w)([[e_\beta,b], [e_\gamma,a]](w) + [e_\beta, [[e_\gamma,a],b]](w)), \\
    \sum_{\alpha \in \Delta_+^\mfrak{u}} [u(w),[u(w),b]]_\alpha [e_\alpha,a](w) &= \sum_{\beta,\gamma \in \Delta_+^\mfrak{u}} a_\beta^*(w)a_\gamma^*(w)[[e_\beta, [e_\gamma,b]],a](w) \\
    &=  \sum_{\beta,\gamma \in \Delta_+^\mfrak{u}} a_\beta^*(w)a_\gamma^*(w)([[e_\beta,a], [e_\gamma,b]](w) + [e_\beta, [[e_\gamma,b],a]](w)),
\end{align*}
it implies that
\begin{gather*}
  {1 \over 2} \sum_{\alpha \in \Delta_+^\mfrak{u}} [u(w),[u(w),a]]_\alpha [e_\alpha,b](w) - {1 \over 2} \sum_{\alpha \in \Delta_+^\mfrak{u}} [u(w),[u(w),b]]_\alpha [e_\alpha,a](w) \\
  + \sum_{\alpha,\beta \in \Delta_+^\mfrak{u}} a_\alpha^*(w)a_\beta^*(w) [[e_\alpha,a], [e_\beta,b]](w)
\end{gather*}
is equal to $0$. By using
\begin{align*}
  \sum_{\alpha,\beta \in \Delta_+^\mfrak{u}} [u(w),[e_\alpha,b]]_\beta  [u(z),[e_\beta,a]]_\alpha &= \sum_{\beta \in \Delta_+^\mfrak{u}} [u(w),[[u(z),[e_\beta,a]],b]]_\beta \\
  &= \sum_{\beta \in \Delta_+^\mfrak{u}} [[e_\beta,[u(z),a]],[u(w),b]]_\beta
  = -\kappa_c^\mfrak{p}([u(z),a],[u(w),b]),
\end{align*}
we obtain that
\begin{align*}
  - \sum_{\alpha,\beta \in \Delta_+^\mfrak{u}} [u(w),[e_\alpha,b]]_\beta  [u(z),[e_\beta,a]]_\alpha -\sum_{\alpha,\beta \in \Delta_+^\mfrak{u}} a_\alpha^*(z)a_\beta^*(w) \kappa_c^\mfrak{p}([e_\alpha,a],[e_\beta,b])=0.
\end{align*}
We may write
\begin{align*}
  {1 \over 2}\sum_{\alpha \in \Delta_+^\mfrak{u}} [u(z),[u(z),a]]_\alpha \kappa(e_\alpha,b) \partial_w \delta(z-w) &= {1 \over 2} \kappa([u(z),[u(z),a]],b)\partial_w \delta(z-w) \\
  & = -{1 \over 2}\kappa([u(z),a],[u(z),b])\partial_w \delta(z-w), \\[2mm]
  {1 \over 2}\sum_{\alpha \in \Delta_+^\mfrak{u}} [u(w),[u(w),b]]_\alpha \kappa(e_\alpha,a) \partial_w \delta(z-w) &= {1 \over 2}\kappa([u(w),[u(w),b]],a)\partial_w \delta(z-w) \\
  & = -{1 \over 2}\kappa([u(w),a],[u(w),b])\partial_w \delta(z-w), \\[2mm]
  \sum_{\alpha,\beta \in \Delta_+^\mfrak{u}} a_\alpha^*(z)a_\beta^*(w)\kappa([e_\alpha,a], [e_\beta,b]) \partial_w\delta(z-w) &= \kappa([u(z),a],[u(w),b])\partial_w\delta(z-w).
\end{align*}
By adding all three equations above together, we get for the right hand side
\begin{align*}
  \bigg(\!\!-{1 \over 2}\kappa([u(z),a],[u(z),b]) -{1 \over 2}\kappa([u(w),a],[u(w),b]) + \kappa([u(z),a],[u(w),b])\!\!\bigg)\partial_w\delta(z-w),
\end{align*}
which is equal to $0$. The last step is to show that
\begin{align*}
   \sum_{\alpha \in \Delta_+^\mfrak{u}} [u(w),[e_\alpha,b]]_\beta [u(w),[u(w),a]]_\alpha -  \sum_{\alpha \in \Delta_+^\mfrak{u}} [u(w),[e_\alpha,a]]_\beta [u(w),[u(w),b]]_\alpha =0.
\end{align*}
Indeed, we have
\begin{align*}
  \sum_{\alpha \in \Delta_+^\mfrak{u}} [u(w),[e_\alpha,b]]_\beta [u(w),[u(w),a]]_\alpha &= [u(w),[[u(w),[u(w),a]],b]]_\beta \\[-5mm]
  &= [[u(w),[u(w),a]],[u(w),b]]_\beta, \\[2mm]
  -\sum_{\alpha \in \Delta_+^\mfrak{u}} [u(w),[e_\alpha,a]]_\beta [u(w),[u(w),b]]_\alpha &= -[u(w),[[u(w),[u(w),b]],a]]_\beta \\[-5mm]
  &= [[u(w),a],[u(w),[u(w),b]]]_\beta,
\end{align*}
which together with the fact
\begin{align*}
  0 = (\ad(u(w)))^3([a,b]) = 3[[u(w),[u(w),a]],[u(w),b]]_\beta + 3[[u(w),a],[u(w),[u(w),b]]]_\beta
\end{align*}
gives us the required statement. This finishes the proof.}

In the next, we use the explicit form of the Feigin--Frenkel homomorphism
\begin{align*}
  w_{\kappa,\mfrak{g}}^\mfrak{p} \colon \mcal{V}^\kappa(\mfrak{g}) \rarr \mcal{M}_\mfrak{u} \otimes_\C \mcal{V}^{\kappa - \kappa_c^\mfrak{p}}(\mfrak{p})
\end{align*}
of $\N_0$-graded vertex algebras to get some information about the simple affine vertex algebra $\mcal{V}_\kappa(\mfrak{g})$. In particular, we will focus on determining the associated variety $X_{\smash{\mcal{V}_\kappa}(\mfrak{g})}$ of $\mcal{V}_\kappa(\mfrak{g})$. As the decomposition $\mfrak{p} = \mfrak{l} \oplus \widebar{\mfrak{u}}$ is orthogonal with respect to $\kappa - \kappa_c^\mfrak{p}$, we obtain a surjective homomorphism $\mcal{V}^{\kappa-\kappa_c^\mfrak{p}}(\mfrak{p}) \rarr  \mcal{V}^{\kappa-\kappa_c^\mfrak{p}}(\mfrak{l})$
of $\N_0$-graded vertex algebras. In order to simplify the following considerations, we will assume that a $\mfrak{g}$-invariant symmetric bilinear form $\kappa$ on $\mfrak{g}$ satisfies the additional condition $(\kappa-\kappa_c^\mfrak{p})_{|[\mfrak{l},\mfrak{l}]} = 0$. This gives us a canonical surjective homomorphism
\begin{align*}
  \mcal{V}^{\kappa-\kappa_c^\mfrak{p}}(\mfrak{l}) \rarr \mcal{V}^{\kappa-\kappa_c^\mfrak{p}}(\mfrak{z}(\mfrak{l}))
\end{align*}
of $\N_0$-graded vertex algebras and together with the equality $\kappa_c^\mfrak{p}{}_{|\mfrak{z}(\mfrak{l})}= \kappa_c{}_{|\mfrak{z}(\mfrak{l})} = -{1\over 2}\kappa_\mfrak{g}{}_{|\mfrak{z}(\mfrak{l})}$ gives rise to a homomorphism
\begin{align}
  \widetilde{w}_{\kappa,\mfrak{g}}^\mfrak{p} \colon \mcal{V}^\kappa(\mfrak{g}) \rarr \mcal{M}_\mfrak{u} \otimes_\C \mcal{V}^{\kappa+{1\over 2} \kappa_\mfrak{g}}\!(\mfrak{z}(\mfrak{l}))
\end{align}
of $\N_0$-graded vertex algebras, which coincides with the vertex algebra homomorphism defined by \eqref{eq:gen.Wakimoto} as follows from Proposition \ref{prop:cdo formula} and Theorem \ref{thm:FF homomorphism 1-graded}. The existence of such homomorphism ensures that the diagrams
\begin{align}\label{eq:Zhu algebra diagram}
\bfig
\square(0,0)|alrb|<900,500>[\mcal{V}^\kappa(\mfrak{g})`\mcal{M}_\mfrak{u} \otimes_\C \mcal{V}^{\kappa+{1\over 2}\kappa_\mfrak{g}}\! (\mfrak{z}(\mfrak{l}))`U(\mfrak{g})`\eus{A}_\mfrak{u} \otimes_\C U(\mfrak{z}(\mfrak{l}));\widetilde{w}_{\kappa,\mfrak{g}}^\mfrak{p}`\pi_{\rm Zhu}`\pi_{\rm Zhu}`\widetilde{\pi}_\mfrak{g}^\mfrak{p}]
\efig
\qquad \text{and} \qquad
\bfig
\square(0,0)|alrb|<900,500>[\mcal{V}^\kappa(\mfrak{g})`\mcal{M}_\mfrak{u} \otimes_\C \mcal{V}^{\kappa+{1\over 2}\kappa_\mfrak{g}}\! (\mfrak{z}(\mfrak{l}))`S(\mfrak{g})`\C{[\mfrak{u} \oplus \mfrak{u}^*]} \otimes_\C S(\mfrak{z}(\mfrak{l}));\widetilde{w}_{\kappa,\mfrak{g}}^\mfrak{p}`\sigma_{\rm Zhu}`\sigma_{\rm Zhu}`\gr\widetilde{\pi}_\mfrak{g}^\mfrak{p}]
\efig
\end{align}
are commutative. The homomorphism $\widetilde{\pi}_\mfrak{g}^\mfrak{p} \colon U(\mfrak{g}) \rarr \eus{A}_\mfrak{u} \otimes_\C U(\mfrak{z}(\mfrak{l}))$ of Zhu's algebras is given by
\begin{align}
  \widetilde{\pi}_\mfrak{g}^\mfrak{p}(a) = - \sum_{\alpha\in \Delta_+^\mfrak{u}}  \bigg[{\ad(u(x)) e^{\ad(u(x))} \over e^{\ad(u(x))}- \id}\, (e^{-\ad(u(x))}a)_\mfrak{u}\bigg]_\alpha \partial_{x_\alpha} + (e^{-\ad(u(x))} a)_{\mfrak{z}(\mfrak{l})} \label{eq:pi_X homomorphism}
\end{align}
for $a \in \mfrak{g}$, where $[a]_\alpha$ denotes the $\alpha$-th coordinate of $a \in \mfrak{u}$ with respect to the basis $\{e_\alpha;\, \alpha \in \Delta_+^\mfrak{u}\}$ of $\mfrak{u}$, the elements $a_\mfrak{u}$ and $a_{\mfrak{z}(\mfrak{l})}$ are the $\mfrak{u}$-part and $\mfrak{z}(\mfrak{l})$-part of $a \in \mfrak{g}$ with respect to the decomposition $\mfrak{g}=\widebar{\mfrak{u}} \oplus \mfrak{z}(\mfrak{l}) \oplus [\mfrak{l},\mfrak{l}] \oplus \mfrak{u}$, and $u(x) \in \mfrak{u} \otimes_\C \C[\mfrak{u}]$ is given by
\begin{align*}
u(x)=\sum_{\alpha \in \Delta_+^\mfrak{u}}  x_\alpha e_\alpha.
\end{align*}
Let us note that $\mfrak{g} \otimes_\C \C[\mfrak{u}]$ has the natural structure of a Lie algebra. Hence, we have a well-defined linear mapping $\ad(u(x)) \colon \mfrak{g} \otimes_\C \C[\mfrak{u}] \rarr \mfrak{g} \otimes_\C \C[\mfrak{u}]$. Let us recall that the expression \eqref{eq:pi_X homomorphism} coincides with the formula derived in \cite{Krizka-Somberg2017} which will be important later on.
\medskip

Since $\widetilde{w}_{\kappa,\mfrak{g}}^\mfrak{p}$ is a homomorphism of $\N_0$-graded vertex algebras, we may introduce a new $\N_0$-graded vertex algebra by
\begin{align*}
  \mcal{N}^\kappa_\mfrak{p}(\mfrak{g}) = \mcal{V}^\kappa(\mfrak{g})/\mcal{I}^\kappa_\mfrak{p}(\mfrak{g}),
\end{align*}
where $\mcal{I}^\kappa_\mfrak{p}(\mfrak{g}) =  \ker \widetilde{w}_{\kappa, \mfrak{g}}^\mfrak{p}$ is a vertex algebra ideal of $\mcal{V}^\kappa(\mfrak{g})$. We denote by
\begin{align*}
  I_\mfrak{p}(\mfrak{g}) = \ker \widetilde{\pi}^\mfrak{p}_\mfrak{g} \qquad \text{and} \qquad J_\mfrak{p}(\mfrak{g}) = \ker \gr \widetilde{\pi}^\mfrak{p}_\mfrak{g}
\end{align*}
ideals of $U(\mfrak{g})$ and $S(\mfrak{g})$, respectively. Then by using the commutative diagrams in \eqref{eq:Zhu algebra diagram} we obtain easily that
\begin{align*}
  \pi_{\rm Zhu}(\mcal{I}^\kappa_\mfrak{p}(\mfrak{g})) \subset I_\mfrak{p}(\mfrak{g}) \qquad \text{and} \qquad \sigma_{\rm Zhu}(\mcal{I}^\kappa_\mfrak{p}(\mfrak{g})) \subset J_\mfrak{p}(\mfrak{g}).
\end{align*}
Besides, we have
\begin{align*}
  A(\mcal{N}^\kappa_\mfrak{p}(\mfrak{g})) \simeq U(\mfrak{g})/ \pi_{\rm Zhu}(\mcal{I}^\kappa_\mfrak{p}(\mfrak{g})) \qquad \text{and} \qquad R_{\smash{\mcal{N}^\kappa_\mfrak{p}}\!(\mfrak{g})} \simeq S(\mfrak{g})/ \sigma_{\rm Zhu}(\mcal{I}^\kappa_\mfrak{p}(\mfrak{g})).
\end{align*}
In the following theorem we determine the associated variety $X_{\smash{\mcal{N}^\kappa_\mfrak{p}}(\mfrak{g})}$ of $\mcal{N}^\kappa_\mfrak{p}(\mfrak{g})$.
\medskip

\theorem{\label{thm:associated variety}Let $\kappa$ be a $\mfrak{g}$-invariant symmetric bilinear form on $\mfrak{g}$ such that $(\kappa-\kappa_c^\mfrak{p})_{|[\mfrak{l},\mfrak{l}]}=0$. Then we have
\begin{align*}
  X_{\smash{\mcal{N}^\kappa_\mfrak{p}}\!(\mfrak{g})} \simeq \widebar{\mcal{S}^*_{\mfrak{p}}},
\end{align*}
where $\mcal{S}_{\mfrak{p}}$ is the Dixmier sheet of $\mfrak{g}$ attached to $\mfrak{p}$. In addition, the quotient affine vertex algebra $\mcal{N}^\kappa_\mfrak{p}(\mfrak{g})$ is a quantization of the infinite jet scheme $J_\infty(\widebar{\mcal{S}^*_{\mfrak{p}}})$ of $\widebar{\mcal{S}^*_{\mfrak{p}}}$ in the sense that
\begin{align*}
  SS(\mcal{N}^\kappa_\mfrak{p}(\mfrak{g})) \simeq J_\infty(\widebar{\mcal{S}^*_{\mfrak{p}}})
\end{align*}
as topological spaces, i.e.\ we have $SS(\mcal{N}^\kappa_\mfrak{p}(\mfrak{g}))_{\rm red} \simeq J_\infty(\widebar{\mcal{S}^*_{\mfrak{p}}})_{\rm red}$.}

\proof{Let us consider the conformal weight filtrations on  $\mcal{V}^\kappa(\mfrak{g})$ and $\smash{\mcal{M}_\mfrak{u} \otimes_\C \mcal{V}^{\kappa+{1\over 2} \kappa_\mfrak{g}}\!(\mfrak{z}(\mfrak{l}))}$, which are preserved by $\widetilde{w}_{\kappa,\mfrak{g}}^\mfrak{p}$. Hence, we obtain a homomorphism
\begin{align*}
   \gr^G\! \widetilde{w}_{\kappa,\mfrak{g}}^\mfrak{p} \colon \gr^G\! \mcal{V}^\kappa(\mfrak{g}) \rarr \gr^G\! \mcal{M}_\mfrak{u} \otimes_\C \gr^G\! \mcal{V}^{\kappa+{1\over 2} \kappa_\mfrak{g}}\!(\mfrak{z}(\mfrak{l}))
\end{align*}
of Poisson vertex algebras. As we have the isomorphisms
\begin{align*}
   \gr^G\!\mcal{V}^\kappa(\mfrak{g}) &\simeq \C[J_\infty(\mfrak{g}^*)], \\
   \gr^G\! \mcal{M}_\mfrak{u} \otimes_\C \gr^G\! \mcal{V}^{\kappa+{1\over 2} \kappa_\mfrak{g}}\!(\mfrak{z}(\mfrak{l})) &\simeq  \C[J_\infty(\mfrak{u} \oplus \mfrak{u}^*)] \otimes_\C \C[J_\infty(\mfrak{z}(\mfrak{l})^*)]
\end{align*}
of Poisson vertex algebras, the homomorphism  $\gr^G\! \widetilde{w}_{\kappa,\mfrak{g}}^\mfrak{p}$ corresponds to a morphism of infinite jet schemes
\begin{align*}
 J_\infty(j \circ \widetilde{\mu}_P) \colon J_\infty(U \times [\mfrak{p},\mfrak{p}]^\perp)  \rarr J_\infty(\mfrak{g}),
\end{align*}
where we identified $\mfrak{g}$ with $\mfrak{g}^*$ by using the non-degenerate $\mfrak{g}$-invariant symmetric bilinear form $\kappa_\omega$ on $\mfrak{g}$. The morphism $j \colon \mfrak{g} \rarr \mfrak{g}$ is given by $j(a)=-a$ for $a \in \mfrak{g}$ and the morphism
\begin{align*}
  \widetilde{\mu}_P \colon U \times [\mfrak{p},\mfrak{p}]^\perp \rarr \mfrak{g}
\end{align*}
is the restriction of the generalized Grothendieck's simultaneous resolution $\mu_P \colon G \times_P [\mfrak{p},\mfrak{p}]^\perp \rarr \mfrak{g}$ defined by $\mu_P(g,a) = \Ad(g)(a)$ for $g \in G$ and $a \in [\mfrak{p},\mfrak{p}]^\perp$ to the open dense subset $UP \times_P [\mfrak{p},\mfrak{p}]^\perp$ of $G \times_P [\mfrak{p},\mfrak{p}]^\perp$, see the formula \eqref{eq:Grothendieck resolusion}. Since the image of $\mu_P$ is the closure of the Dixmier sheet $\mcal{S}_\mfrak{p}$,
we may consider two morphisms
\begin{align*}
  p \colon U \times [\mfrak{p},\mfrak{p}]^\perp \rarr \widebar{\mcal{S}_{\mfrak{p}}} \qquad \text{and} \qquad i \colon \widebar{\mcal{S}_{\mfrak{p}}} \rarr \mfrak{g}
\end{align*}
of algebraic varieties defined by $p(g,a) = \Ad(g)(a)$ for $g \in G$ and $a \in [\mfrak{p},\mfrak{p}]^\perp$ and $i(a) = -a$ for $a \in \smash{\widebar{\mcal{S}_{\mfrak{p}}}}$.
It easily follows that $i$ is a closed embedding and that $p$ is a dominant morphism which is generically one-to-one. As we have $j \circ \widetilde{\mu} = i \circ p$, we may write
\begin{align*}
  J_\infty(j \circ \widetilde{\mu}) = J_\infty(i) \circ J_\infty(p),
\end{align*}
where $J_\infty(i) \colon J_\infty(\widebar{\mcal{S}_{\mfrak{p}}}) \rarr J_\infty(\mfrak{g})$ is a closed embedding and $J_\infty(p) \colon J_\infty(U \times [\mfrak{p},\mfrak{p}]^\perp) \rarr J_\infty(\widebar{\mcal{S}_{\mfrak{p}}})$ is a dominant morphism by \cite[Lemma 1.8]{Arakawa-Moreau2023-book}. Therefore, we have that
\begin{align*}
  \im \gr^G\! \widetilde{w}_{\kappa,\mfrak{g}}^\mfrak{p} \simeq \C[J_\infty(\widebar{\mcal{S}^*_{\mfrak{p}}})]
\end{align*}
as Poisson vertex algebras.

Further, let us recall that on the $\N_0$-graded vertex subalgebra $\im \widetilde{w}_{\kappa,\mfrak{g}}^\mfrak{p}$ of $\smash{\mcal{M}_\mfrak{u} \otimes_\C \mcal{V}^{\kappa+{1\over 2} \kappa_\mfrak{g}}\!(\mfrak{z}(\mfrak{l}))}$ we have two different filtrations defined as
\begin{align*}
  G_p \im \widetilde{w}_{\kappa,\mfrak{g}}^\mfrak{p} = \widetilde{w}_{\kappa,\mfrak{g}}^\mfrak{p}(G_p \mcal{V}^\kappa(\mfrak{g}))
  \qquad \text{and} \qquad
  H_p \im \widetilde{w}_{\kappa,\mfrak{g}}^\mfrak{p} = \im \widetilde{w}_{\kappa,\mfrak{g}}^\mfrak{p} \cap G_p(\mcal{M}_\mfrak{u} \otimes_\C \mcal{V}^{\kappa+{1\over 2} \kappa_\mfrak{g}}\!(\mfrak{z}(\mfrak{l})))
\end{align*}
for $p \in \Z$. Both filtrations are good in the sense of \cite{Li2004} and $\gr^H\! \im \widetilde{w}_{\kappa,\mfrak{g}}^\mfrak{p} \simeq \im \gr^G\! \widetilde{w}_{\kappa,\mfrak{g}}^\mfrak{p}$ as Poisson vertex algebras. Besides, we have a surjective homomorphism
\begin{align*}
  \gr^G\! \im \widetilde{w}_{\kappa,\mfrak{g}}^\mfrak{p} \rarr \im \gr^G\! \widetilde{w}_{\kappa,\mfrak{g}}^\mfrak{p}
\end{align*}
of Poisson vertex algebras, which means that $\im \gr^G\! \widetilde{w}_{\kappa,\mfrak{g}}^\mfrak{p}$ is a finitely generated $\gr^G\! \im \widetilde{w}_{\kappa,\mfrak{g}}^\mfrak{p}$-module. Therefore, the corresponding singular supports satisfy
\begin{align*}
  SS(\gr^G\! \im \widetilde{w}_{\kappa,\mfrak{g}}^\mfrak{p}) = SS(\im \gr^G\! \widetilde{w}_{\kappa,\mfrak{g}}^\mfrak{p})
\end{align*}
since the singular support does not depend on the choice of a good filtration. As we have $\mcal{N}^\kappa_\mfrak{p}(\mfrak{g}) \simeq \im \widetilde{w}_{\kappa,\mfrak{g}}^\mfrak{p}$, we obtain
\begin{align*}
  SS(\mcal{N}^\kappa_\mfrak{p}(\mfrak{g})) \simeq J_\infty(\widebar{\mcal{S}^*_{\mfrak{p}}})
\end{align*}
as topological spaces, which gives us that
\begin{align*}
   X_{\smash{\mcal{N}^\kappa_\mfrak{p}}\!(\mfrak{g})} \simeq \pi_{\mfrak{g}^*\!,\infty}(J_\infty(\widebar{\mcal{S}^*_{\mfrak{p}}})) = \widebar{\mcal{S}^*_{\mfrak{p}}},
\end{align*}
where $\pi_{\mfrak{g}^*\!,\infty} \colon J_\infty(\mfrak{g}^*) \rarr \mfrak{g}^*$ is the canonical projection.}

\corollary{\label{co:Ass-var}Let $\kappa$ be a $\mfrak{g}$-invariant symmetric bilinear form on $\mfrak{g}$ such that $(\kappa-\kappa_c^\mfrak{p})_{|[\mfrak{l},\mfrak{l}]}=0$. Then we have $X_{\mcal{V}_\kappa(\mfrak{g})}\subset \smash{\widebar{\mcal{S}^*_{\mfrak{p}}}}$.}

\proposition{\label{prop:ideal Ip}We have
\begin{align*}
  I_\mfrak{p}(\mfrak{g}) = \ker \Phi_Y = \bigcap_{\lambda \in \smash{\mfrak{z}(\mfrak{l})^*}} \Ann_{U(\mfrak{g})}\! M^\mfrak{g}_\mfrak{p}(\lambda)
\end{align*}
and $\mcal{V}(I_\mfrak{p}(\mfrak{g})) = \widebar{\mcal{S}^*_{\mfrak{p}}}$, where $\mcal{S}_\mfrak{p}$ is the Dixmier sheet attached to $\mfrak{p}$.}

\proof{By the results of \cite{Krizka-Somberg2017} we have
\begin{align*}
  \ker \Phi^\lambda_X = \ker\, ((\id_{\eus{A}_\mfrak{u}}\!
  \otimes f_{\lambda+\rho_\mfrak{p}}) \circ \widetilde{\pi}^\mfrak{p}_\mfrak{g}),
\end{align*}
where $f_{\lambda+\rho_\mfrak{p}} \colon U(\mfrak{z}(\mfrak{l})) \rarr \C$ for $\lambda \in \mfrak{z}(\mfrak{l})^*$ is the evaluation homomorphism at $\lambda+\rho_\mfrak{p}$ since we may identify $U(\mfrak{z}(\mfrak{l})) \simeq S(\mfrak{z}(\mfrak{l}))$ with the algebra of polynomial functions on $\mfrak{z}(\mfrak{l})^*$.
Hence, we easily get that
\begin{align*}
  \bigcap_{\lambda \in \smash{\mfrak{z}(\mfrak{l})^*}} \ker \Phi_X^\lambda = \bigcap_{\lambda \in \smash{\mfrak{z}(\mfrak{l})^*}} \ker\, ((\id_{\eus{A}_\mfrak{u}}\!
  \otimes f_{\lambda+\rho_\mfrak{p}}) \circ \widetilde{\pi}^\mfrak{p}_\mfrak{g}) = \ker \widetilde{\pi}_\mfrak{g}^\mfrak{p} = I_\mfrak{p}(\mfrak{g}),
\end{align*}
which by using \eqref{eq:ker Phi_Y} gives us
\begin{align*}
  I_\mfrak{p}(\mfrak{g}) = \ker \Phi_Y = \bigcap_{\lambda \in \smash{\mfrak{z}(\mfrak{l})^*}} \Ann_{U(\mfrak{g})}\! M^\mfrak{g}_\mfrak{p}(\lambda).
\end{align*}
The rest of the proof follows from the formula \eqref{eq:ker Phi_Y variety}.}

\proposition{Let $\kappa$ be a $\mfrak{g}$-invariant symmetric bilinear form on $\mfrak{g}$ such that $(\kappa-\kappa_c^\mfrak{p})_{|[\mfrak{l},\mfrak{l}]}=0$. Then we have
\begin{align*}
  X_{\smash{\mcal{N}^\kappa_\mfrak{p}\!(\mfrak{g})}} \simeq \Specm (\gr A(\mcal{N}^\kappa_\mfrak{p}(\mfrak{g}))),
\end{align*}
or in other words $\sqrt{\sigma_{\rm Zhu}(\mcal{I}^\kappa_\mfrak{p}(\mfrak{g}))} = \sqrt{\gr \pi_{\rm Zhu}(\mcal{I}^\kappa_\mfrak{p}(\mfrak{g}))}$.}

\proof{Let us recall that the following diagram
\begin{align*}
\bfig
\square(0,0)|alrb|<1200,500>[S(\mfrak{g})`S(\mfrak{g})`S(\mfrak{g})/ \sigma_{\rm Zhu}(\mcal{I}^\kappa_\mfrak{p}(\mfrak{g}))`S(\mfrak{g})/\gr \pi_{\rm Zhu}(\mcal{I}^\kappa_\mfrak{p}(\mfrak{g})); \eta_{\mcal{V}^\kappa\!(\mfrak{g})}=\id_{S(\mfrak{g})}``` \eta_{\mcal{N}^\kappa_\mfrak{p}\!(\mfrak{g})}]
\efig
\end{align*}
is commutative, which implies that $\sigma_{\rm Zhu}(\mcal{I}^\kappa_\mfrak{p}(\mfrak{g})) \subset \gr \pi_{\rm Zhu}(\mcal{I}^\kappa_\mfrak{p}(\mfrak{g}))$. On the other hand, as a consequence of Theorem \ref{thm:associated variety} and Proposition \ref{prop:ideal Ip} we have
\begin{align*}
  \sqrt{\sigma_{\rm Zhu}(\mcal{I}^\kappa_\mfrak{p}(\mfrak{g}))} = \sqrt{\gr I_\mfrak{p}(\mfrak{g})},
\end{align*}
which together with the fact $\gr \pi_{\rm Zhu}(\mcal{I}^\kappa_\mfrak{p}(\mfrak{g})) \subset \gr I_\mfrak{p}(\mfrak{g})$ implies that
\begin{align*}
  \sqrt{\gr \pi_{\rm Zhu}(\mcal{I}^\kappa_\mfrak{p}(\mfrak{g}))} \subset \sqrt{\gr I_\mfrak{p}(\mfrak{g})} =  \sqrt{\sigma_{\rm Zhu}(\mcal{I}^\kappa_\mfrak{p}(\mfrak{g}))}.
\end{align*}
This gives us the required statement.}

In addition, if $\mfrak{g}$ is a simple Lie algebra, then we have $\kappa = k\kappa_0$ for some $k \in \C$. Therefore, by using Proposition \ref{prop:critical levels} and the direct sum decomposition \eqref{eq:othogonal decomposition} of the Levi subalgebra $\mfrak{l}$ we can easily determine all levels $\kappa$ for which $(\kappa-\kappa_c^\mfrak{p})_{|[\mfrak{l},\mfrak{l}]} = 0$. The classification of simple Lie algebras and their opposite standard parabolic subalgebras with commutative nilradical is in Table \ref{tab:parabolic commutative nilradical}. For a simple root $\alpha \in \Pi$, we denote by $\mfrak{p}_\alpha$ the maximal opposite standard parabolic subalgebra of $\mfrak{g}$ associated to the subset $\Sigma = \Pi \setminus \{\alpha\}$ of $\Pi$.

\begin{table}[ht]
\centering
\renewcommand{\arraystretch}{2.5}
\begin{tabular}{|c|c|c|c|c|c|}
  \hline
    $\mfrak{g}$ & Dynkin diagram & $\theta$ & $\mfrak{p}$ & $[\mfrak{l}, \mfrak{l}]$ & $\mfrak{u}$ \\
  \hline
    $\mfrak{a}_n$, $n\geq 1$ & $\dynkin[edge length=7mm, labels*={1,2,n-1,n}, label macro*/.code={\alpha_{#1}}, x/.style={thin}]{A}{oo..oo}$ & $\omega_1 + \omega_n$ & $\mfrak{p}_{\alpha_k}$ & $\mfrak{a}_{k-1} \oplus \mfrak{a}_{n-k}$ & $(\omega_1, \omega_{n-k})$ \\
  \hline
    $\mfrak{b}_2$ & $\dynkin[edge length=7mm, labels*={1,2}, label macro*/.code={\alpha_{#1}}, x/.style={thin}]{B}{oo}$ & $2\omega_2$ & $\mfrak{p}_{\alpha_1}$ & $\mfrak{a}_1$ & $2\omega_1$ \\
  \hline
    $\mfrak{b}_n$, $n\geq 3$ & $\dynkin[edge length=7mm, labels*={1,2,n-1,n}, label macro*/.code={\alpha_{#1}}, x/.style={thin}]{B}{oo.oo}$ & $\omega_2$ & $\mfrak{p}_{\alpha_1}$ & $\mfrak{b}_{n-1}$ & $\omega_1$ \\
  \hline
    $\mfrak{c}_n$, $n \geq 2$ & $\dynkin[edge length=7mm, labels*={1,2,n-1,n}, label macro*/.code={\alpha_{#1}}, x/.style={thin}]{C}{oo.oo}$ & $2\omega_1$ & $\mfrak{p}_{\alpha_n}$ & $\mfrak{a}_{n-1}$ & $2\omega_1$ \\
  \hline
    $\mfrak{d}_n$, $n\geq 4$ & $\dynkin[edge length=7mm, labels*={1,2}, label macro*/.code={\alpha_{#1}}, labels={,,n-2,n-1,n}, label macro/.code={\alpha_{#1}}, x/.style={thin}]{D}{oo.ooo}$ & $\omega_2$ & $\mfrak{p}_{\alpha_1}$ & $\mfrak{d}_{n-1}$ & $\omega_1$ \\
  \hline
    $\mfrak{d}_n$, $n\geq 4$ & $\dynkin[edge length=7mm, labels*={1,2}, label macro*/.code={\alpha_{#1}}, labels={,,n-2,n-1,n}, label macro/.code={\alpha_{#1}}, x/.style={thin}]{D}{oo.ooo}$ & $\omega_2$ & $\mfrak{p}_{\alpha_n}$, $\mfrak{p}_{\alpha_{n-1}}$ & $\mfrak{a}_{n-1}$ & $\omega_2$ \\
  \hline
    $\mfrak{e}_6$ & $\dynkin[edge length=7mm, labels={1,2,3,4,5,6}, label macro/.code={\alpha_{#1}}, ordering=Dynkin, x/.style={thin}]{E}{oooooo}$ & $\omega_6$ & $\mfrak{p}_{\alpha_1}$, $\mfrak{p}_{\alpha_5}$ & $\mfrak{d}_5$ & $\omega_5$ \\
  \hline
    $\mfrak{e}_7$ & $\dynkin[edge length=7mm, labels={1,2,3,4,5,6,7}, label macro/.code={\alpha_{#1}}, ordering=Dynkin, x/.style={thin}]{E}{ooooooo}$ & $\omega_1$ & $\mfrak{p}_{\alpha_6}$ & $\mfrak{e}_6$ & $\omega_1$ \\
  \hline
\end{tabular}
\caption{Parabolic subalgebras with commutative nilradical}
\label{tab:parabolic commutative nilradical}
\vspace{-2mm}
\end{table}

If $[\mfrak{l},\mfrak{l}] = 0$, then there is no restriction and we obtain $k \in \C$. This corresponds to $\mfrak{g}=\mfrak{a}_1$ with the opposite standard Borel subalgebra. This case will be excluded from our next considerations. On the other hand, if $[\mfrak{l},\mfrak{l}] = \bigoplus_{i=1}^r\! \mfrak{l}_i$ is a semisimple Lie algebra, we have the compatibility condition
\begin{align*}
   - {\ind_D(\mfrak{l}_i,\mfrak{u})h^\vee \over \ind_D(\mfrak{l}_i,\mfrak{u}) + h^\vee_i} = - {\ind_D(\mfrak{l}_j,\mfrak{u})h^\vee \over \ind_D(\mfrak{l}_j,\mfrak{u}) + h^\vee_j},
\end{align*}
or equivalently $\ind_D(\mfrak{l}_j,\mfrak{u}) h_i^\vee = \ind_D(\mfrak{l}_i,\mfrak{u})h_j^\vee$, for $1 \leq i,j \leq r$, which gives us
\begin{align*}
  k = - {\ind_D(\mfrak{l}_i,\mfrak{u})h^\vee \over \ind_D(\mfrak{l}_i,\mfrak{u}) + h^\vee_i}
\end{align*}
for any $i=1,2,\dots,r$. Hence, we need to know the Dynkin index of the $\mfrak{l}_i$-module $\mfrak{u}$ whose highest weight is given in Table \ref{tab:parabolic commutative nilradical}. By using the formula \eqref{eq:Dynkin index} for the Dynkin index, we get the complete list of simple Lie algebras and their opposite standard parabolic subalgebras with commutative nilradical satisfying $(\kappa - \kappa_c^\mfrak{p})_{|[\mfrak{l},\mfrak{l}]} = 0$ for some level $\kappa=k\kappa_0$, see Table \ref{tab:parabolic commutative nilradical level}. Let us note that the critical level is given by $\kappa_c = -h^\vee\kappa_0$.

\begin{table}[ht]
\centering
\renewcommand{\arraystretch}{2.5}
\begin{tabular}{|c|c|c|c|c|c|c|}
  \hline
    $\mfrak{g}$ & Dynkin diagram & $\mfrak{p}$ & $h^\vee$ & $k$ & $2\dim \mfrak{u}$ \\
  \hline
    $\mfrak{a}_n$, $n\geq 2$ & $\dynkin[edge length=7mm, labels*={1,2,n-1,n}, label macro*/.code={\alpha_{#1}}, x/.style={thin}]{A}{oo..oo}$ & $\mfrak{p}_{\alpha_1}$, $\mfrak{p}_{\alpha_n}$ & $n+1$ & $-1$ & $2n$ \\
  \hline
    $\mfrak{a}_{2n-1}$, $n\geq 2$ & $\dynkin[edge length=7mm, labels*={1,2,2n-2,2n-1}, label macro*/.code={\alpha_{#1}}, x/.style={thin}]{A}{oo..oo}$ & $\mfrak{p}_{\alpha_n}$ & $2n$ & $-n$ & $2n^2$ \\
  \hline
    $\mfrak{b}_n$, $n\geq 2$ & $\dynkin[edge length=7mm, labels*={1,2,n-1,n}, label macro*/.code={\alpha_{#1}}, x/.style={thin}]{B}{oo.oo}$ & $\mfrak{p}_{\alpha_1}$ & $2n-1$ & $-2$ & $4n-2$ \\
  \hline
    $\mfrak{c}_n$, $n \geq 2$ & $\dynkin[edge length=7mm, labels*={1,2,n-1,n}, label macro*/.code={\alpha_{#1}}, x/.style={thin}]{C}{oo.oo}$ & $\mfrak{p}_{\alpha_n}$ & $n+1$ & $-{n \over 2}-1$ & $n^2+n$ \\
  \hline
    $\mfrak{d}_n$, $n\geq 4$ & $\dynkin[edge length=7mm, labels*={1,2}, label macro*/.code={\alpha_{#1}}, labels={,,n-2,n-1,n}, label macro/.code={\alpha_{#1}}, x/.style={thin}]{D}{oo.ooo}$ & $\mfrak{p}_{\alpha_1}$ & $2n-2$ & $-2$ & $4n-4$ \\
  \hline
    $\mfrak{d}_n$, $n\geq 4$ & $\dynkin[edge length=7mm, labels*={1,2}, label macro*/.code={\alpha_{#1}}, labels={,,n-2,n-1,n}, label macro/.code={\alpha_{#1}}, x/.style={thin}]{D}{oo.ooo}$ & $\mfrak{p}_{\alpha_n}$, $\mfrak{p}_{\alpha_{n-1}}$ & $2n-2$ & $-n+2$ & $n^2-n$ \\
  \hline
    $\mfrak{e}_6$ & $\dynkin[edge length=7mm, labels={1,2,3,4,5,6}, label macro/.code={\alpha_{#1}}, ordering=Dynkin, x/.style={thin}]{E}{oooooo}$ & $\mfrak{p}_{\alpha_1}$, $\mfrak{p}_{\alpha_5}$ & $12$ & $-4$ & $32$ \\
  \hline
    $\mfrak{e}_7$ & $\dynkin[edge length=7mm, labels={1,2,3,4,5,6,7}, label macro/.code={\alpha_{#1}}, ordering=Dynkin, x/.style={thin}]{E}{ooooooo}$ & $\mfrak{p}_{\alpha_6}$ & $18$ & $-6$ & $54$ \\
  \hline
\end{tabular}
\caption{Collapsing levels of parabolic subalgebras with commutative nilradical}
\label{tab:parabolic commutative nilradical level}
\vspace{-2mm}
\end{table}

Therefore, by Theorem \ref{thm:FF homomorphism 1-graded} and Table \ref{tab:parabolic commutative nilradical level} we may construct a free field realization of the corresponding vertex algebra $\mcal{N}^\kappa_\mfrak{p}(\mfrak{g})$.


\subsection{Zhu's algebra modules}

Let $\mfrak{g}$ be a semisimple Lie algebra and let $\mfrak{p}$ be an opposite standard parabolic subalgebra of $\mfrak{g}$. Based on the previous section, we define a two-sided ideal $I_\mfrak{p}(\mfrak{g})$ of $U(\mfrak{g})$ by
\begin{align*}
  I_\mfrak{p}(\mfrak{g}) = \bigcap_{\lambda \in \smash{\mfrak{z}(\mfrak{l})^*}} \Ann_{U(\mfrak{g})}\! M^\mfrak{g}_\mfrak{p}(\lambda).
\end{align*}
We want to determine simple modules over the associative algebra $U(\mfrak{g})/I_\mfrak{p}(\mfrak{g})$, which is equivalent to describe simple $\mfrak{g}$-modules $E$ such that the ideal $I_\mfrak{p}(\mfrak{g})$ is contained in the annihilator $\Ann_{U(\mfrak{g})}\!E$. Let us assume that $E$ is a simple $\mfrak{g}$-module. Then by a theorem of Duflo \cite{Duflo1977} there exists $\lambda \in \mfrak{h}^*$ such that $\Ann_{U(\mfrak{g})}\!E =J_\lambda$, where
\begin{align*}
  J_\lambda = \Ann_{U(\mfrak{g})}\!L^\mfrak{g}_\mfrak{b}(\lambda)
\end{align*}
is a primitive ideal of $U(\mfrak{g})$. Therefore, we may restrict to the cases when $E$ is a simple highest weight $\mfrak{g}$-module. Let us denote by
\begin{align*}
  \Upsilon \colon Z(\mfrak{g}) \rarr S(\mfrak{h})^W
\end{align*}
the Harish-Chandra isomorphism. It is known that $\Upsilon$ is an isomorphism of algebras. Further, for $\lambda \in \mfrak{h}^*$ we define the algebra homomorphism $\chi_\lambda \colon Z(\mfrak{g}) \rarr \C$ by
\begin{align*}
  \chi_\lambda(z) = \Upsilon(z)(\lambda)
\end{align*}
for $z \in Z(\mfrak{g})$. Let us note that $\chi_{\lambda+\rho}$ is the central character of $M^\mfrak{g}_\mfrak{b}(\lambda)$ for $\lambda \in \mfrak{h}^*$. By definition of the ideal $I_\mfrak{p}(\mfrak{g})$, we obtain immediately that
\begin{align*}
  I_\mfrak{p}(\mfrak{g}) \cap Z(\mfrak{g}) = \bigcap_{\lambda \in \smash{\mfrak{z}(\mfrak{l})^*}} \ker \chi_{\lambda-\rho} = \bigcap_{\lambda \in \smash{\mfrak{z}(\mfrak{l})^*}} \ker \chi_{w_0(\lambda)+\rho},
\end{align*}
where $w_0$ is the longest element of $W$, which implies that $\Upsilon(z)(\lambda+\rho)=0$ for all $z \in I_\mfrak{p}(\mfrak{g}) \cap Z(\mfrak{g})$ provided $L^\mfrak{g}_\mfrak{b}(\lambda)$ with $\lambda \in \mfrak{h}^*$ satisfies $I_\mfrak{p}(\mfrak{g}) \subset J_\lambda$.
\medskip

For $\alpha \in \Pi$, let us consider the opposite standard parabolic subalgebra $\mfrak{p}$ of $\mfrak{g}$ associated to the subset $\Sigma = \Pi \setminus \{\alpha\}$ of $\Pi$. Hence, we have $\dim \mfrak{z}(\mfrak{l})=1$ and $\omega_\alpha \in \mfrak{z}(\mfrak{l})^*$, which implies that
\begin{align}
  I_\mfrak{p}(\mfrak{g}) \cap Z(\mfrak{g}) = \{z \in Z(\mfrak{g});\, \Upsilon(z)(a\omega_\alpha-\rho)=0\ \text{for all}\ a \in \C\}.  \label{eq:center 1-dim}
\end{align}
Further, by \cite{Borho-Jantzen1977} we have that
\begin{align*}
  I_\mfrak{p}(\mfrak{g}) = \bigcap_{\lambda \in \Lambda} \Ann_{U(\mfrak{g})}\! M^\mfrak{g}_\mfrak{p}(\lambda)
\end{align*}
for any Zariski dense subset $\Lambda$ of $\mfrak{z}(\mfrak{l})^*$. Therefore, we set
\begin{align*}
\Lambda = \{\lambda \in \mfrak{z}(\mfrak{l})^*;\, \langle \lambda, \gamma^\vee \rangle \notin \Z\ \text{for all $\gamma \in \Delta \setminus \Delta_\Sigma$}\},
\end{align*}
which is a Zariski dense subset of $\mfrak{z}(\mfrak{l})^*$. For $\lambda \in \Lambda$, we have that $M^\mfrak{g}_\mfrak{p}(\lambda)$ and $M^\mfrak{g}_{\smash{w_0(\mfrak{p})}}(w_0(\lambda))$ are simple $\mfrak{g}$-modules by \cite{Jantzen1977}. Besides, the annihilators $\Ann_{U(\mfrak{g})}\!M^\mfrak{g}_\mfrak{p}(\lambda)$ and $\Ann_{U(\mfrak{g})}\!M^\mfrak{g}_{\smash{w_0(\mfrak{p})}}(w_0(\lambda))$ are the unique maximal two-sided ideals of $U(\mfrak{g})$ containing $U(\mfrak{g})\ker \chi_{\lambda-\rho}$ and $U(\mfrak{g})\ker \chi_{w_0(\lambda)+\rho}$, respectively, by \cite[Section 8.5.8]{Dixmier1977-book}. However, as we have $\chi_{\lambda-\rho} = \chi_{w_0(\lambda)+\rho}$, we get that
\begin{align*}
  \Ann_{U(\mfrak{g})}\!M^\mfrak{g}_\mfrak{p}(\lambda) = \Ann_{U(\mfrak{g})}\!M^\mfrak{g}_{\smash{w_0(\mfrak{p})}}(w_0(\lambda))
\end{align*}
for $\lambda \in \Lambda$. Hence, we may write
\begin{align*}
  I_\mfrak{p}(\mfrak{g}) &= \bigcap_{\lambda \in \smash{\mfrak{z}(\mfrak{l})^*}} \Ann_{U(\mfrak{g})}\! M^\mfrak{g}_\mfrak{p}(\lambda) = \bigcap_{\lambda \in \Lambda} \Ann_{U(\mfrak{g})}\! M^\mfrak{g}_\mfrak{p}(\lambda) = \bigcap_{\lambda \in \Lambda} \Ann_{U(\mfrak{g})}\!M^\mfrak{g}_{\smash{w_0(\mfrak{p})}}(w_0(\lambda)) \\
  & = \bigcap_{\lambda \in \smash{\mfrak{z}(\mfrak{l})^*}} \Ann_{U(\mfrak{g})}\!M^\mfrak{g}_{\smash{w_0(\mfrak{p})}}(w_0(\lambda)).
\end{align*}
Let us note that $w_0(\mfrak{p})$ is the standard parabolic subalgebra of $\mfrak{g}$ associated to the subset $w_0(-\Sigma)$ of $\Pi$. As an immediate consequence we obtain that $L^\mfrak{g}_\mfrak{b}(w_0(\lambda))$ is a simple $U(\mfrak{g})/I_\mfrak{p}(\mfrak{g})$-module for all $\lambda \in \mfrak{z}(\mfrak{l})^*$.
\medskip

\theorem{\label{thm:modules over Zhu algebra}Let $\kappa$ be a $\mfrak{g}$-invariant symmetric bilinear form on $\mfrak{g}$ such that $(\kappa-\kappa_c^\mfrak{p})_{|[\mfrak{l},\mfrak{l}]}=0$. Then the $\mfrak{g}$-module $L^\mfrak{g}_\mfrak{b}(w_0(\lambda))$ is a simple $A(\mcal{N}^\kappa_\mfrak{p}(\mfrak{g}))$-module for $\lambda \in \mfrak{z}(\mfrak{l})^*$.}

\proof{Since we have $\pi_{\rm Zhu}(\mcal{I}^\kappa_\mfrak{p}(\mfrak{g})) \subset I_\mfrak{p}(\mfrak{g})$, it gives rise to a surjective homomorphism
\begin{align*}
   A(\mcal{N}^\kappa_\mfrak{p}(\mfrak{g})) \rarr U(\mfrak{g})/I_\mfrak{p}(\mfrak{g})
\end{align*}
of associative algebras. The rest of the proof follows from the discussion above.}

Let us recall that by Theorem \ref{thm:Zhu correspondence} we know that simple positive energy $\mcal{N}^\kappa_\mfrak{p}(\mfrak{g})$-modules are in a one-to-one correspondence with simple $A(\mcal{N}^\kappa_\mfrak{p}(\mfrak{g}))$-modules. Hence, Theorem \ref{thm:modules over Zhu algebra} provides a class of simple positive energy $\mcal{N}^\kappa_\mfrak{p}(\mfrak{g})$-modules. Besides, we have a conjecture that the Zhu's algebra $A(\mcal{N}^\kappa_\mfrak{p}(\mfrak{g}))$ is isomorphic to $U(\mfrak{g})/I_\mfrak{p}(\mfrak{g})$.


\subsection{Affine $\mcal{W}$-algebras and collapsing levels}
\label{subsection:affineW}

Let $f$ be a non-zero nilpotent element of a simple Lie algebra $\mfrak{g}$ and let $\kappa$ be a $\mfrak{g}$-invariant symmetric bilinear form on $\mfrak{g}$. By the Jacobson--Morozov theorem there exists an $\mfrak{sl}_2$-triple $(e,h,f)$ associated with $f$. Since $h$ is a semisimple element of $\mfrak{g}$, we have a ${1 \over 2}\Z$-grading on $\mfrak{g}$ defined by
\begin{align*}
  \mfrak{g}_j = \{u\in \mfrak{g};\, [h,u]=2ju\}
\end{align*}
for $j\in {1 \over 2}\Z$ and called the \emph{Dynkin grading} associated with $h$. We assume that $\mfrak{b} \subset \smash{\bigoplus_{j \geq 0} \mfrak{g}_j}$ and that $h \in \mfrak{h}$. We set also
\begin{align*}
  \Delta_j = \{\alpha \in \Delta;\, \mfrak{g}_\alpha \subset \mfrak{g}_j\}
\end{align*}
for $j \in {1 \over 2}\Z$. Further, we introduce nilpotent Lie subalgebras $\mfrak{g}_{\geq 1}$ and $\mfrak{g}_{>0}$ of $\mfrak{g}$ by
\begin{align*}
  \mfrak{g}_{\geq 1} = \bigoplus_{j \geq 1} \mfrak{g}_j \qquad \text{and} \qquad \mfrak{g}_{>0} = \bigoplus_{j > 0} \mfrak{g}_j.
\end{align*}
Let $(\cdot\,,\cdot)$ be a non-degenerate $\mfrak{g}$-invariant symmetric bilinear form on $\mfrak{g}$ normalized in such a way that $(e,f)=1$. We define a character $\chi \colon \mfrak{g}_{\geq 1} \rarr \C$ by
\begin{align*}
  \chi(u) = (f,u)
\end{align*}
for $u \in \mfrak{g}_{\geq 1}$. The character $\chi$ gives rise to a non-degenerate skew-symmetric bilinear form $\omega_\chi$ on $\smash{\mfrak{g}_{{1 \over 2}}}$ defined through
\begin{align*}
  \omega_\chi(u,v) = \chi([u,v])
\end{align*}
for $u,v \in \smash{\mfrak{g}_{{1\over 2}}}$. In addition, the bilinear form $\omega_\chi$ is $\mfrak{g}^\natural$-invariant, where $\mfrak{g}^\natural$ is the centralizer of the $\mfrak{sl}_2$-subalgebra $\mfrak{s} = \vspan_\C \{e,h,f\}$ in $\mfrak{g}$.

For a $\mcal{V}^\kappa(\mfrak{g})$-module $M$, we set
\begin{align*}
  C^\bullet(M) = M \otimes_\C \mcal{F}^{\omega_\chi}_{\mfrak{g}_{1/2}} \otimes_\C {\textstyle \bigwedge}^{{\infty \over 2}+\bullet}_{\mfrak{g}_{>0}}.
\end{align*}
Then it follows easily that $C^\bullet(\mcal{V}^\kappa(\mfrak{g}))$ has a natural structure of an $\N_0$-graded vertex superalgebra and $C^\bullet(M)$ is a $C^\bullet(\mcal{V}^\kappa(\mfrak{g}))$-module for any $\mcal{V}^\kappa(\mfrak{g})$-module $M$. In addition, there is an odd vector field $Q(z)$ in $C^\bullet(\mcal{V}^\kappa(\mfrak{g}))$ such that $\smash{Q_{(0)}^2}=0$, see \cite{Kac-Roan-Wakimoto2003}. Hence, we have that $(C^\bullet(M),Q_{(0)})$ is a chain complex for any $\mcal{V}^\kappa(\mfrak{g})$-module $M$. We define the \emph{generalized quantized Drinfeld--Sokolov reduction} $H^0_{DS,f}(M)$ associated with $(\mfrak{g},f)$ with coefficients in a $\mcal{V}^\kappa(\mfrak{g})$-module $M$ as
\begin{align*}
  H^0_{DS,f}(M) = H^0(C^\bullet(M),Q_{(0)}).
\end{align*}
The \emph{universal affine $\mcal{W}$-algebra} $\mcal{W}^\kappa(\mfrak{g},f)$ associated with $(\mfrak{g},f)$ of level $\kappa$ is then given by
\begin{align*}
  \mcal{W}^\kappa(\mfrak{g},f) = H_{DS,f}^0(\mcal{V}^\kappa(\mfrak{g})).
\end{align*}
Let us note that $\mcal{W}^\kappa(\mfrak{g},f)$ is naturally a ${1\over 2}\N_0$-graded vertex algebra. In fact, the grading of $\mcal{W}^\kappa(\mfrak{g},f)$ is induced by the ${1 \over 2}\Z$-grading of $C^\bullet(\mcal{V}^\kappa(\mfrak{g}))$ defined by
\begin{align}\label{eq:new-grading}
\begin{gathered}
  \deg e_{\alpha,n} = -n - {1\over 2}\alpha(h), \qquad \deg h_{\alpha,n} = -n, \qquad \deg f_{\alpha,n} = -n + {1 \over 2}\alpha(h) \qquad (\alpha \in \Delta_+), \\
  \deg \phi_{\alpha,n} = -n - {1\over 2}\alpha(h) \qquad (\alpha \in \Delta_{1\over 2}), \\
  \deg \psi_{\alpha,n} = -n - {1 \over 2}\alpha(h), \qquad \deg \psi_{\alpha,n}^* = -n + {1\over 2} \alpha(h) \qquad (\alpha \in \Delta_{>0})
\end{gathered}
\end{align}
for $n \in \Z$. Besides, we denote by $\mcal{W}_\kappa(\mfrak{g},f)$ the unique simple graded quotient of $\mcal{W}^\kappa(\mfrak{g},f)$.

In addition, we define a $\mfrak{g}^\natural$-invariant symmetric bilinear form $\kappa^\natural$ on $\mfrak{g}^\natural$ by
\begin{align}
  \kappa^\natural = \kappa + {\textstyle {1 \over 2}}(\kappa_\mfrak{g} - \kappa_{\mfrak{g}_0} - \kappa_{\mfrak{g}_{1/2}}), \label{eq:DS level}
\end{align}
where $\kappa_{\mfrak{g}_0}$ and $\kappa_{\mfrak{g}_{1/2}}$ are the trace forms of the $\mfrak{g}^\natural$-modules $\mfrak{g}_0$ and $\mfrak{g}_{1/2}$, respectively. Then there exists an injective homomorphism
\begin{align}
  \mcal{V}^{\kappa^\natural}\!(\mfrak{g}^\natural) \rarr \mcal{W}^\kappa(\mfrak{g},f) \label{eq:DS embedding}
\end{align}
of graded vertex algebras, see \cite{Kac-Roan-Wakimoto2003}, which leads to the following notion.
\medskip

\definition{Let $\kappa$ be a $\mfrak{g}$-invariant symmetric bilinear form on $\mfrak{g}$ and $f$ be a non-zero nilpotent element of $\mfrak{g}$. We say that the level $\kappa$ is \emph{collapsing} for $f$ if $\mcal{W}_\kappa(\mfrak{g},f) \simeq \mcal{V}_{\kappa^\natural}\!(\mfrak{g}^\natural)$. Equivalently, the level $\kappa$ is collapsing for $f$ if $\smash{\mcal{W}_\kappa(\mfrak{g},f)^{\mfrak{g}^\natural[[t]]}} \simeq \C$.}

The notion of collapsing levels for the case $f \in \mcal{O}_{\rm min}$ goes back to \cite{Adamovic-Kac-Frajria-Papi-Perse2018}. There is a full classification of collapsing levels for $f=f_{\rm min}$.
\medskip

Let $\kappa$ be a $\mfrak{g}$-invariant symmetric bilinear form on $\mfrak{g}$. Since $\mfrak{g}$ is a simple Lie algebra, we may write $\kappa=k\kappa_0$ for some $k \in \C$. Further, let us consider the minimal nilpotent element $f_\theta \in \mcal{O}_{\rm min}$, where $\theta$ is the maximal root of $\mfrak{g}$. Then $(e_\theta,h_\theta,f_\theta)$ is an $\mfrak{sl}_2$-triple associated with $f_\theta$ and the corresponding ${1\over 2}\Z$-grading on $\mfrak{g}$ is a \emph{minimal grading}, i.e.\ we have
\begin{align}
  \mfrak{g} = \mfrak{g}_{-1} \oplus \mfrak{g}_{-{1\over 2}} \oplus \mfrak{g}_0 \oplus \mfrak{g}_{{1\over 2}} \oplus \mfrak{g}_1
  \label{eq:minimal-Dynkin}
\end{align}
with $\mfrak{g}_{-1} = \C f_\theta$ and $\mfrak{g}_1 =\C e_\theta$. Let us note that the $\mfrak{g}$-invariant symmetric bilinear form $(\cdot\,,\cdot)$ on $\mfrak{g}$ is equal to $\kappa_0$. Moreover, we have
\begin{align*}
  \mfrak{g}^\natural = \{a \in \mfrak{g}_0;\, (h_\theta,a)=0\} \qquad \text{and} \qquad \mfrak{g}_0 = \mfrak{g}^\natural \oplus \C h_\theta.
\end{align*}
Therefore, if we set
\begin{align*}
  \mfrak{h}^\natural = \{h \in \mfrak{h};\, (h_\theta,h)=0\},
\end{align*}
then $\mfrak{h}^\natural$ is a Cartan subalgebra of $\mfrak{g}^\natural$ and satisfies
\begin{align*}
  \mfrak{h} = \mfrak{h}^\natural \oplus \C h_\theta.
\end{align*}
Besides, the root system of $\mfrak{g}^\natural$ with respect to $\mfrak{h}^\natural$ is naturally isomorphic to $\Delta_0$. We will use this identification in the next. As a positive root system in $\Delta_0$ we will consider $\Delta_+ \cap \Delta_0$ together with the set of simple roots $\Pi^\natural = \Pi \cap \Delta_0$.

The nilpotent Lie subalgebras $\mfrak{g}_{\geq 1}$ and $\mfrak{g}_{>0}$ of $\mfrak{g}$ are given by
\begin{align*}
  \mfrak{g}_{\geq 1} = \mfrak{g}_1 \qquad \text{and} \qquad \mfrak{g}_{>0} = \mfrak{g}_{{1\over 2}} \oplus \mfrak{g}_1.
\end{align*}
Let us recall also that the character $\chi \colon \mfrak{g}_{\geq 1} \rarr \C$ is defined by
\begin{align*}
  \chi(u) = (f_\theta,u)
\end{align*}
for $u \in \mfrak{g}_{\geq 1}$, which gives us $\chi(e_\theta)=1$. Moreover, the associated $\mfrak{g}^\natural$-invariant symplectic form $\omega_\chi$ on $\smash{\mfrak{g}_{{1\over 2}}}$ satisfies
\begin{align*}
  [u,v] = \omega_\chi(u,v)e_\theta
\end{align*}
for $u,v \in \smash{\mfrak{g}_{{1\over 2}}}$. The corresponding odd field $Q(z)$ in the ${1 \over 2}\Z$-graded vertex superalgebra $C^\bullet(\mcal{V}^\kappa(\mfrak{g}))$ is given by
\begin{multline*}
  Q(z) = \sum_{n \in \Z} Q_{(n)}z^{-n-1} = (e_\theta(z)+1)\psi_\theta^*(z) + \sum_{\alpha \in \Delta_{1/2}} (e_\alpha(z)+\phi_\alpha(z))\psi_\alpha^*(z) \\
  - {1 \over 2} \sum_{\alpha \in \Delta_{1/2}} \omega_\chi(e_\alpha,e_{\theta-\alpha}) \psi_\alpha^*(z)\psi_{\theta-\alpha}^*(z)\psi_\theta(z),
\end{multline*}
where the field $\phi_\alpha(z)$ in $\mcal{F}^{\omega_\chi}_{\mfrak{g}_{1/2}}$ corresponds to the root vector $e_\alpha$ in $\smash{\mfrak{g}_{1 \over 2}}$ for $\alpha \in \smash{\Delta_{1\over 2}}$, i.e.\ we have
\begin{align*}
  [\phi_\alpha(z),\phi_\beta(w)] = \omega_\chi(e_\alpha,e_\beta)\delta(z-w)
\end{align*}
for $\alpha, \beta \in \smash{\Delta_{1\over 2}}$. Since $Q_{(0)}^2=0$ and $(C^\bullet(M), Q_{(0)})$ is a cochain complex for any $\mcal{V}^\kappa(\mfrak{g})$-module $M$, we define
\begin{align*}
  H_{DS,f_\theta}^0(M) = H^0(C^\bullet(M),Q_{(0)}).
\end{align*}
The universal affine $\mcal{W}$-algebra $\mcal{W}^\kappa(\mfrak{g},f_\theta)$ associated with $(\mfrak{g},f_\theta)$ of level $\kappa$ is then given by
\begin{align*}
  \mcal{W}^\kappa(\mfrak{g},f_\theta) = H_{DS,f_\theta}^0(\mcal{V}^\kappa(\mfrak{g})).
\end{align*}
It is a conformal vertex algebra with central charge
\begin{align*}
  c_\kappa(\mfrak{g},f_\theta) = {k \dim \mfrak{g} \over k+h^\vee} -6k+h^\vee-4
\end{align*}
provided $k+h^\vee \neq 0$.
\medskip

\theorem{\label{thm:DS functor min}\cite{Arakawa2005, Arakawa2015}
\begin{enumerate}[topsep=3pt,itemsep=0pt]
  \item[i)] The functor $H_{DS,f_\theta}^0 \colon \mcal{O}_\kappa(\mfrak{g}) \rarr \mcal{W}^\kappa(\mfrak{g},f_\theta)\text{-}\Mod$ is exact.
  \item[ii)] Let $\mathbb{L}_{\kappa,\mfrak{g}}(L^\mfrak{g}_\mfrak{b}(\lambda))$ for $\lambda \in \mfrak{h}^*$ be a simple $\widehat{\mfrak{g}}_\kappa$-module. Then we have $H_{DS,f_\theta}^0(\mathbb{L}_{\kappa,\mfrak{g}}(L^\mfrak{g}_\mfrak{b}(\lambda))) \neq 0$ if and only if $k- \langle \lambda, \theta^\vee \rangle \notin \N_0$, where $\kappa = k\kappa_0$. Moreover, $\smash{H_{DS,f_\theta}^0}(\mathbb{L}_{\kappa,\mfrak{g}}(L^\mfrak{g}_\mfrak{b}(\lambda)))$ is a simple $\mcal{W}^\kappa(\mfrak{g},f_\theta)$-module if it is non-zero.
  \item[iii)] For any quotient $\mcal{V}$ of $\mcal{V}^\kappa(\mfrak{g})$, the associated variety of the vertex algebra $H^0_{DS,f_\theta}(\mcal{V})$, provided it is non-zero, is given by
  \begin{align*}
    X_{\smash{H^0_{DS,f_\theta}(\mcal{V})}} = X_\mcal{V} \cap S_{f_\theta}\subset \mfrak{g}^*,
  \end{align*}
  where $S_{f_\theta}$ is the Slodowy slice associated with $(e_\theta,h_\theta,f_\theta)$.
\end{enumerate}
\vspace{-1mm}}

Let us recall that the Slodowy slice $S_{f_\theta}$ associated with $(e_\theta,h_\theta,f_\theta)$ is defined through
\begin{align}
  S_{f_\theta} = \chi + \phi(\mfrak{g}^{e_\theta}), \label{eq:Slodowy slice def}
\end{align}
where $\phi \colon \mfrak{g} \rarr \mfrak{g}^*$ is an isomorphism of $\mfrak{g}$-modules induced from the non-degenerate $\mfrak{g}$-invariant symmetric bilinear form $(\cdot\,,\cdot)$ on $\mfrak{g}$.
\medskip

As an immediate consequence of the previous theorem we get that
\begin{align*}
  H^0_{DS,f_\theta}(\mcal{V}_\kappa(\mfrak{g})) \simeq \mcal{W}_\kappa(\mfrak{g},f_\theta)
\end{align*}
provided $k \notin \N_0$. In addition, if $k \notin \N_0$ and $M$ is a non-zero $\widehat{\mfrak{g}}_\kappa$-submodule of $\mcal{V}^\kappa(\mfrak{g})$, then we have that $\smash{H_{DS,f_\theta}^0(M)} \neq 0$.
\medskip

In the next, we will need a more detailed description of the level $\kappa^\natural$. Since $\mfrak{g}^\natural$ is a reductive Lie subalgebra of $\mfrak{g}_0$, we have the direct sum decomposition
\begin{align}
  \mfrak{g}^\natural = \bigoplus_{i=0}^r \mfrak{g}_i^\natural, \label{eq:natural decomposition}
\end{align}
of $\mfrak{g}^\natural$ into the direct sum of simple Lie subalgebras $\smash{\mfrak{g}^\natural_i}$ for $i=1,2,\dots,r$ and an abelian Lie subalgebra $\smash{\mfrak{g}_0^\natural}$ of $\mfrak{g}^\natural$ such that these direct summands are mutually orthogonal with respect to the Cartan--Killing form $\kappa_\mfrak{g}$ on $\mfrak{g}$. We denote by $\smash{\kappa_i^\natural}$ the normalized $\smash{\mfrak{g}_i^\natural}$-invariant symmetric bilinear form on $\smash{\mfrak{g}_i^\natural}$ for $i=1,2,\dots,r$ and $\smash{\kappa_0^\natural} = \kappa_{0|\smash{\mfrak{g}_0^\natural}}$.
\medskip

\lemma{\label{lem:form restriction}For $k\in \C$, we have
\begin{align*}
 \kappa^\natural{}_{|\smash{\mfrak{g}_i^\natural}} & = {(2k+h^\vee)\ind_D(\mfrak{g}_i^\natural,\mfrak{g}_{1\over 2}) + k\ind_D(\mfrak{g}_i^\natural,\mfrak{g}_i^\natural) \over 2h^\vee}\, \kappa_i^\natural, \\
  \kappa^\natural{}_{|\smash{\mfrak{g}_0^\natural}} &= (k+ {\textstyle {1\over 2}} h^\vee)\kappa_0^\natural
\end{align*}
for $i=1,2,\dots,r$.}

\proof{From the decomposition \eqref{eq:minimal-Dynkin} we obtain immediately that
\begin{align*}
  \kappa^\natural = \kappa + {\textstyle {1 \over 2}} \kappa_{\mfrak{g}_{1/2}} = {\textstyle {k \over 2h^\vee}} \kappa_\mfrak{g} + {\textstyle {1 \over 2}} \kappa_{\mfrak{g}_{1/2}}.
\end{align*}
Since $\kappa_\mfrak{g}$ and $\kappa_{\mfrak{g}_{1/2}}$ are trace forms of $\mfrak{g}^\natural_i$-modules $\mfrak{g}$ and $\mfrak{g}_{1\over 2}$ for $i=1,2,\dots,r$, we may write
\begin{align*}
  \kappa^\natural{}_{|\smash{\mfrak{g}_i^\natural}} &= {k \over 2h^\vee}\ind_D(\mfrak{g}_i^\natural,\mfrak{g}) \kappa_i^\natural + {1\over 2} \ind_D(\mfrak{g}_i^\natural,\mfrak{g}_{1\over 2}) \kappa_i^\natural \\
  &= {(2k+h^\vee)\ind_D(\mfrak{g}_i^\natural,\mfrak{g}_{1\over 2}) + k\ind_D(\mfrak{g}_i^\natural,\mfrak{g}_i^\natural) \over 2h^\vee}\, \kappa_i^\natural,
\end{align*}
where used that $\ind_D(\mfrak{g}_i^\natural,\mfrak{g}) = \ind_D(\mfrak{g}_i^\natural,\mfrak{g}_i^\natural) + 2 \ind_D(\mfrak{g}_i^\natural,\mfrak{g}_{1 \over 2})$ as follows from \eqref{eq:natural decomposition} and \eqref{eq:minimal-Dynkin}.
On the other hand, we have $\kappa_{\mfrak{g}}{}_{|\smash{\mfrak{g}_0^\natural}} = 2\kappa_{\mfrak{g}_{1/2}}{}_{|\smash{\mfrak{g}_0^\natural}}$, which gives us
\begin{align*}
  \kappa^\natural{}_{|\smash{\mfrak{g}_0^\natural}} = {k + {1 \over 2}h^\vee \over 2h^\vee}\, \kappa_{\mfrak{g}}{}_{|\smash{\mfrak{g}_0^\natural}} = (k + {\textstyle {1 \over 2}}h^\vee) \kappa_0{}_{|\smash{\mfrak{g}_0^\natural}} = (k + {\textstyle {1 \over 2}}h^\vee) \kappa_0^\natural,
\end{align*}
This finishes the proof.}

Based on the previous lemma, we know that
\begin{align*}
  \kappa^\natural{}_{|\smash{\mfrak{g}_i^\natural}} = k_i^\natural \kappa_i^\natural,
\end{align*}
where $k_i^\natural$ is a polynomial of degree $1$ in $k$, for $i=0,1,\dots,r$. The reductive Lie algebra $\mfrak{g}^\natural = \bigoplus_{i=0}^r \mfrak{g}_i^\natural$ and the polynomials $\smash{k_i^\natural}$ for $i=0,1,\dots,r$ are described in Table \ref{tab:g natural classical} and Table \ref{tab:g natural exceptional}, see \cite[Lemma 7.5]{Arakawa-Moreau2018}.

\begin{table}[ht]
\centering
\renewcommand{\arraystretch}{1.4}
\begin{tabular}{|c|c|c|c|c|c|c|}
  \hline
  $\mfrak{g}$ & $\mfrak{sl}_3$ & $\mfrak{sl}_n$, $n \geq 4$ & $\mfrak{sp}_{2n}$, $n \geq 2$ & $\mfrak{so}_7$ & $\mfrak{so}_8$ & $\mfrak{so}_n$, $n\geq 9$ \\
  \hline
  $\mfrak{g}^\natural$ & $\mfrak{g}_0^\natural \simeq \mfrak{gl}_1$ & $\begin{gathered} \mfrak{g}_0^\natural \simeq \mfrak{gl}_1 \\ \mfrak{g}_1^\natural \simeq \mfrak{sl}_{n-2} \end{gathered}$ & $\mfrak{g}_1^\natural \simeq \mfrak{sp}_{2n-2}$ & $\begin{gathered} \mfrak{g}_1^\natural \simeq \mfrak{sl}_2 \\ \mfrak{g}_2^\natural \simeq \mfrak{sl}_2 \end{gathered}$ & $\begin{gathered} \mfrak{g}_1^\natural \simeq \mfrak{sl}_2 \\ \mfrak{g}_2^\natural \simeq \mfrak{sl}_2 \\ \mfrak{g}_3^\natural \simeq \mfrak{sl}_2 \end{gathered}$ & $\begin{gathered} \mfrak{g}_1^\natural \simeq \mfrak{sl}_2 \\ \mfrak{g}_2^\natural \simeq \mfrak{so}_{n-4} \end{gathered}$ \\
  \hline
  $\kappa^\natural$ & $k_0^\natural = k + {3 \over 2}$ & $\begin{gathered} k_0^\natural = k + {\textstyle {n \over 2}} \\ k_1^\natural = k+1
  \end{gathered}$ & $k_1^\natural = k+{1\over 2}$ & $\begin{gathered} k_1^\natural = k + {\textstyle {3 \over 2}} \\ k_2^\natural = 2k+4 \end{gathered}$ & $\begin{gathered} k_1^\natural = k + 2 \\ k_2^\natural = k+2 \\ k_3^\natural = k+2
  \end{gathered}$ & $\begin{gathered} k_1^\natural = k + {\textstyle {n \over 2}} - 2 \\ k_2^\natural = k+2
  \end{gathered}$ \\
  \hline
\end{tabular}
\caption{$\mfrak{g}^\natural$ and $\kappa^\natural$ for the classical types}
\label{tab:g natural classical}
\vspace{-2mm}
\end{table}

\begin{table}[ht]
\centering
\renewcommand{\arraystretch}{1.4}
\begin{tabular}{|c|c|c|c|c|c|}
  \hline
  $\mfrak{g}$ & $\mfrak{g}_2$ & $\mfrak{f}_4$ & $\mfrak{e}_6$ & $\mfrak{e}_7$ & $\mfrak{e}_8$ \\
  \hline
  $\mfrak{g}^\natural$ & $\mfrak{g}_1^\natural \simeq \mfrak{sl}_2$ & $\mfrak{g}_1^\natural \simeq \mfrak{sp}_6$ & $\mfrak{g}_1^\natural \simeq \mfrak{sl}_6$ & $\mfrak{g}_1^\natural \simeq \mfrak{so}_{12}$ & $\mfrak{g}_1^\natural \simeq \mfrak{e}_7$  \\
  \hline
  $\kappa^\natural$ & $k_1^\natural = 3k+5$ & $k_1^\natural = k+{5 \over 2}$ & $k_1^\natural = k+3$ & $k_1^\natural = k+4$ & $k_1^\natural = k+6$ \\
  \hline
\end{tabular}
\caption{$\mfrak{g}^\natural$ and $\kappa^\natural$ for the exceptional types}
\label{tab:g natural exceptional}
\vspace{-2mm}
\end{table}

Let us recall that we have the injective homomorphism \eqref{eq:DS embedding} of graded vertex algebras, which gives rise to a surjective morphism
\begin{align*}
  X_{\smash{\mcal{W}^\kappa(\mfrak{g},f_\theta)}} \rarr X_{\smash{\mcal{V}^{\kappa^\natural}\!(\mfrak{g}^\natural)}}
\end{align*}
of the corresponding associated varieties. Since $X_{\smash{\mcal{W}^\kappa(\mfrak{g},f_\theta)}} \simeq S_{f_\theta}$ and $X_{\smash{\mcal{V}^{\kappa^\natural}\!(\mfrak{g}^\natural)}} \simeq (\mfrak{g}^\natural)^*$, we have a surjective morphism $S_{f_\theta} \rarr (\mfrak{g}^\natural)^*$
of Poisson varieties. By using the $\mfrak{g}$-invariant symmetric bilinear form $(\cdot\,,\cdot)$ on $\mfrak{g}$, we identify subsets of $\mfrak{g}^*$ with subsets of $\mfrak{g}$. Hence, we get a surjective morphism
\begin{align}
 \varphi_{f_\theta} \colon f_\theta + \mfrak{g}^{e_\theta} \rarr \mfrak{g}^\natural \label{eq:Slodowy slice}
\end{align}
of varieties given by
\begin{align*}
  f_\theta + a \mapsto a^\natural
\end{align*}
for $a \in \mfrak{g}^{e_\theta}$, where $a^\natural$ is the $\mfrak{g}^\natural$-part of $a$ with respect to the decomposition $\mfrak{g}^{e_\theta} = \mfrak{g}^\natural \oplus \smash{\mfrak{g}_{1\over 2}} \oplus \mfrak{g}_1$.


\subsection{Associated varieties of simple affine vertex algebras}

Let $\mfrak{g}$ be a simple Lie algebra and let us fix a minimal grading of $\mfrak{g}$ given by the element $f_\theta \in \mcal{O}_{\rm min}$ as in Section \ref{subsection:affineW}. We assume that $\mfrak{b} \subset \smash{\bigoplus_{j\geq 0}\mfrak{g}_j}$. Let us consider an opposite standard parabolic subalgebra $\mfrak{p}=\mfrak{l} \oplus \widebar{\mfrak{u}}$ of $\mfrak{g}$ with the Levi subalgebra $\mfrak{l}$, the nilradical $\widebar{\mfrak{u}}$ and the opposite nilradical $\mfrak{u}$. Let us assume that $\mfrak{u}$ is a commutative Lie algebra and $[\mfrak{l},\mfrak{l}] \neq 0$. Further, let $\kappa_\mfrak{p}$ be a $\mfrak{g}$-invariant symmetric bilinear form on $\mfrak{g}$ satisfying the additional condition  $(\kappa_\mfrak{p}-\kappa_c^\mfrak{p})_{|[\mfrak{l},\mfrak{l}]} = 0$. So we are in the setting of Table \ref{tab:parabolic commutative nilradical level}. Let us note that such $\kappa_\mfrak{p}$ is uniquely determined by $\mfrak{p}$ when it exists. Then, as follows from Section \ref{subsec:FF homomorphism}, there is the homomorphism
\begin{align*}
  \widetilde{w}_{\kappa_\mfrak{p},\mfrak{g}}^\mfrak{p} \colon \mcal{V}^{\kappa_\mfrak{p}}\!(\mfrak{g}) \rarr \mcal{M}_\mfrak{u} \otimes_\C \mcal{V}^{\kappa_\mfrak{p}+{1\over 2} \kappa_\mfrak{g}}\!(\mfrak{z}(\mfrak{l}))
\end{align*}
of $\N_0$-graded vertex algebras. The image of this homomorphism, which is naturally an $\N_0$-graded vertex algebra, we denote by $\smash{\widetilde{\mcal{V}}}^{\kappa_\mfrak{p}}\!(\mfrak{g})$.

We shall consider the universal affine $\mcal{W}$-algebra $\mcal{W}^{\kappa_\mfrak{p}}\!(\mfrak{g},f_\theta)$ of level $\kappa_\mfrak{p}$ associated to the minimal nilpotent element $f_\theta$. Since $\mfrak{p}$ is an opposite standard parabolic subalgebra of $\mfrak{g}$, we may introduce an opposite standard parabolic subalgebra $\mfrak{p}^\natural$ of $\mfrak{g}^\natural$ by
\begin{align*}
  \mfrak{p}^\natural = \mfrak{p} \cap \mfrak{g}^\natural.
\end{align*}
In addition, the parabolic subalgebra $\mfrak{p}^\natural$ determines a triangular decomposition
\begin{align*}
  \mfrak{g}^\natural = \widebar{\mfrak{u}}^\natural \oplus \mfrak{l}^\natural \oplus \mfrak{u}^\natural
\end{align*}
of $\mfrak{g}^\natural$, where
\begin{align*}
   \widebar{\mfrak{u}}^\natural = \widebar{\mfrak{u}} \cap \mfrak{g}^\natural, \qquad  \mfrak{l}^\natural = \mfrak{l} \cap \mfrak{g}^\natural, \qquad \mfrak{u}^\natural = \mfrak{u} \cap \mfrak{g}^\natural.
\end{align*}
Obviously, the Lie algebra $\mfrak{u}^\natural$ is commutative.

We denote by $\Sigma$ and $\Sigma^\natural$ the sets of simple roots of $\mfrak{l}$ and $\mfrak{l}^\natural$, respectively. Therefore, we obtain that $\Sigma = \Pi \cap \Delta_+^\mfrak{l}$ and $\Sigma^\natural = \Sigma \cap \Pi^\natural$. Further, let us recall that $\Pi \setminus \Sigma = \{\alpha\}$ for some simple root $\alpha \in \Pi$, which implies that $\Pi^\natural \setminus \Sigma^\natural = \Pi^\natural \cap \{\alpha\}$.
\medskip

\lemma{\label{lem:for_induction} There exists a Langrangian subspace $\mfrak{l}_\chi$ of $\smash{\mfrak{g}_{1\over 2}}$ such that $\mfrak{u}=\mfrak{u}^\natural \oplus \mfrak{l}_\chi \oplus \mfrak{g}_1$.}

\proof{As we have $\mfrak{u} \subset \mfrak{g}^\natural \oplus \smash{\mfrak{g}_{1\over 2}} \oplus \mfrak{g}_1$, we may write
\begin{align*}
  \mfrak{u} = (\mfrak{u} \cap \mfrak{g}^\natural) \oplus (\mfrak{u} \cap \smash{\mfrak{g}_{1 \over 2}}) \oplus (\mfrak{u} \cap \mfrak{g}_1) = \mfrak{u}^\natural \oplus \mfrak{l}_\chi \oplus \mfrak{g}_1,
\end{align*}
where $\mfrak{l}_\chi = \mfrak{u} \cap \smash{\mfrak{g}_{1 \over 2}}$. Since $\mfrak{u}$ is a commutative Lie algebra, we obtain that $\mfrak{l}_\chi$ is an isotropic subspace of $\smash{\mfrak{g}_{1 \over 2}}$. Let us suppose that $\mfrak{l}_\chi$ is not a Lagrangian subspace of $\smash{\mfrak{g}_{1 \over 2}}$. Then there exists $\alpha \in \smash{\Delta_{1 \over 2}}$ such that $\alpha, \theta-\alpha \in \Delta_+^\mfrak{l}$, which is a contradiction with the fact that $\theta \in \Delta_+^\mfrak{u}$.}

We denote by $\kappa_{\mfrak{p}^\natural}$ the $\mfrak{g}^\natural$-invariant symmetric bilinear form on $\mfrak{g}^\natural$ defined through the formula \eqref{eq:DS level}.
\medskip

\lemma{\label{lem:level condition} We have $(\kappa_{\mfrak{p}^\natural}-\kappa_c^{\mfrak{p}^\natural})_{|\mfrak{l}^\natural} = (\kappa_\mfrak{p} - \kappa_c^\mfrak{p})_{|\mfrak{l}^\natural}$.}

\proof{As $\mfrak{u}$, $\mfrak{l}$ and $\smash{\mfrak{g}_{1 \over 2}}$ are $\mfrak{l}^\natural$-modules, we have that $\smash{\mfrak{g}_{1 \over 2}} = (\mfrak{u} \cap \smash{\mfrak{g}_{1 \over 2}}) \oplus (\mfrak{l} \cap \smash{\mfrak{g}_{1 \over 2}})$ as $\mfrak{l}^\natural$-modules. Besides, the symplectic form $\omega_\chi$ on $\smash{\mfrak{g}_{1 \over 2}}$ is $\mfrak{l}^\natural$-invariant, which implies that $\smash{\kappa_{\mfrak{g}_{1/2}}{}_{|\mfrak{l}^\natural}} = 2\kappa_{\mfrak{l}_\chi}$. Further, from the decomposition \eqref{eq:minimal-Dynkin} and the definition of $\kappa_{\mfrak{p}^\natural}$, we have that $\kappa_{\mfrak{p}^\natural} = \kappa_\mfrak{p} + \smash{{1 \over 2}\kappa_{\mfrak{g}_{1/2}}}$. Hence, we may write
\begin{align*}
  (\kappa_{\mfrak{p}^\natural}-\kappa_c^{\mfrak{p}^\natural})_{|\mfrak{l}^\natural} = \kappa_\mfrak{p}{}_{|\mfrak{l}^\natural} + \kappa_{\mfrak{l}_\chi} + \kappa_{\mfrak{u}^\natural} = \kappa_\mfrak{p}{}_{|\mfrak{l}^\natural} + \kappa_\mfrak{u},
\end{align*}
where we used that $\mfrak{u} = \mfrak{u}^\natural \oplus \mfrak{l}_\chi \oplus \mfrak{g}_1$. On the other hand, we have $(\kappa_\mfrak{p} - \kappa_c^\mfrak{p})_{|\mfrak{l}} = \kappa_\mfrak{p}{}_{|\mfrak{l}} + \kappa_\mfrak{u}$, which gives us
\begin{align*}
  (\kappa_{\mfrak{p}^\natural}-\kappa_c^{\mfrak{p}^\natural})_{|\mfrak{l}^\natural} = (\kappa_\mfrak{p} - \kappa_c^\mfrak{p})_{|\mfrak{l}^\natural}.
\end{align*}
This finishes the proof.}

For an element $h \in \mfrak{z}(\mfrak{l}^\natural)$, we denote by $h_{\mfrak{z}(\mfrak{l})}$ the $\mfrak{z}(\mfrak{l})$-part of $h$ with respect to the decomposition $\mfrak{l} = \mfrak{z}(\mfrak{l}) \oplus [\mfrak{l},\mfrak{l}]$.
\medskip

\lemma{\label{lem:vertex algebra embedding}There is an isomorphism
\begin{align*}
  \iota \colon \mcal{V}_{\kappa_{\smash{\mfrak{p}^\natural}}+ {1\over 2}\kappa_{\smash{\mfrak{g}^\natural}}}\!(\mfrak{z}(\mfrak{l}^\natural)) \rarr \mcal{V}^{\kappa_\mfrak{p}+{1\over 2} \kappa_\mfrak{g}}\!(\mfrak{z}(\mfrak{l}))
\end{align*}
of $\N_0$-graded vertex algebras given by
\begin{align*}
  h(z) \mapsto h_{\mfrak{z}(\mfrak{l})}(z)
\end{align*}
for $h \in \mfrak{z}(\mfrak{l}^\natural)$.}

\proof{For elements $a, b \in \mfrak{z}(\mfrak{l}^\natural)$, we may write
\begin{align*}
  (\kappa_{\smash{\mfrak{p}^\natural}}+ {\textstyle {1\over 2}}\kappa_{\smash{\mfrak{g}^\natural}})(a,b) &= (\kappa_{\mfrak{p}^\natural}-\kappa_c^{\mfrak{p}^\natural})(a,b) = (\kappa_\mfrak{p} - \kappa_c^\mfrak{p})(a,b) = (\kappa_\mfrak{p} - \kappa_c^\mfrak{p})(a_{\mfrak{z}(\mfrak{l})},b_{\mfrak{z}(\mfrak{l})}) \\
  &= (\kappa_\mfrak{p}+{\textstyle {1\over 2}} \kappa_\mfrak{g})(a_{\mfrak{z}(\mfrak{l})},b_{\mfrak{z}(\mfrak{l})}),
\end{align*}
where we used Lemma \ref{lem:level condition} in the second equality, which implies that the mapping $\iota$ is a homomorphism of $\N_0$-graded vertex algebras. To show that $\iota$ is an isomorphism, it is enough to find an element $a \in \mfrak{z}(\mfrak{l}^\natural)$ such that $a_{\mfrak{z}(\mfrak{l})} \neq 0$.

Let us consider the grading element $h \in \mfrak{z}(\mfrak{l})$ for the decomposition $\mfrak{g} = \widebar{\mfrak{u}} \oplus \mfrak{l} \oplus \mfrak{u}$, which means that $[h,a]=a$ for $a \in \mfrak{u}$ and $[h,a]=-a$ for $a \in \widebar{\mfrak{u}}$. Then we have
\begin{align*}
   (h_\theta,h) = {\kappa_\mfrak{g}(h_\theta,h) \over 2h^\vee} = {4+ 2\dim \mfrak{l}_\chi \over 2h^\vee} = {4 + \dim \mfrak{g}_{1 \over 2} \over 2h^\vee}=1,
\end{align*}
where we used the fact that $\dim \smash{\mfrak{g}_{1\over 2}} = 2h^\vee -4$ following from \cite[Lemma 3]{Wang1999}. Therefore, we get a decomposition
\begin{align*}
  h = {\textstyle {1 \over 2}} h_\theta + h^\natural,
\end{align*}
where $h^\natural \in \mfrak{z}(\mfrak{l}^\natural)$. Moreover, we have that $h^\natural$ is a non-zero element. We show that $\smash{(h^\natural)_{\mfrak{z}(\mfrak{l})}} \neq 0$. We may write
\begin{align*}
  (h^\natural)_{\mfrak{z}(\mfrak{l})} = {\kappa_\mfrak{g}(h,h^\natural) \over \kappa_\mfrak{g}(h,h)}\, h = {\kappa_\mfrak{g}(h,h) - {1 \over 2}\kappa_\mfrak{g}(h,h_\theta) \over \kappa_\mfrak{g}(h,h)}\, h = \bigg(\!1 - {h^\vee \over \kappa_\mfrak{g}(h,h)}\!\bigg) h = \bigg(\!1 - {h^\vee \over 2 \dim \mfrak{u}}\!\bigg) h,
\end{align*}
where we used that $\kappa_\mfrak{g}(h,h) = 2\dim \mfrak{u}$. Further, from Table \ref{tab:parabolic commutative nilradical level} we have that $h^\vee \neq 2\dim \mfrak{u}$, which implies that $\smash{(h^\natural)_{\mfrak{z}(\mfrak{l})}} \neq 0$.}

By using Lemma \ref{lem:level condition} and the fact that $\smash{(\kappa_\mfrak{p}-\kappa_c^\mfrak{p})_{|[\mfrak{l},\mfrak{l}]}} = 0$, we obtain that $\smash{(\kappa_{\mfrak{p}^\natural}-\kappa_c^{\mfrak{p}^\natural})_{|[\mfrak{l}^\natural, \mfrak{l}^\natural]}} = 0$. Moreover, we have that $\mfrak{u}^\natural$ is a commutative Lie algebra. Therefore, as follows from Section \ref{subsec:FF homomorphism}, there is the homomorphism
\begin{align*}
  \widehat{w}^{\mfrak{p}^\natural}_{\smash{\kappa_{\smash{\mfrak{p}^\natural}},\mfrak{g}^\natural}} \colon \mcal{V}^{\kappa_{\smash{\mfrak{p}^\natural}}}\!(\mfrak{g}^\natural) \rarr \mcal{M}_{\mfrak{u}^\natural}\! \otimes_\C \mcal{V}_{\kappa_{\smash{\mfrak{p}^\natural}}+ {1\over 2}\kappa_{\smash{\mfrak{g}^\natural}}}\!(\mfrak{z}(\mfrak{l}^\natural))
\end{align*}
of $\N_0$-graded vertex algebras. The image of this homomorphism, which is naturally an $\N_0$-graded vertex algebra, we denote by $\smash{\widehat{\mcal{V}}}^{\kappa_{\smash{\mfrak{p}^\natural}}}\! (\mfrak{g}^\natural)$.
By applying the functor $\smash{H_{DS,f_\theta}^0}(?)$ to $\widetilde{w}^\mfrak{p}_{\kappa_\mfrak{p},\mfrak{g}}$ and by using the fact that there is an embedding of $\smash{\mcal{V}^{\kappa_{\smash{\mfrak{p}^\natural}}}}\!(\mfrak{g}^\natural)$ into $\mcal{W}^{\kappa_\mfrak{p}}\!(\mfrak{g},f_\theta)$, we get a sequence
\begin{align*}
  \mcal{V}^{\kappa_{\smash{\mfrak{p}^\natural}}}\!(\mfrak{g}^\natural) \xrarr{\varphi} \mcal{W}^{\kappa_\mfrak{p}}\!(\mfrak{g}, f_\theta) \xrarr{\psi} H_{DS,f_\theta}^0(\mcal{M}_{\mfrak{u}} \otimes_\C \mcal{V}^{\kappa_\mfrak{p}+{1\over 2}\kappa_\mfrak{g}}\!(\mfrak{z}(\mfrak{l})))
\end{align*}
of graded vertex algebras. Let us note that the homomorphism $\varphi$ is given by
\begin{align}
  \varphi(a(z)) = a(z) + \sum_{\alpha,\beta \in \Delta_{>0}} [a,e_\beta]_\alpha \normOrd{\psi_\alpha(z)\psi^*_\beta(z)} - {1 \over 2} \sum_{\alpha,\beta \in \Delta_{1/2}} [a,e_\alpha]_\beta \normOrd{\phi^\alpha(z)\phi_\beta(z)} \label{eq:natural embedding}
\end{align}
for $a \in \mfrak{g}^\natural$, where the fields $\{\phi_\alpha(z);\, \alpha \in \smash{\Delta_{1 \over 2}}\}$ in $\mcal{F}^{\omega_\chi}_{\mfrak{g}_{1/2}}$ correspond to the basis $\{e_\alpha;\, \alpha \in \smash{\Delta_{1 \over 2}}\}$ of $\smash{\mfrak{g}_{1 \over 2}}$ and the fields $\{\phi^\alpha(z);\, \alpha \in \smash{\Delta_{1 \over 2}}\}$ in $\smash{\mcal{F}^{\omega_\chi}_{\mfrak{g}_{1/2}}}$ are determined by the dual basis $\{e^\alpha;\, \alpha \in \smash{\Delta_{1 \over 2}}\}$ of $\smash{\mfrak{g}_{1 \over 2}}$ with respect to $\omega_\chi$, see \cite{Kac-Wakimoto2004}. For $a \in \mfrak{g}^\natural$, we set
\begin{align*}
  \varphi_0(a(z)) = \varphi(a(z)) - a(z).
\end{align*}
In fact, we have the following important assertion.
\medskip

\theorem{\label{thm:DS reduction}We have
\begin{enumerate}[topsep=3pt,itemsep=0pt]
  \item[i)] $H_{DS,f_\theta}^0(\mcal{M}_{\mfrak{u}} \otimes_\C   \mcal{V}^{\kappa_\mfrak{p}+{1\over 2}\kappa_\mfrak{g}}\!(\mfrak{z}(\mfrak{l}))) \simeq \mcal{M}_{\mfrak{u}^\natural}\! \otimes_\C   \smash{\mcal{V}_{\kappa_{\smash{\mfrak{p}^\natural}}+ {1\over 2}\kappa_{\smash{\mfrak{g}^\natural}}}\!(\mfrak{z}(\mfrak{l}^\natural))}$,
  \item[ii)] $\im(\psi) = \im(\psi \circ \varphi)$,
  \item[iii)] $H^0_{DS,f_\theta}(\smash{\widetilde{\mcal{V}}}^{\kappa_\mfrak{p}}\!(\mfrak{g})) \simeq \smash{\widehat{\mcal{V}}}^{\kappa_{\smash{\mfrak{p}^\natural}}}\! (\mfrak{g}^\natural)$.
\end{enumerate}}

\proof{i) We set
\begin{align*}
  J_a(z) = - \sum_{\alpha \in \smash{\Delta_+^{\mfrak{u}^\natural}}} [a]_\alpha a_\alpha(z) + \sum_{\alpha,\beta \in \Delta_{>0}} [a,e_\beta]_\alpha \normOrd{\psi_\alpha(z)\psi^*_\beta(z)} - {1 \over 2} \sum_{\alpha,\beta \in \Delta_{1/2}} [a,e_\alpha]_\beta \normOrd{\phi^\alpha(z)\phi_\beta(z)}
\end{align*}
for $a \in \mfrak{u}^\natural$. Based on formula \eqref{eq:natural embedding} and Lemma \ref{lem:vertex algebra embedding}, there is a vertex algebra homomorphism
\begin{align*}
   \Psi\colon \mcal{M}_{\mfrak{u}^\natural}\! \otimes_\C   \mcal{V}_{\kappa_{\smash{\mfrak{p}^\natural}} + {1\over 2} \kappa_{\smash{\mfrak{g}^\natural}}}\!(\mfrak{z}(\mfrak{l}^\natural))
   \rarr
   H_{DS,f_\theta}^0(\mcal{M}_{\mfrak{u}} \otimes_\C   \mcal{V}^{\kappa_\mfrak{p}+{1\over 2}\kappa_\mfrak{g}}\!(\mfrak{z}(\mfrak{l})))
\end{align*}
given by
\begin{align*}
  a_\alpha(z) \mapsto -J_{e_\alpha}\!(z), \qquad
  a_\alpha^*(z) \mapsto a_\alpha^*(z), \qquad
  h(z) \mapsto h_{\mfrak{z}(\mfrak{l})}(z)
\end{align*}
for $\alpha \in \Delta^{\mfrak{u}^\natural}_+$ and $h \in \mfrak{z}(\mfrak{l}^\natural)$. Since $Q_{(0)}(\Psi(a_\alpha(z))) = 0$, $Q_{(0)}(\Psi(a^*_\alpha(z))) = 0$ and $Q_{(0)}(\Psi(h(z))) = 0$, the homomorphism $\Psi$ is well defined.

We wish to show  that the above mapping is an isomorphism. In order to do so, let us consider the Hochschild--Serre spectral sequence
\begin{align*}
E_r^{p,q} \Rightarrow H_{DS,f_\theta}^{p+q}(\mcal{M}_\mfrak{u} \otimes_\C  \mcal{V}^{\kappa_\mfrak{p}+{1\over 2}\kappa_\mfrak{g}}\!(\mfrak{z}(\mfrak{l}))) = H^{{\infty\over 2}+p+q}(\mfrak{g}_{>0}(\!(t)\!), \mcal{M}_\mfrak{u} \otimes_\C   \mcal{V}^{\kappa_\mfrak{p}+{1\over 2}\kappa_\mfrak{g}}\!(\mfrak{z}(\mfrak{l}))
\otimes_\C \mcal{F}^{\omega_\chi}_{\mfrak{g}_{1/2}})
\end{align*}
for the ideal $(\mfrak{l}_\chi \oplus \mfrak{g}_{\geq 1})(\!(t)\!)$ of $\mfrak{g}_{>0}(\!(t)\!)$. By definition, we have
\begin{align*}
E_1^{\bullet,q} = H^{{\infty \over 2}+q}((\mfrak{l}_\chi \oplus \mfrak{g}_{\geq 1})(\!(t)\!), \mcal{M}_\mfrak{u} \otimes_\C \mcal{V}^{\kappa_\mfrak{p}+ {1\over 2} \kappa_\mfrak{g}}\!(\mfrak{z}(\mfrak{l})) \otimes_\C \mcal{F}^{\omega_\chi}_{\mfrak{g}_{1/2}} \otimes_\C {\textstyle \bigwedge}^{{\infty \over 2}+\bullet}_{\mfrak{g}_{>0}/(\mfrak{l}_\chi \oplus \mfrak{g}_{\geq 1})}).
\end{align*}
Further, by Lemma \ref{lem:for_induction} we get that
\begin{align*}
E_1^{\bullet,q} \simeq \delta_{q,0}\, \mcal{M}_{\mfrak{u}^\natural} \otimes_\C  \mcal{V}^{\kappa_\mfrak{p}+{1\over 2}\kappa_\mfrak{g}}\!(\mfrak{z}(\mfrak{l}))\otimes  \mcal{F}^{\omega_\chi}_{\mfrak{g}_{1/2}} \otimes_\C {\textstyle \bigwedge}^{{\infty \over 2}+\bullet}_{\mfrak{g}_{>0}/(\mfrak{l}_\chi \oplus \mfrak{g}_{\geq 1})}.
\end{align*}
Since $\smash{\mcal{F}^{\omega_\chi}_{\mfrak{g}_{1/2}}}$ is free over $\mfrak{g}_{>0}/(\mfrak{l}_\chi \oplus \mfrak{g}_{\geq 1})\otimes_\C t^{-1}\C[t^{-1}]$ and cofree over $\mfrak{g}_{>0}/(\mfrak{l}_\chi \oplus \mfrak{g}_{\geq 1})\otimes_\C\C[[t]]$, we obtain that
\begin{align*}
E_2^{p,q} &=\delta_{q,0}\,H^{{\infty \over 2}+p}(\mfrak{g}_{>0}/(\mfrak{l}_\chi \oplus \mfrak{g}_{\geq 1})(\!(t)\!), \mcal{M}_{\mfrak{u}^\natural} \otimes_\C \mcal{V}^{\kappa_\mfrak{p}+{1\over 2}\kappa_\mfrak{g}}\!(\mfrak{z}(\mfrak{l})) \otimes_\C  \mcal{F}^{\omega_\chi}_{\mfrak{g}_{1/2}}) \\
&\simeq \delta_{p,0}\delta_{q,0}\, \mcal{M}_{\mfrak{u}^\natural} \otimes_\C  \mcal{V}^{\kappa_\mfrak{p}+{1\over 2}\kappa_\mfrak{g}}\!(\mfrak{z}(\mfrak{l})) \simeq \delta_{p,0}\delta_{q,0}\, \mcal{M}_{\mfrak{u}^\natural} \otimes_\C  \mcal{V}_{\kappa_{\smash{\mfrak{p}^\natural}}+ {1\over 2}\kappa_{\smash{\mfrak{g}^\natural}}}\!(\mfrak{z}(\mfrak{l}^\natural))
\end{align*}
as vector spaces, where we used Lemma \ref{lem:vertex algebra embedding}. It follows that $\Psi$ is an isomorphism. In addition, since we have $\mfrak{g}^\natural \subset \mfrak{g}_0$, the vertex algebra $\smash{H_{DS,f_\theta}^0(\mcal{M}_{\mfrak{u}} \otimes_\C   \mcal{V}^{\kappa_\mfrak{p}+{1\over 2}\kappa_\mfrak{g}}\!(\mfrak{z}(\mfrak{l})))}$ is $\N_0$-graded, see \eqref{eq:new-grading}. Hence, it implies that $\Psi$ is an isomorphism of $\N_0$-graded vertex algebras.
\smallskip

\noindent
ii) As the vertex algebra $\smash{H_{DS,f_\theta}^0(\mcal{M}_{\mfrak{u}} \otimes_\C   \mcal{V}^{\kappa_\mfrak{p}+{1\over 2}\kappa_\mfrak{g}}\!(\mfrak{z}(\mfrak{l})))}$ is $\N_0$-graded by (i), there are no elements of weight ${3 \over 2}$ in the image of $\psi$. Since $\mcal{W}^{\kappa_\mfrak{p}}\!(\mfrak{g},f_\theta)$ is generated by the elements of $\mfrak{g}^\natural$ of weight $1$ and by the elements of $\smash{\mfrak{g}_{1 \over 2}}$ of weight ${3 \over 2}$, see \cite{Kac-Wakimoto2004}, it follows that
the image of $\mcal{W}^{\kappa_\mfrak{p}}\!(\mfrak{g},f_\theta)$ in $\smash{H_{DS,f_\theta}^0(\mcal{M}_{\mfrak{u}} \otimes_\C \smash{\mcal{V}^{\kappa_\mfrak{p}+{1\over 2}\kappa_\mfrak{g}}\!(\mfrak{z}(\mfrak{l}))})}$ is generated by the image of $\mcal{V}^{\kappa_{\smash{\mfrak{p}^\natural}}}\!(\mfrak{g}^\natural)$.
\smallskip

\noindent
iii) We claim that the following diagram
\begin{align*}
\bfig
\Vtriangle(0,0)|alr|<650,500>[\mcal{V}^{\kappa_{\smash{\mfrak{p}^\natural}}}\!(\mfrak{g}^\natural)`
H_{DS,f_\theta}^0(\mcal{M}_{\mfrak{u}} \otimes_\C \mcal{V}^{\kappa_\mfrak{p}+{1\over 2}\kappa_\mfrak{g}}\!(\mfrak{z}(\mfrak{l})))`
\mcal{M}_{\mfrak{u}^\natural}\! \otimes_\C \mcal{V}_{\kappa_{\smash{\mfrak{p}^\natural}}+ {1\over 2}\kappa_{\smash{\mfrak{g}^\natural}}}\!(\mfrak{z}(\mfrak{l}^\natural))
;\psi\circ \varphi` \smash{\widehat{w}^{\mfrak{p}^\natural}_{\smash{\kappa_{\smash{\mfrak{p}^\natural}},\mfrak{g}^\natural}}}\;`
\Psi^{-1}]
\efig
\end{align*}
is commutative. To show this, let us consider the following homomorphism
\begin{align*}
  \iota=(\Psi^{-1}\circ \psi  \circ \varphi) - \widehat{w}^{\mfrak{p}^\natural}_{\smash{\kappa_{\smash{\mfrak{p}^\natural}},\mfrak{g}^\natural}}
\end{align*}
of $\N_0$-graded vertex algebras. By the definition of $\Psi$, we have $\iota(a(z))=0$ for $a\in \mfrak{u}^\natural$. Further, for $a \in \mfrak{l}^\natural$ we may write
\begin{align*}
   (\Psi \circ \iota)(a(z)) &= (\psi \circ \varphi)(a(z)) - (\Psi \circ \widehat{w}^{\mfrak{p}^\natural}_{\smash{\kappa_{\smash{\mfrak{p}^\natural}},\mfrak{g}^\natural}})(a(z)) \\
   & = \widetilde{w}^\mfrak{p}_{\kappa_\mfrak{p},\mfrak{g}}(a(z)) + \varphi_0(a(z)) - \sum_{\alpha, \beta \in \smash{\Delta_+^{\mfrak{u}^\natural}}} [e_\beta,a]_\alpha \normOrd{a_\beta^*(z) \Psi(a_\alpha(z))} - \Psi(a_{\mfrak{z}(\mfrak{l}^\natural)}(z)) \\
   & = \sum_{\alpha,\beta \in \Delta_+^\mfrak{u}} [e_\beta,a]_\alpha\normOrd{a_\beta^*(z) a_\alpha(z)} + a_{\mfrak{z}(\mfrak{l})}(z) + \varphi_0(a(z)) - \sum_{\alpha,\beta \in \smash{\Delta_+^{\mfrak{u}^\natural}}} [e_\beta,a]_\alpha \normOrd{a_\beta^*(z)a_\alpha(z)} \\
   &\quad +\sum_{\alpha,\beta \in \smash{\Delta_+^{\mfrak{u}^\natural}}} [e_\beta,a]_\alpha \normOrd{a_\beta^*(z)\varphi_0(e_\alpha(z))} - a_{\mfrak{z}(\mfrak{l})}(z) \\
   &= \sum\limits_{\alpha,\beta \in \Delta_+^\mfrak{u} \setminus \smash{\Delta_+^{\mfrak{u}^\natural}}} [e_\beta,a]_\alpha \normOrd{a_\beta^*(z)a_\alpha(z)} + \varphi_0(a(z)) +\sum_{\alpha,\beta \in \smash{\Delta_+^{\mfrak{u}^\natural}}} [e_\beta,a] \normOrd{a_\beta^*(z) \varphi_0(e_\alpha(z))},
\end{align*}
where we used Theorem \ref{thm:FF homomorphism 1-graded} and Lemma \ref{lem:for_induction}. Therefore, it follows that $\iota(a(z))$ commutes with $a_\alpha^*(z)$ for $\alpha \in \smash{\Delta^{\mfrak{u}^\natural}_+}$ and also with $h(z)$ for $h \in \mfrak{z}(\mfrak{l}^\natural)$, which means that $\iota(a(z))$ is contained in the vertex subalgebra generated by $a^*_\alpha(z)$ for $\alpha \in \smash{\Delta^{\mfrak{u}^\natural}_+}$. By a similar calculation for $a \in \widebar{\mfrak{u}}^\natural$, we obtain that $\iota(a(z))$ is also contained in the vertex subalgebra generated by $a^*_\alpha(z)$ for $\alpha \in \smash{\Delta^{\mfrak{u}^\natural}_+}$. Hence, the homomorphism $\iota$ of $\N_0$-graded vertex algebras maps $\mcal{V}^{\kappa_{\smash{\mfrak{p}^\natural}}}\!(\mfrak{g}^\natural)$ into the commutative vertex subalgebra generated by $a_\alpha^*(z)$ for $\alpha \in \smash{\Delta^{\mfrak{u}^\natural}_+}$, which implies that $\iota(a(z))=0$ for $a \in [\mfrak{g}^\natural,\mfrak{g}^\natural]$. Moreover, for $a \in \mfrak{z}(\mfrak{g}^\natural)$ we have that $a \in \mfrak{l}^\natural$ and by the calculation above we obtain that $\iota(a(z))$ commutes also with $a_\alpha(z)$ for $\alpha \in \smash{\Delta^{\mfrak{u}^\natural}_+}$, which gives us that $\iota(a(z))=0$. We have shown that the above diagram commutes. Now, the exactness of the functor $H_{DS,f_\theta}^0(?)$ implies that we have an embedding
\begin{align*}
H_{DS,f_\theta}^0(\widetilde{\mcal{V}}^{\kappa_{\mfrak{p}}}\!(\mfrak{g})) \rightarrow
H_{DS,f_\theta}^0(\mcal{M}_{\mfrak{u}} \otimes_\C \mcal{V}^{\kappa_\mfrak{p}+{1\over 2}\kappa_\mfrak{g}}\!(\mfrak{z}(\mfrak{l}))) \simeq \mcal{M}_{\mfrak{u}^\natural}\! \otimes_\C   \mcal{V}_{\kappa_{\smash{\mfrak{p}^\natural}}+ {1\over 2}\kappa_{\smash{\mfrak{g}^\natural}}}\!(\mfrak{z}(\mfrak{l}^\natural)).
\end{align*}
Therefore, by (ii) and the above commutative diagram we conclude that
\begin{align*}
  H^0_{DS,f_\theta}(\smash{\widetilde{\mcal{V}}}^{\kappa_\mfrak{p}}\!(\mfrak{g})) \simeq \im \widehat{w}^{\mfrak{p}^\natural}_{\smash{\kappa_{\smash{\mfrak{p}^\natural}},\mfrak{g}^\natural}} = \smash{\widehat{\mcal{V}}}^{\kappa_{\smash{\mfrak{p}^\natural}}}\! (\mfrak{g}^\natural),
\end{align*}
which was our goal.}

\corollary{The level $\kappa_\mfrak{p}$ is collapsing for $f_\theta$.}

Let us recall that $\Pi \setminus \Sigma = \{\alpha\}$ for some root $\alpha \in \Pi$ and that $\Pi^\natural \setminus \Sigma^\natural = \Pi^\natural \cap \{\alpha\}$. Therefore, we have two possibilities. If $\alpha \notin \Pi^\natural$, then $\Sigma^\natural = \Pi^\natural$. We set $\tilde{\mfrak{g}}^\natural = \smash{\mfrak{g}^\natural_0}$, see Table \ref{tab:associated varieties}. On the other hand, if $\alpha \in \Pi^\natural$, then $\Sigma^\natural = \Pi^\natural \setminus \{\alpha\}$, which means that there exists an index $j\in \{1,2,\dots,r\}$ such that $\mfrak{g}_\alpha \subset \smash{\mfrak{g}^\natural_j}$. Moreover, we have that $\smash{(\kappa_{\mfrak{p}^\natural} + {1\over 2}\kappa_{\mfrak{g}^\natural})_{|\smash{\mfrak{g}^\natural_0}}} = 0$. We set $\tilde{\mfrak{g}}^\natural = \smash{\mfrak{g}^\natural_j}$, see Table \ref{tab:associated varieties}. It easily follows that $\mfrak{u}^\natural \subset \tilde{\mfrak{g}}^\natural$ and $\widebar{\mfrak{u}}^\natural \subset \tilde{\mfrak{g}}^\natural$. We introduce an opposite standard parabolic subalgebra $\tilde{\mfrak{p}}^\natural$ of $\tilde{\mfrak{g}}^\natural$ by
\begin{align*}
  \tilde{\mfrak{p}}^\natural = \mfrak{p}^\natural \cap \tilde{\mfrak{g}}^\natural.
\end{align*}
We have also a triangular decomposition
\begin{align*}
  \tilde{\mfrak{g}}^\natural = \widebar{\mfrak{u}}^\natural \oplus \tilde{\mfrak{l}}^\natural \oplus \mfrak{u}^\natural
\end{align*}
of $\tilde{\mfrak{g}}^\natural$, where $\smash{\tilde{\mfrak{l}}^\natural} = \mfrak{l}^\natural \cap \tilde{\mfrak{g}}^\natural$ is the Levi subalgebra of $\tilde{\mfrak{p}}^\natural$. Further, by using Theorem \ref{thm:FF homomorphism 1-graded} we get that the homomorphism
\begin{align*}
  \widehat{w}^{\mfrak{p}^\natural}_{\smash{\kappa_{\smash{\mfrak{p}^\natural}},\mfrak{g}^\natural}} \colon \mcal{V}^{\kappa_{\smash{\mfrak{p}^\natural}}}\!(\mfrak{g}^\natural) \rarr \mcal{M}_{\mfrak{u}^\natural}\! \otimes_\C \mcal{V}_{\kappa_{\smash{\mfrak{p}^\natural}}+ {1\over 2}\kappa_{\smash{\mfrak{g}^\natural}}}\!(\mfrak{z}(\mfrak{l}^\natural))
\end{align*}
of $\N_0$-graded vertex algebras coincides with the homomorphism
\begin{align*}
  \mcal{V}^{\kappa_{\smash{\mfrak{p}^\natural}}}\!(\mfrak{g}^\natural) \rarr \mcal{V}^{\kappa_{\smash{\mfrak{p}^\natural}}}\!(\tilde{\mfrak{g}}^\natural) \rarr \mcal{M}_{\mfrak{u}^\natural}\! \otimes_\C \mcal{V}^{\kappa_{\smash{\mfrak{p}^\natural}}+ {1\over 2}\kappa_{\smash{\mfrak{g}^\natural}}}\!(\mfrak{z}(\tilde{\mfrak{l}}^\natural))
\end{align*}
of $\N_0$-graded vertex algebras. Let us note that the vertex algebra  $\mcal{V}^{\kappa_{\smash{\mfrak{p}^\natural}}+ {1\over 2}\kappa_{\smash{\mfrak{g}^\natural}}}\!(\mfrak{z}(\tilde{\mfrak{l}}^\natural))$ is simple. If we denote by $\kappa_{\tilde{\mfrak{p}}^\natural}$ the restriction of $\kappa_{\mfrak{p}^\natural}$ to $\tilde{\mfrak{g}}^\natural$, we get the following consequence of Theorem \ref{thm:DS reduction}.
\medskip

\proposition{\label{prop:DS reduction}We have $H^0_{DS,f_\theta}(\smash{\widetilde{\mcal{V}}}^{\kappa_\mfrak{p}}\!(\mfrak{g})) \simeq \smash{\widetilde{\mcal{V}}}^{\kappa_{\smash{\tilde{\mfrak{p}}^\natural}}}\! (\tilde{\mfrak{g}}^\natural)$.}

As $(\cdot\,,\cdot)$ is a non-degenerate $\mfrak{g}$-invariant symmetric bilinear form on $\mfrak{g}$, it induces an isomorphism $\phi \colon \mfrak{g} \rarr \mfrak{g}^*$ of $\mfrak{g}$-modules. The vertex algebras $\smash{\widetilde{\mcal{V}}}^{\kappa_\mfrak{p}}\!(\mfrak{g})$ and $\smash{\widetilde{\mcal{V}}}^{\kappa_{\smash{\tilde{\mfrak{p}}^\natural}}}\! (\tilde{\mfrak{g}}^\natural)$ are isomorphic to $\smash{\mcal{N}^{\kappa_\mfrak{p}}_\mfrak{p}}\!(\mfrak{g})$ and $\smash{\mcal{N}^{\kappa_{\smash{\tilde{\mfrak{p}}^\natural}}}_{\tilde{\mfrak{p}}^\natural}} \!(\tilde{\mfrak{g}}^\natural)$, respectively. Hence, by Theorem \ref{thm:associated variety} we have that
\begin{align*}
  X_{\smash{\widetilde{\mcal{V}}^{\kappa_\mfrak{p}}\!(\mfrak{g})}} \simeq \phi(\widebar{\mcal{S}_{\mfrak{p}}}) \qquad \text{and} \qquad X_{\smash{\widetilde{\mcal{V}}^{\kappa_{\smash{\tilde{\mfrak{p}}^\natural}}}\!(\tilde{\mfrak{g}}^\natural)}} \simeq \phi(\overline{\mcal{S}_{\smash{\tilde{\mfrak{p}}^\natural}}}),
\end{align*}
which together with Theorem \ref{thm:DS functor min} and Proposition \ref{prop:DS reduction} gives rise to an isomorphism
\begin{align*}
  \phi(\widebar{\mcal{S}_{\mfrak{p}}}) \cap S_{f_\theta} \simeq X_{\smash{H^0_{DS,f_\theta}\!(\widetilde{\mcal{V}}^{\kappa_\mfrak{p}}\!(\mfrak{g}))}} \simeq X_{\smash{\widetilde{\mcal{V}}^{\kappa_{\smash{\tilde{\mfrak{p}}^\natural}}}\!(\tilde{\mfrak{g}}^\natural)}} \simeq \phi(\overline{\mcal{S}_{\smash{\tilde{\mfrak{p}}^\natural}}})
\end{align*}
of Poisson varieties, where $\smash{S_{f_\theta}}$ is the Slodowy slice associated to the $\mfrak{sl}_2$-triple $(e_\theta,h_\theta,f_\theta)$, see the formula \eqref{eq:Slodowy slice def}. Further, by using the fact that $S_{f_\theta} = \phi(f_\theta + \mfrak{g}^{e_\theta})$, we get an equality of sets
\begin{align}
  \widebar{\mcal{S}_{\mfrak{p}}} \cap (f_\theta + \mfrak{g}^{e_\theta}) = f_\theta + \overline{\mcal{S}_{\smash{\tilde{\mfrak{p}}^\natural}}}, \label{eq:Slodowy slice condition}
\end{align}
where we used \eqref{eq:Slodowy slice} to identify $\overline{\mcal{S}_{\smash{\tilde{\mfrak{p}}^\natural}}}$ with the corresponding subset of $f_\theta + \mfrak{g}^{e_\theta}$. Moreover, since the associated varieties of $X_{\smash{\mcal{V}_{\kappa_\mfrak{p}}\!(\mfrak{g})}}$ and $X_{\smash{\mcal{V}_{\kappa_{\smash{\tilde{\mfrak{p}}^\natural}}}\!(\tilde{\mfrak{g}}^\natural)}}$ are subsets of $\phi(\widebar{\mcal{S}_{\mfrak{p}}})$ and $\phi(\overline{\mcal{S}_{\smash{\tilde{\mfrak{p}}^\natural}}})$, respectively, they are linked through the condition
\begin{align*}
  X_{\mfrak{p}} \cap (f_\theta + \mfrak{g}^{e_\theta}) = f_\theta + X_{\smash{\tilde{\mfrak{p}}^\natural}},
\end{align*}
where $X_{\smash{\mcal{V}_{\kappa_\mfrak{p}}\!(\mfrak{g})}} = \phi(X_\mfrak{p})$ and $X_{\smash{\mcal{V}_{\kappa_{\smash{\tilde{\mfrak{p}}^\natural}}}\!(\tilde{\mfrak{g}}^\natural)}} = \phi(X_{\smash{\tilde{\mfrak{p}}^\natural}})$. Besides, we have that $X_\mfrak{p}$ is a $G$-invariant closed conical subset of $\mfrak{g}$. In fact, this condition is quite restrictive as we may see from the subsequent lemma.
\medskip

\lemma{\label{lem:G-invariant subset}Let $X$ be a $G$-invariant closed conical subset of $\widebar{\mcal{S}_{\mfrak{p}}}$. Then we have either $X = \widebar{\mcal{S}_{\mfrak{p}}}$ or $X \subset \widebar{\mcal{O}_{\mfrak{p}}}$.}

\proof{First, let us recall that
\begin{align*}
  \widebar{\mcal{S}_{\mfrak{p}}} = \widebar{G.{\C^*}h}  = G.\C^*h \cup \widebar{\mcal{O}_{\mfrak{p}}}
\end{align*}
by \cite[Lemma 2.1]{Arakawa-Moreau2017}, where $\mfrak{z}(\mfrak{l})=\C h$ for some $h \in \mfrak{h}$. Further, let us assume that $X \not\subset \widebar{\mcal{O}_{\mfrak{p}}}$. Then there exist $g\in G$ and $\xi \in \C^*$ such that $\xi\Ad(g)(h) \in X$. Since $X$ is a $G$-invariant conical subset of
$\widebar{\mcal{S}_{\mfrak{p}}}$, we have that $G.\C^*h \subset X$. Moreover, $X$ is a closed subset of $\widebar{\mcal{S}_{\mfrak{p}}}$, which implies that $X =\widebar{\mcal{S}_{\mfrak{p}}}$.}

\lemma{\label{lem:nilpotent}Let $f_\theta + x \in \mcal{N}(\mfrak{g})$ with $x \in \mfrak{g}^\natural$. Then we have $x \in \mcal{N}(\mfrak{g}^\natural)$.}

\proof{Let us assume that $f_\theta + x \in \mcal{N}(\mfrak{g})$ and $x \in \mfrak{g}^\natural$. Since $f_\theta$ and $x$ commute, we may write
\begin{align*}
  (\ad(f_\theta+x)_{|\mfrak{g}^\natural}\!)^n = \sum_{k=0}^n \binom{n}{k} (\ad(x)_{|\mfrak{g}^\natural}\!)^{n-k}(\ad(f_\theta)|_{\mfrak{g}^\natural}\!)^k = (\ad(x)_{|\mfrak{g}^\natural}\!)^n
\end{align*}
for $n \in \N_0$, which gives us that $\smash{\ad(x)_{|\mfrak{g}^\natural}}$ is nilpotent. Further, the Lie algebra $\mfrak{g}^\natural$ is reductive, which enables us to write $x = x_s + x_n$ with $x_s \in \mfrak{z}(\mfrak{g}^\natural)$ and $x_n \in [\mfrak{g}^\natural, \mfrak{g}^\natural]$. Hence, we have $x_n \in \mcal{N}(\mfrak{g}^\natural)$. Since $\mfrak{g}$ is a finite-dimensional $[\mfrak{g}^\natural,\mfrak{g}^\natural]$-module and $x_n$ is a nilpotent element of $[\mfrak{g}^\natural,\mfrak{g}^\natural]$, we obtain that $x_n$ is also a nilpotent element of $\mfrak{g}$. In addition, by \cite[Lemma 2.1.2]{Collingwood-McGovern1993-book} we have $\mfrak{z}(\mfrak{g}^\natural) \subset \mfrak{h}$, which implies that $x_s$ is a semisimple element of $\mfrak{g}$. As $f_\theta + x_n$ is also nilpotent element of $\mfrak{g}$, $x_s$ is a semisimple element of $\mfrak{g}$ and $[f_\theta + x_n, x_s]=0$, we have that $x_s$ is the semisimple part in the Jordan decomposition of $f_\theta + x$. Hence, we have immediately that $x_s=0$ and $x \in \mcal{N}(\mfrak{g}^\natural)$.}

\lemma{\label{lem:Slodowy Richardson}We have
\begin{align*}
  \widebar{\mcal{O}_{\mfrak{p}}} \cap (f_\theta + \mfrak{g}^{e_\theta}) = f_\theta + \overline{\mcal{O}_{\smash{\tilde{\mfrak{p}}^\natural}}} ,
\end{align*}
where $\mcal{O}_\mfrak{p}$ and $\mcal{O}_{\smash{\tilde{\mfrak{p}}^\natural}}$ are the corresponding Richardson orbits in $\mcal{S}_\mfrak{p}$ and $\mcal{S}_{\smash{\tilde{\mfrak{p}}^\natural}}$, respectively.}

\proof{Since we have $\overline{\mcal{S}_{\smash{\tilde{\mfrak{p}}^\natural}}} \subset \tilde{\mfrak{g}}^\natural$, we may rewrite the equality \eqref{eq:Slodowy slice condition} into the form
\begin{align*}
  \widebar{\mcal{S}_{\mfrak{p}}} \cap (f_\theta + \tilde{\mfrak{g}}^\natural) = f_\theta + \overline{\mcal{S}_{\smash{\tilde{\mfrak{p}}^\natural}}},
\end{align*}
which enables us to write
\begin{align*}
  \widebar{\mcal{O}_{\mfrak{p}}} \cap (f_\theta + \mfrak{g}^{e_\theta}) &= \widebar{\mcal{O}_{\mfrak{p}}} \cap (f_\theta + \tilde{\mfrak{g}}^\natural) = \mcal{N}(\mfrak{g}) \cap \widebar{\mcal{S}_{\mfrak{p}}} \cap (f_\theta + \tilde{\mfrak{g}}^\natural) \\
  &= \mcal{N}(\mfrak{g}) \cap (f_\theta + \overline{\mcal{S}_{\smash{\tilde{\mfrak{p}}^\natural}}}) = f_\theta + (\overline{\mcal{S}_{\smash{\tilde{\mfrak{p}}^\natural}}} \cap \mcal{N}(\tilde{\mfrak{g}}^\natural)) = f_\theta + \overline{\mcal{O}_{\smash{\tilde{\mfrak{p}}^\natural}}},
\end{align*}
where we used Lemma \ref{lem:nilpotent} in the next-to-last equality.}

Let us consider a non-zero nilpotent orbit $\mcal{O}$ of $\mfrak{g}$ such that $\mcal{O} \subset \widebar{\mcal{O}_{\mfrak{p}}}$. Since we have $f_\theta \in \widebar{\mcal{O}}$, it implies that $\mcal{O} \cap (f_\theta + \mfrak{g}^{e_\theta}) \neq \emptyset$. Therefore, by Lemma \ref{lem:Slodowy Richardson} we may assign to $\mcal{O}$ a subset $\chi(\mcal{O})$ of $\overline{\mcal{O}_{\smash{\tilde{\mfrak{p}}^\natural}}}$ through the condition
\begin{align*}
  \mcal{O} \cap (f_\theta + \mfrak{g}^{e_\theta}) = f_\theta + \chi(\mcal{O}).
\end{align*}
The following lemmas give us basic properties of the set $\chi(\mcal{O})$.
\medskip

\lemma{\label{lem:orbit correspondence}Let $\mcal{O}$ be a non-zero nilpotent orbit of $\mfrak{g}$ such that $\mcal{O} \subset \widebar{\mcal{O}_{\mfrak{p}}}$. Then $\chi(\mcal{O})$ is a union of nilpotent orbits of $\tilde{\mfrak{g}}^\natural$. In addition, if $\chi(\mcal{O}_1) \cap \chi(\mcal{O}_2) \neq \emptyset$, then we have $\mcal{O}_1 = \mcal{O}_2$.}

\proof{As $\mcal{O}$ is a non-zero nilpotent orbit of $\mfrak{g}$, we know that $\chi(\mcal{O})$ is not the empty set. We need to show that the set $\chi(\mcal{O})$ is $\smash{\widetilde{G}^\natural}$-invariant. Let $g \in \smash{\widetilde{G}^\natural}$ and $x \in \chi(\mcal{O})$. Then we have $\Ad(g)(f_\theta)=f_\theta$, which follows from the fact that $\tilde{\mfrak{g}}^\natural \subset \mfrak{g}^{f_\theta}$. Therefore, we may write $\Ad(g)(f_\theta+x) = f_\theta + \Ad(g)(x)$, which implies that $\Ad(g)(x) \in \chi(\mcal{O})$, and so $\chi(\mcal{O})$ is $\smash{\widetilde{G}^\natural}$-invariant. In addition, by using the fact $\chi(\mcal{O}) \subset \overline{\mcal{O}_{\smash{\tilde{\mfrak{p}}^\natural}}}$, we get that $\chi(\mcal{O})$ is a union of nilpotent orbits of $\tilde{\mfrak{g}}^\natural$.

Now, let us assume that $\chi(\mcal{O}_1) \cap \chi(\mcal{O}_2) \neq \emptyset$ for non-zero nilpotent orbits $\mcal{O}_1$ and $\mcal{O}_2$ of $\mfrak{g}$ such that $\mcal{O}_1, \mcal{O}_2 \subset \widebar{\mcal{O}_{\mfrak{p}}}$. Hence, for an element $x \in \chi(\mcal{O}_1) \cap \chi(\mcal{O}_2)$ we get $f_\theta + x \in \mcal{O}_1 \cap \mcal{O}_2$, which means that $\mcal{O}_1 = \mcal{O}_2$.}

\lemma{\label{lem:orbit correspondence inverse}Let $\mcal{O}^\natural$ be a nilpotent orbit of $\tilde{\mfrak{g}}^\natural$ such that $\mcal{O}^\natural \subset \overline{\mcal{O}_{\smash{\tilde{\mfrak{p}}^\natural}}}$. Then there exists a unique nilpotent orbit $\mcal{O}$ of $\mfrak{g}$ satisfying $\mcal{O}^\natural \subset \chi(\mcal{O})$.}

\proof{For $x \in \mcal{O}^\natural$, we have $f_\theta + x \in \widebar{\mcal{O}_{\mfrak{p}}}$ by Lemma \ref{lem:Slodowy Richardson}. If we set $\mcal{O}=G.(f_\theta+x)$, we obtain that $x \in \chi(\mcal{O})$, which gives us $\mcal{O}^\natural \subset \chi(\mcal{O})$. The uniqueness follows immediately from Lemma \ref{lem:orbit correspondence}.}

Let $\mcal{O}$ be a nilpotent orbit of $\mfrak{g}$. We denote by $\mcal{N}_\mfrak{g}(\mcal{O})$ the set of all nilpotent orbits of $\mfrak{g}$ which are contained in $\widebar{\mcal{O}}$. Analogously, we define $\mcal{N}_{\smash{\tilde{\mfrak{g}}^\natural}}(\mcal{O}^\natural)$ for a nilpotent orbit $\mcal{O}^\natural$ of $\tilde{\mfrak{g}}^\natural$. As a direct consequence of Lemma \ref{lem:orbit correspondence} and Lemma \ref{lem:orbit correspondence inverse} we have a decomposition
\begin{align*}
   \overline{\mcal{O}_{\smash{\tilde{\mfrak{p}}^\natural}}} = \bigsqcup_{\mcal{O} \in \mcal{N}_\mfrak{g}(\mcal{O}_\mfrak{p}) \setminus \{\mcal{O}_{\rm zero}\}} \chi(\mcal{O}).
\end{align*}
Further, let us suppose that $\#\mcal{N}_\mfrak{g}(\mcal{O}_\mfrak{p}) = \#\mcal{N}_{\smash{\tilde{\mfrak{g}}^\natural}}(\mcal{O}_{\smash{\tilde{\mfrak{p}}^\natural}}) + 1$. Let us note also that the assumption is satisfied in all cases. It implies that $\chi(\mcal{O})$ is a nilpotent orbit of $\tilde{\mfrak{g}}^\natural$ for any non-zero nilpotent orbit $\mcal{O}$ of $\mfrak{g}$ satisfying $\mcal{O} \subset \smash{\widebar{\mcal{O}_{\mfrak{p}}}}$. Hence, we get a one-to-one correspondence between  $\mcal{N}_\mfrak{g}(\mcal{O}) \setminus \{\mcal{O}_{\rm zero}\}$ and $\mcal{N}_{\smash{\tilde{\mfrak{g}}^\natural}}(\mcal{O}_{\smash{\tilde{\mfrak{p}}^\natural}})$. In addition, the mapping $\chi$ gives rise to an isomorphism of the Hasse diagrams for $\mcal{N}_\mfrak{g}(\mcal{O}_\mfrak{p}) \setminus \{\mcal{O}_{\rm zero}\}$ and $\mcal{N}_{\smash{\tilde{\mfrak{g}}^\natural}}(\mcal{O}_{\smash{\tilde{\mfrak{p}}^\natural}})$, see Figure \ref{fig:hasse correspondence}.
\medskip

\begin{table}[ht]
\centering
\subcaptionbox{$\mfrak{sl}_{n+1} \rarr \mfrak{gl}_1, n \geq 2$\label{fig:sl}}[0.49\textwidth]
{\begin{tikzpicture}
[yscale=1.3,xscale=2,vector/.style={circle,draw=white,fill=black,ultra thick, inner sep=0.8mm},vector2/.style={circle,draw=white,fill=white,ultra thick, inner sep=1mm}]
\begin{scope}
  \node (A) at (0,1) {$\mcal{O}_{{\rm zero}}$};
  \node (B) at (1,1) {$\mcal{O}_{[2,1^{n-1}]}$};
  \node (C) at (1,0) {$\mcal{O}_{{\rm zero}}$};
  \draw [thin, -] (A) -- (B);
  \draw [thin, ->] (B) -- (C);
  \node at (-0.6,1) {$\mfrak{sl}_{n+1}$};
  \node at (-0.6,0) {$\mfrak{gl}_1$};
\end{scope}
\end{tikzpicture}}
\subcaptionbox{$\mfrak{e}_6 \rarr \mfrak{sl}_6$\label{fig:e_6}}
{\begin{tikzpicture}
[yscale=1.3,xscale=2,vector/.style={circle,draw=white,fill=black,ultra thick, inner sep=0.8mm},vector2/.style={circle,draw=white,fill=white,ultra thick, inner sep=1mm}]
\begin{scope}
  \node (A) at (0,1) {$\mcal{O}_{{\rm zero}}$};
  \node (B) at (1,1) {$\mcal{O}_{A_1}$};
  \node (C) at (2,1) {$\mcal{O}_{2A_1}$};
  \node (D) at (1,0) {$\mcal{O}_{{\rm zero}}$};
  \node (E) at (2,0) {$\mcal{O}_{[2,1^4]}$};
  \draw [thin, -] (A) -- (B);
  \draw [thin, -] (B) -- (C);
  \draw [thin, -] (D) -- (E);
  \draw [thin, ->] (B) -- (D);
  \draw [thin, ->] (C) -- (E);
  \node at (-0.6,1) {$\mfrak{e}_6$};
  \node at (-0.6,0) {$\mfrak{sl}_6$};
\end{scope}
\end{tikzpicture}}
\vspace{6mm}

\subcaptionbox{$\mfrak{sl}_{2n} \rarr \mfrak{gl}_{2n-2}, n \geq 2$\label{fig:sl_2n}}
{\begin{tikzpicture}
[yscale=1.3,xscale=2,vector/.style={circle,draw=white,fill=black,ultra thick, inner sep=0.8mm},vector2/.style={circle,draw=white,fill=white,ultra thick, inner sep=1mm}]
\begin{scope}
  \node (A) at (0,1) {$\mcal{O}_{{\rm zero}}$};
  \node (B) at (1,1) {$\mcal{O}_{[2,1^{2n-2}]}$};
  \node (C) at (2,1) {\dots};
  \node (D) at (3,1) {$\mcal{O}_{[2^{n-1},1^2]}$};
  \node (E) at (4,1) {$\mcal{O}_{[2^n]}$};
  \node (F) at (1,0) {$\mcal{O}_{{\rm zero}}$};
  \node (G) at (2,0) {\dots};
  \node (H) at (3,0) {$\mcal{O}_{[2^{n-2},1^2]}$};
  \node (I) at (4,0) {$\mcal{O}_{[2^{n-1}]}$};
  \draw [thin, -] (A) -- (B);
  \draw [thin, -] (B) -- (C);
  \draw [thin, -] (C) -- (D);
  \draw [thin, -] (D) -- (E);
  \draw [thin, -] (F) -- (G);
  \draw [thin, -] (G) -- (H);
  \draw [thin, -] (H) -- (I);
  \draw [thin, ->] (B) -- (F);
  \draw [thin, ->] (C) -- (G);
  \draw [thin, ->] (D) -- (H);
  \draw [thin, ->] (E) -- (I);
  \node at (-0.6,1) {$\mfrak{sl}_{2n}$};
  \node at (-0.6,0) {$\mfrak{sl}_{2n-2}$};
\end{scope}
\end{tikzpicture}}
\vspace{6mm}

\subcaptionbox{$\mfrak{so}_{2n+1} \rarr \mfrak{sl}_2, n \geq 3$\label{fig:so 2n+1}}[0.49\textwidth]
{\begin{tikzpicture}
[yscale=1.3,xscale=2,vector/.style={circle,draw=white,fill=black,ultra thick, inner sep=0.8mm},vector2/.style={circle,draw=white,fill=white,ultra thick, inner sep=1mm}]
\begin{scope}
  \node (A) at (0,1) {$\mcal{O}_{{\rm zero}}$};
  \node (B) at (1,1) {$\mcal{O}_{[2^2,1^{2n-3}]}$};
  \node (C) at (2,1) {$\mcal{O}_{[3,1^{2n-2}]}$};
  \node (E) at (1,0) {$\mcal{O}_{{\rm zero}}$};
  \node (F) at (2,0) {$\mcal{O}_{[2]}$};
  \draw [thin, -] (A) -- (B);
  \draw [thin, -] (B) -- (C);
  \draw [thin, -] (E) -- (F);
  \draw [thin, ->] (B) -- (E);
  \draw [thin, ->] (C) -- (F);
  \node at (-0.6,1) {$\mfrak{so}_{2n+1}$};
  \node at (-0.6,0) {$\mfrak{sl}_2$};
\end{scope}
\end{tikzpicture}}
\subcaptionbox{$\mfrak{so}_{2n} \rarr \mfrak{sl}_2, n \geq 4$\label{fig:so_2n}}[0.49\textwidth]
{\begin{tikzpicture}
[yscale=1.3,xscale=2,vector/.style={circle,draw=white,fill=black,ultra thick, inner sep=0.8mm},vector2/.style={circle,draw=white,fill=white,ultra thick, inner sep=1mm}]
\begin{scope}
  \node (A) at (0,1) {$\mcal{O}_{{\rm zero}}$};
  \node (B) at (1,1) {$\mcal{O}_{[2^2,1^{2n-4}]}$};
  \node (C) at (2,1) {$\mcal{O}_{[3,1^{2n-3}]}$};
  \node (E) at (1,0) {$\mcal{O}_{{\rm zero}}$};
  \node (F) at (2,0) {$\mcal{O}_{[2]}$};
  \draw [thin, -] (A) -- (B);
  \draw [thin, -] (B) -- (C);
  \draw [thin, -] (E) -- (F);
  \draw [thin, ->] (B) -- (E);
  \draw [thin, ->] (C) -- (F);
  \node at (-0.6,1) {$\mfrak{so}_{2n}$};
  \node at (-0.6,0) {$\mfrak{sl}_2$};
\end{scope}
\end{tikzpicture}}
\vspace{6mm}

\subcaptionbox{$\mfrak{sp}_{2n} \rarr \mfrak{sp}_{2n-2}, n \geq 2$\label{fig:sp_2n}}
{\begin{tikzpicture}
[yscale=1.3,xscale=2,vector/.style={circle,draw=white,fill=black,ultra thick, inner sep=0.8mm},vector2/.style={circle,draw=white,fill=white,ultra thick, inner sep=1mm}]
\begin{scope}
  \node (A) at (0,1) {$\mcal{O}_{{\rm zero}}$};
  \node (B) at (1,1) {$\mcal{O}_{[2,1^{2n-2}]}$};
  \node (C) at (2,1) {\dots};
  \node (D) at (3,1) {$\mcal{O}_{[2^{n-1},1^2]}$};
  \node (E) at (4,1) {$\mcal{O}_{[2^n]}$};
  \node (F) at (1,0) {$\mcal{O}_{{\rm zero}}$};
  \node (G) at (2,0) {\dots};
  \node (H) at (3,0) {$\mcal{O}_{[2^{n-2},1^2]}$};
  \node (I) at (4,0) {$\mcal{O}_{[2^{n-1}]}$};
  \draw [thin, -] (A) -- (B);
  \draw [thin, -] (B) -- (C);
  \draw [thin, -] (C) -- (D);
  \draw [thin, -] (D) -- (E);
  \draw [thin, -] (F) -- (G);
  \draw [thin, -] (G) -- (H);
  \draw [thin, -] (H) -- (I);
  \draw [thin, ->] (B) -- (F);
  \draw [thin, ->] (C) -- (G);
  \draw [thin, ->] (D) -- (H);
  \draw [thin, ->] (E) -- (I);
  \node at (-0.6,1) {$\mfrak{sp}_{2n}$};
  \node at (-0.6,0) {$\mfrak{sp}_{2n-2}$};
\end{scope}
\end{tikzpicture}}
\vspace{6mm}

\subcaptionbox{$\mfrak{so}_{2n} \rarr \mfrak{so}_{2n-4}, n \geq 4, n \in 2\N+1$\label{fig:so_2n odd}}
{\begin{tikzpicture}
[yscale=1.3,xscale=2,vector/.style={circle,draw=white,fill=black,ultra thick, inner sep=0.8mm},vector2/.style={circle,draw=white,fill=white,ultra thick, inner sep=1mm}]
\begin{scope}
  \node (A) at (0,1) {$\mcal{O}_{{\rm zero}}$};
  \node (B) at (1,1) {$\mcal{O}_{[2^2,1^{2n-4}]}$};
  \node (C) at (2,1) {\dots};
  \node (D) at (3,1) {$\mcal{O}_{[2^{n-3},1^6]}$};
  \node (E) at (4,1) {$\mcal{O}_{[2^{n-1},1^2]}$};
  \node (F) at (1,0) {$\mcal{O}_{{\rm zero}}$};
  \node (G) at (2,0) {\dots};
  \node (H) at (3,0) {$\mcal{O}_{[2^{n-5},1^6]}$};
  \node (I) at (4,0) {$\mcal{O}_{[2^{n-3},1^2]}$};
  \draw [thin, -] (A) -- (B);
  \draw [thin, -] (B) -- (C);
  \draw [thin, -] (C) -- (D);
  \draw [thin, -] (D) -- (E);
  \draw [thin, -] (F) -- (G);
  \draw [thin, -] (G) -- (H);
  \draw [thin, -] (H) -- (I);
  \draw [thin, ->] (B) -- (F);
  \draw [thin, ->] (C) -- (G);
  \draw [thin, ->] (D) -- (H);
  \draw [thin, ->] (E) -- (I);
  \node at (-0.6,1) {$\mfrak{so}_{2n}$};
  \node at (-0.6,0) {$\mfrak{so}_{2n-4}$};
\end{scope}
\end{tikzpicture}}
\vspace{6mm}

\subcaptionbox{$\mfrak{so}_{2n} \rarr \mfrak{so}_{2n-4}, n \geq 4, n \in 2\N$\label{fig:so_2n even}}
{\begin{tikzpicture}
[yscale=1.3,xscale=2,vector/.style={circle,draw=white,fill=black,ultra thick, inner sep=0.8mm},vector2/.style={circle,draw=white,fill=white,ultra thick, inner sep=1mm}]
\begin{scope}
  \node (A) at (0,1) {$\mcal{O}_{{\rm zero}}$};
  \node (B) at (1,1) {$\mcal{O}_{[2^2,1^{2n-4}]}$};
  \node (C) at (2,1) {\dots};
  \node (D) at (3,1) {$\mcal{O}_{[2^{n-2},1^4]}$};
  \node (E) at (4,1) {$\mcal{O}^{II}_{[2^n]}$};
  \node (F) at (1,0) {$\mcal{O}_{{\rm zero}}$};
  \node (G) at (2,0) {\dots};
  \node (H) at (3,0) {$\mcal{O}_{[2^{n-4},1^4]}$};
  \node (I) at (4,0) {$\mcal{O}^{II}_{[2^{n-2}]}$};
  \draw [thin, -] (A) -- (B);
  \draw [thin, -] (B) -- (C);
  \draw [thin, -] (C) -- (D);
  \draw [thin, -] (D) -- (E);
  \draw [thin, -] (F) -- (G);
  \draw [thin, -] (G) -- (H);
  \draw [thin, -] (H) -- (I);
  \draw [thin, ->] (B) -- (F);
  \draw [thin, ->] (C) -- (G);
  \draw [thin, ->] (D) -- (H);
  \draw [thin, ->] (E) -- (I);
  \node at (-0.6,1) {$\mfrak{so}_{2n}$};
  \node at (-0.6,0) {$\mfrak{so}_{2n-4}$};
\end{scope}
\end{tikzpicture}}
\vspace{6mm}

\subcaptionbox{$\mfrak{e}_7 \rarr \mfrak{so}_{12}$\label{fig:e_7}}
{\begin{tikzpicture}
[yscale=1.3,xscale=2,vector/.style={circle,draw=white,fill=black,ultra thick, inner sep=0.8mm},vector2/.style={circle,draw=white,fill=white,ultra thick, inner sep=1mm}]
\begin{scope}
  \node (A) at (0,1) {$\mcal{O}_{{\rm zero}}$};
  \node (B) at (1,1) {$\mcal{O}_{A_1}$};
  \node (C) at (2,1) {$\mcal{O}_{2A_1}$};
  \node (D) at (3,1) {$\mcal{O}_{(3A_1)''}$};
  \node (E) at (1,0) {$\mcal{O}_{{\rm zero}}$};
  \node (F) at (2,0) {$\mcal{O}_{[2^2,1^8]}$};
  \node (G) at (3,0) {$\mcal{O}_{[3,1^9]}$};
  \draw [thin, -] (A) -- (B);
  \draw [thin, -] (B) -- (C);
  \draw [thin, -] (C) -- (D);
  \draw [thin, -] (E) -- (F);
  \draw [thin, -] (F) -- (G);
  \draw [thin, ->] (B) -- (E);
  \draw [thin, ->] (C) -- (F);
  \draw [thin, ->] (D) -- (G);
  \node at (-0.6,1) {$\mfrak{e}_7$};
  \node at (-0.6,0) {$\mfrak{so}_{12}$};
\end{scope}
\end{tikzpicture}}
\caption{Correspondence of Hasse diagrams}
\label{fig:hasse correspondence}
\end{table}

\clearpage

Now, we are able to determine the associated variety  $X_{\smash{\mcal{V}_{\kappa_\mfrak{p}}\!(\mfrak{g})}}$ of the simple affine vertex algebra $\mcal{V}_{\kappa_\mfrak{p}}\!(\mfrak{g})$, which is our main theorem.
\medskip

\theorem{\label{MainTheorem}The associated variety $X_{\smash{\mcal{V}_{\kappa_\mfrak{p}}\!(\mfrak{g})}}$ of the simple affine vertex algebra $\mcal{V}_{\kappa_\mfrak{p}}\!(\mfrak{g})$ is given in Table \ref{tab:associated varieties}, where $\kappa_\mfrak{p} = k_\mfrak{p}\kappa_0$. Moreover, if $X_{\smash{\mcal{V}_{\kappa_\mfrak{p}}\!(\mfrak{g})}}=\smash{\widebar{\mcal{S}^*_{\mfrak{p}}}}$, then $\mcal{V}_{\kappa_\mfrak{p}}\!(\mfrak{g}) \simeq \mcal{N}_\mfrak{p}^{\kappa_\mfrak{p}}\!(\mfrak{g})$.}

\proof{Since we have $X_{\smash{\mcal{V}_{\kappa_\mfrak{p}}\!(\mfrak{g})}} = \phi(X_\mfrak{p})$, where $X_\mfrak{p}$ is a $G$-invariant closed conical subset of $\widebar{\mcal{S}_{\mfrak{p}}}$, by Lemma \ref{lem:G-invariant subset} we obtain that either $X_\mfrak{p} = \widebar{\mcal{S}_{\mfrak{p}}}$ or $X_\mfrak{p} \subset \widebar{\mcal{O}_{\mfrak{p}}}$. Besides, we have
\begin{align*}
  X_{\mfrak{p}} \cap (f_\theta + \mfrak{g}^{e_\theta}) = f_\theta + X_{\smash{\tilde{\mfrak{p}}^\natural}},
\end{align*}
where $X_{\smash{\mcal{V}_{\kappa_{\smash{\tilde{\mfrak{p}}^\natural}}} \!(\tilde{\mfrak{g}}^\natural)}} = \phi(X_{\smash{\tilde{\mfrak{p}}^\natural}})$. If $X_\mfrak{p} = \widebar{\mcal{S}_{\mfrak{p}}}$, the equality \eqref{eq:Slodowy slice condition} implies that $X_{\smash{\tilde{\mfrak{p}}^\natural}} = \overline{\mcal{S}_{\smash{\tilde{\mfrak{p}}^\natural}}}$. On the other hand, if $X_\mfrak{p} \subset \widebar{\mcal{O}_{\mfrak{p}}}$, we get that
\begin{align*}
X_\mfrak{p} = \bigsqcup_{\mcal{O} \in \Lambda} \mcal{O}
\end{align*}
for some non-empty subset $\Lambda$ of $\mcal{N}_\mfrak{g}(\mcal{O}_\mfrak{p})$, which by using the correspondence $\chi$ gives us that
\begin{align*}
X_{\smash{\tilde{\mfrak{p}}^\natural}} = \bigsqcup_{\mcal{O} \in \Lambda} \chi(\mcal{O}).
\end{align*}
This gives us a simple criterion how to determine $X_\mfrak{p}$ when we already know  $X_{\smash{\tilde{\mfrak{p}}^\natural}}$.

Further, let us consider a short exact sequence
\begin{align*}
  0 \rarr \mcal{I}_{\kappa_\mfrak{p}}\!(\mfrak{g}) \rarr \widetilde{\mcal{V}}^{\kappa_\mfrak{p}}\!(\mfrak{g}) \rarr \mcal{V}_{\kappa_\mfrak{p}}\!(\mfrak{g}) \rarr 0
\end{align*}
of $\mcal{V}^{\kappa_\mfrak{p}}\!(\mfrak{g})$-modules. Since the functor $H^0_{DS,f_\theta}(?)$ is exact, we obtain a short exact sequence
\begin{align*}
  0 \rarr H^0_{DS,f_\theta}(\mcal{I}_{\kappa_\mfrak{p}}\!(\mfrak{g})) \rarr H^0_{DS,f_\theta}(\widetilde{\mcal{V}}^{\kappa_\mfrak{p}}\!(\mfrak{g})) \rarr H^0_{DS,f_\theta}(\mcal{V}_{\kappa_\mfrak{p}}\!(\mfrak{g})) \rarr 0
\end{align*}
of $\mcal{W}^{\kappa_\mfrak{p}}\!(\mfrak{g},f_\theta)$-modules and hence also $\mcal{V}^{\kappa_{\smash{\tilde{\mfrak{p}}^\natural}}}\!(\tilde{\mfrak{g}}^\natural)$-modules. Moreover, Proposition \ref{prop:DS reduction} and the consequence of Theorem \ref{thm:DS functor min} give us a short exact sequence
\begin{align*}
  0 \rarr H^0_{DS,f_\theta}(\mcal{I}_{\kappa_\mfrak{p}}\!(\mfrak{g})) \rarr \smash{\widetilde{\mcal{V}}}^{\kappa_{\smash{\tilde{\mfrak{p}}^\natural}}}\!(\tilde{\mfrak{g}}^\natural) \rarr \mcal{V}_{\kappa_{\smash{\tilde{\mfrak{p}}^\natural}}}\!(\tilde{\mfrak{g}}^\natural) \rarr 0
\end{align*}
of $\mcal{V}^{\kappa_{\smash{\tilde{\mfrak{p}}^\natural}}}\!(\tilde{\mfrak{g}}^\natural)$-modules. If $\smash{\widetilde{\mcal{V}}}^{\kappa_{\smash{\tilde{\mfrak{p}}^\natural}}} \!(\tilde{\mfrak{g}}^\natural)$ is a simple vertex algebra, we get that $\smash{H^0_{DS,f_\theta}(\mcal{I}_{\kappa_\mfrak{p}}\!(\mfrak{g}))} = 0$, which implies $\mcal{I}_{\kappa_\mfrak{p}}\!(\mfrak{g})=0$ by Theorem \ref{thm:DS functor min}. Therefore, the vertex algebra $\smash{\widetilde{\mcal{V}}}^{\kappa_\mfrak{p}}\!(\mfrak{g})$ is also simple and we get
\begin{align*}
  X_\mfrak{p} = \widebar{\mcal{S}_{\mfrak{p}}}.
\end{align*}
As we have the decomposition $\mfrak{g}^\natural = \smash{\bigoplus_{i=0}^r \mfrak{g}_i^\natural}$ and $\tilde{\mfrak{g}}^\natural = \smash{\mfrak{g}_j^\natural}$ for some $j \in \{0,1,\dots,r\}$, we may write $\kappa_{\tilde{\mfrak{p}}^\natural} = k_{\tilde{\mfrak{p}}^\natural} \smash{\kappa_j^\natural}$, where $k_{\tilde{\mfrak{p}}^\natural}$ is given by Table \ref{tab:g natural classical} and Table \ref{tab:g natural exceptional}, see Table \ref{tab:associated varieties}.

\smallskip

\noindent
i) For $\mfrak{g}=\mfrak{sl}_{n+1}$ and $k_\mfrak{p}=-1$ with $n\geq 2$, we have $\tilde{\mfrak{g}}^\natural \simeq \mfrak{gl}_1$ and $k_{\tilde{\mfrak{p}}^\natural}=\smash{{n-1 \over 2}}$. Since $\smash{\widetilde{\mcal{V}}}^{\kappa_{\smash{\tilde{\mfrak{p}}^\natural}}}\! (\widetilde{\mfrak{g}}^\natural)$ is a simple vertex algebra, it means that also $\smash{\widetilde{\mcal{V}}}^{\kappa_\mfrak{p}}\!(\mfrak{g})$ is a simple vertex algebra, which gives us $X_\mfrak{p} = \smash{\widebar{\mcal{S}_{\mfrak{p}}}}$ with $\mfrak{p} = \mfrak{p}_{\alpha_1}$.
\smallskip

\noindent
ii) For $\mfrak{g} = \mfrak{sl}_{2n}$ and $k_\mfrak{p}=-n$ with $n \geq 2$, we have $\tilde{\mfrak{g}}^\natural \simeq \mfrak{sl}_{2n-2}$ and $k_{\tilde{\mfrak{p}}^\natural}= -n+1$. If $n=2$, the vertex algebra $\smash{\widetilde{\mcal{V}}}^{\kappa_{\smash{\tilde{\mfrak{p}}^\natural}}}\! (\tilde{\mfrak{g}}^\natural)$ is simple. Hence, by induction we obtain that $\smash{\widetilde{\mcal{V}}}^{\kappa_\mfrak{p}}\!(\mfrak{g})$ is a simple vertex algebra, which implies that $X_\mfrak{p} = \smash{\widebar{\mcal{S}_{\mfrak{p}}}}$ with $\mfrak{p} = \mfrak{p}_{\alpha_n}$.
\smallskip

\noindent
iii) For $\mfrak{g}=\mfrak{so}_{2n+1}$ and $k_\mfrak{p}= -2$ with $n\geq 3$, we have $\tilde{\mfrak{g}}^\natural \simeq \mfrak{sl}_2$ and $k_{\tilde{\mfrak{p}}^\natural}=n-\smash{{7 \over 2}}$. The vertex algebra $\smash{\widetilde{\mcal{V}}}^{\kappa_{\smash{\tilde{\mfrak{p}}^\natural}}}\! (\tilde{\mfrak{g}}^\natural)$ is not simple and $X_{\smash{\tilde{\mfrak{p}}^\natural}} = \smash{\widebar{{\mcal{O}}_{[2]}}}$, which by using the correspondence $\chi$ in Figure \ref{fig:so 2n+1} gives us $X_\mfrak{p} = \smash{\widebar{{\mcal{O}}_{[3,1^{2n-2}]}}}$.
\smallskip

\noindent
iv) For $\mfrak{g}=\mfrak{sp}_{2n}$ and $k_\mfrak{p}= -\smash{{n\over 2}}-1$ with $n\geq 2$, we have $\tilde{\mfrak{g}}^\natural \simeq \mfrak{sp}_{2n-2}$ and $k_{\tilde{\mfrak{p}}^\natural}= - \smash{{n-1 \over 2}} -1 $. If $n=2$, the vertex algebra $\smash{\widetilde{\mcal{V}}}^{\kappa_{\smash{\tilde{\mfrak{p}}^\natural}}}\! (\tilde{\mfrak{g}}^\natural)$ is simple. Hence, by induction we obtain that $\smash{\widetilde{\mcal{V}}}^{\kappa_\mfrak{p}}\!(\mfrak{g})$ is a simple vertex algebra, which implies that $X_\mfrak{p} = \smash{\widebar{\mcal{S}_{\mfrak{p}}}}$ with $\mfrak{p} = \mfrak{p}_{\alpha_n}$.
\smallskip

\noindent
v) For $\mfrak{g}=\mfrak{so}_{2n}$ and $k_\mfrak{p}= -2$ with $n\geq 4$, we have $\tilde{\mfrak{g}}^\natural \simeq \mfrak{sl}_2$ and $k_{\tilde{\mfrak{p}}^\natural}=n-4$. The vertex algebra $\smash{\widetilde{\mcal{V}}}^{\kappa_{\smash{\tilde{\mfrak{p}}^\natural}}}\! (\tilde{\mfrak{g}}^\natural)$ is not simple and $X_{\smash{\tilde{\mfrak{p}}^\natural}} = \smash{\widebar{{\mcal{O}}_{\rm zero}}}$, which by using the correspondence $\chi$ in Figure \ref{fig:so_2n} gives us $X_\mfrak{p} = \smash{\widebar{{\mcal{O}}_{\rm min}}}$.
\smallskip

\noindent
vi) For $\mfrak{g}=\mfrak{so}_{2n}$ and $k_\mfrak{p}= -n+2$ with $n\geq 4$, $n \in 2\N+1$, we have $\tilde{\mfrak{g}}^\natural \simeq \mfrak{so}_{2n-4}$ and $k_{\tilde{\mfrak{p}}^\natural}=-n+4$. If $n=5$, the vertex algebra $\smash{\widetilde{\mcal{V}}}^{\kappa_{\smash{\tilde{\mfrak{p}}^\natural}}}\! (\tilde{\mfrak{g}}^\natural)$ is simple by (i). Hence, by induction we obtain that $\smash{\widetilde{\mcal{V}}}^{\kappa_\mfrak{p}}\!(\mfrak{g})$ is a simple vertex algebra, which implies that $X_\mfrak{p} = \smash{\widebar{\mcal{S}_{\mfrak{p}}}}$ with $\mfrak{p} = \mfrak{p}_{\alpha_n}$.
\smallskip

\noindent
vii) For $\mfrak{g}=\mfrak{so}_{2n}$ and $k_\mfrak{p}= -n+2$ with $n\geq 4$, $n \in 2\N$, we have $\tilde{\mfrak{g}}^\natural \simeq \mfrak{so}_{2n-4}$ and $k_{\tilde{\mfrak{p}}^\natural}=-n+4$. If $n=4$, the vertex algebra $\smash{\widetilde{\mcal{V}}}^{\kappa_{\smash{\tilde{\mfrak{p}}^\natural}}}\! (\tilde{\mfrak{g}}^\natural)$ is not simple and $X_{\smash{\tilde{\mfrak{p}}^\natural}} = \smash{\widebar{{\mcal{O}}_{\rm zero}}}$, which by induction and by using the correspondence $\chi$ in Figure \ref{fig:so_2n even} gives us $X_\mfrak{p} = \smash{\widebar{{\mcal{O}}_{[2^{n-2},1^4]}}}$.
\smallskip

\noindent
viii) For $\mfrak{g}=\mfrak{e}_6$ and $k_\mfrak{p}= -4$, we have $\tilde{\mfrak{g}}^\natural \simeq \mfrak{sl}_6$ and $k_{\tilde{\mfrak{p}}^\natural}=-1$. As $\smash{\widetilde{\mcal{V}}}^{\kappa_{\smash{\tilde{\mfrak{p}}^\natural}}}\! (\tilde{\mfrak{g}}^\natural)$ is a simple vertex algebra by (i), it means that also $\smash{\widetilde{\mcal{V}}}^{\kappa_\mfrak{p}}\!(\mfrak{g})$ is a simple vertex algebra, which gives us $X_\mfrak{p} = \smash{\widebar{\mcal{S}_{\mfrak{p}}}}$ with $\mfrak{p} = \mfrak{p}_{\alpha_1}$.
\smallskip

\noindent
ix) For $\mfrak{g}=\mfrak{e}_7$ and $k_\mfrak{p}= -6$, we have $\tilde{\mfrak{g}}^\natural \simeq \mfrak{so}_{12}$ and $k_{\tilde{\mfrak{p}}^\natural}=-2$. The vertex algebra $\smash{\widetilde{\mcal{V}}}^{\kappa_{\smash{\tilde{\mfrak{p}}^\natural}}}\! (\tilde{\mfrak{g}}^\natural)$ is not simple and $X_{\smash{\tilde{\mfrak{p}}^\natural}} = \smash{\widebar{{\mcal{O}}_{[2^2,1^8]}}}$ by (v), which by using the correspondence $\chi$ in Figure \ref{fig:e_7} gives us $X_\mfrak{p} = \smash{\widebar{{\mcal{O}}_{2A_1}}}$.
\smallskip

\noindent
This finishes the proof.}

\vspace{-2mm}

\begin{sidewaystable}
\centering
\renewcommand{\arraystretch}{2.5}
\begin{tabular}{|c|c|c|c|c|c|c|c|c|}
  \hline
    $\mfrak{g}$ & Dynkin diagram & $k_\mfrak{p}$ & $\Pi^\natural$ & $\Pi^\natural \setminus \Sigma^\natural$ & $\mfrak{g}^\natural$ & $\tilde{\mfrak{g}}^\natural$ & $k_{\tilde{\mfrak{p}}^\natural}$ & $X_{\mcal{V}_{\kappa_\mfrak{p}}\!(\mfrak{g})}$ \\
  \hline
    $\mfrak{sl}_{n+1}$, $n\geq 2$ & $\dynkin[edge length=7mm, labels*={1,2,n-1,n}, label macro*/.code={\alpha_{#1}}, x/.style={thin}]{A}{oo..oo}$ & $-1$ & $\{\alpha_2,\dots,\alpha_{n-1}\}$ & $\emptyset$ & $\mfrak{gl}_1 \oplus \mfrak{sl}_{n-1}$ & $\mfrak{gl}_1$ & ${n-1 \over 2}$ & $\widebar{\mcal{S}^*_{\mfrak{p}_{\alpha_1}}}$ \\
  \hline
    $\mfrak{sl}_{2n}$, $n\geq 2$ & $\dynkin[edge length=7mm, labels*={1,2,2n-2,2n-1}, label macro*/.code={\alpha_{#1}}, x/.style={thin}]{A}{oo..oo}$ & $-n$ & $\{\alpha_2,\dots,\alpha_{2n-2}\}$ & $\{\alpha_n\}$ & $\mfrak{gl}_1 \oplus \mfrak{sl}_{2n-2}$ & $\mfrak{sl}_{2n-2}$ & $-n+1$ & $\widebar{\mcal{S}^*_{\mfrak{p}_{\alpha_n}}}$ \\
  \hline
    $\mfrak{so}_{2n+1}$, $n\geq 3$ & $\dynkin[edge length=7mm, labels*={1,2,n-1,n}, label macro*/.code={\alpha_{#1}}, x/.style={thin}]{B}{oo.oo}$ & $-2$ & $\{\alpha_1,\alpha_3,\dots,\alpha_n\}$ & $\{\alpha_1\}$ & $\mfrak{sl}_2 \oplus \mfrak{so}_{2n-3}$ & $\mfrak{sl}_2$ & $n-{7 \over 2}$ & $\widebar{{\mcal{O}}^*_{[3,1^{2n-2}]}}$ \\
  \hline
    $\mfrak{sp}_{2n}$, $n \geq 2$ & $\dynkin[edge length=7mm, labels*={1,2,n-1,n}, label macro*/.code={\alpha_{#1}}, x/.style={thin}]{C}{oo.oo}$ & $-{n \over 2}-1$ & $\{\alpha_2,\dots,\alpha_n\}$ & $\{\alpha_n\}$ & $\mfrak{sp}_{2n-2}$ & $\mfrak{sp}_{2n-2}$ & $-{n-1 \over 2}-1$ & $\widebar{\mcal{S}^*_{\mfrak{p}_{\alpha_n}}}$\\
  \hline
    $\mfrak{so}_{2n}$, $n\geq 4$ & $\dynkin[edge length=7mm, labels*={1,2}, label macro*/.code={\alpha_{#1}}, labels={,,n-2,n-1,n}, label macro/.code={\alpha_{#1}}, x/.style={thin}]{D}{oo.ooo}$ & $-2$ & $\{\alpha_1,\alpha_3,\dots,\alpha_n\}$ & $\{\alpha_1\}$ & $\mfrak{sl}_2 \oplus \mfrak{so}_{2n-4}$ & $\mfrak{sl}_2$ & $n-4$ & $\widebar{\mcal{O}^*_{\rm min}}$ \\
  \hline
    $\mfrak{so}_{2n}$, $n\geq 4$, $n \in 2\N+1$ & $\dynkin[edge length=7mm, labels*={1,2}, label macro*/.code={\alpha_{#1}}, labels={,,n-2,n-1,n}, label macro/.code={\alpha_{#1}}, x/.style={thin}]{D}{oo.ooo}$ & $-n+2$ & $\{\alpha_1,\alpha_3,\dots,\alpha_n\}$ & $\{\alpha_n\}$ & $\mfrak{sl}_2 \oplus \mfrak{so}_{2n-4}$ & $\mfrak{so}_{2n-4}$ & $-n+4$ & $\widebar{\mcal{S}^*_{\mfrak{p}_{\alpha_n}}}$ \\
  \hline
    $\mfrak{so}_{2n}$, $n\geq 4$, $n \in 2\N$ & $\dynkin[edge length=7mm, labels*={1,2}, label macro*/.code={\alpha_{#1}}, labels={,,n-2,n-1,n}, label macro/.code={\alpha_{#1}}, x/.style={thin}]{D}{oo.ooo}$ & $-n+2$ & $\{\alpha_1,\alpha_3,\dots,\alpha_n\}$ & $\{\alpha_n\}$ & $\mfrak{sl}_2 \oplus \mfrak{so}_{2n-4}$ & $\mfrak{so}_{2n-4}$ & $-n+4$ & $\widebar{{\mcal{O}}^*_{[2^{n-2},1^4]}}$ \\
  \hline
    $\mfrak{e}_6$ & $\dynkin[edge length=7mm, labels={1,2,3,4,5,6}, label macro/.code={\alpha_{#1}}, ordering=Dynkin, x/.style={thin}]{E}{oooooo}$ & $-4$ & $\{\alpha_1,\dots,\alpha_5\}$ & $\{\alpha_1\}$ & $\mfrak{sl}_6$ & $\mfrak{sl}_6$ & $-1$ & $\widebar{\mcal{S}^*_{\mfrak{p}_{\alpha_1}}}$ \\
  \hline
    $\mfrak{e}_7$ & $\dynkin[edge length=7mm, labels={1,2,3,4,5,6,7}, label macro/.code={\alpha_{#1}}, ordering=Dynkin, x/.style={thin}]{E}{ooooooo}$ & $-6$ & $\{\alpha_2,\dots,\alpha_7\}$ & $\{\alpha_6\}$ & $\mfrak{so}_{12}$ & $\mfrak{so}_{12}$ & $-2$ & $\widebar{{\mcal{O}}^*_{2A_1}}$ \\
  \hline
\end{tabular}
\caption{Associated varieties of simple affine vertex algebras}
\label{tab:associated varieties}
\end{sidewaystable}


\section{Free field realizations of affine vertex algebras}

In this section we provide a free field realization of affine vertex algebras $\mcal{V}^{\kappa_\mfrak{p}}\!(\mfrak{g})$ of level $\kappa_\mfrak{p}$ for $\mfrak{g}$ of the classical type described in Table \ref{tab:parabolic commutative nilradical level}. It is a straightforward application of Theorem \ref{thm:FF homomorphism 1-graded}.
\medskip

For $n \in \N$, we shall denote by $M_{n \times n}(\C)$ the algebra of $(n\times n)$-matrices.
The identity matrix of $M_{n \times n}(\C)$ is $I_n$ and $E_{i,j} \in M_{n\times n}(\C)$ for $1\leq i,j \leq n$ stands for the standard elementary $(n\times n)$-matrix having $1$ at the intersection of the $i$-th row and the $j$-th column and $0$ elsewhere.


\subsection{Lie algebra $\mfrak{sl}_{n+1}$ with the parabolic subalgebra $\mfrak{p}_{\alpha_1}$}

Let us consider the simple Lie algebra $\mfrak{g}=\mfrak{sl}_{n+1}$ with $n\geq 1$. A Cartan subalgebra $\mfrak{h}$ of $\mfrak{g}$ is given by diagonal matrices
\begin{align*}
  \mfrak{h} = \{\diag(a_1,a_2,\dots,a_{n+1});\, a_1,a_2,\dots,a_{n+1} \in \C,\, {\textstyle \sum_{i=1}^{n+1}} a_i=0\}.
\end{align*}
For $i=1,2,\dots,n+1$, we define $\veps_i \in \mfrak{h}^*$ through $\veps_i(\diag(a_1,a_2,\dots,a_{n+1}))=a_i$. The root system of $\mfrak{g}$ with respect to $\mfrak{h}$ is then given by $\Delta=\{\veps_i-\veps_j;\, 1\leq i \neq j \leq n+1\}$. A positive root system in $\Delta$ is $\Delta_+=\{\veps_i - \veps_j;\, 1 \leq i < j \leq n+1\}$ with the set of simple roots $\Pi=\{\alpha_1,\alpha_2,\dots,\alpha_n\}$, $\alpha_i= \veps_i - \veps_{i+1}$ for $i=1,2,\dots,n$. The fundamental weights are $\omega_i = \smash{\sum_{j=1}^i} \veps_j$ for $i=1,2,\dots,n$ and the maximal root of $\mfrak{g}$ is $\theta=\veps_1-\veps_{n+1} = \alpha_1 + \alpha_2 + \dots + \alpha_n$. Further, we define root vectors of $\mfrak{g}$ by
\begin{align*}
  e_{\veps_i-\veps_j} = E_{i,j} \qquad \text{and} \qquad f_{\veps_i-\veps_j} = E_{j,i}
\end{align*}
together with the corresponding coroots
\begin{align*}
  h_{\veps_i-\veps_j} = E_{i,i} - E_{j,j}
\end{align*}
for $1 \leq i < j \leq n+1$, where $E_{i,j} \in M_{n+1\times n+1}(\C)$ for $1\leq i,j \leq n+1$. The normalized $\mfrak{g}$-invariant symmetric bilinear form $\kappa_0$ on $\mfrak{g}$ is given by
\begin{align*}
  \kappa_0(a,b) = \tr(ab)
\end{align*}
for $a,b \in \mfrak{g}$. Let us note that $\kappa_0$ is the trace form of the basic representation of $\mfrak{g}$.
\medskip

Let us denote by $\mfrak{p}$ the maximal opposite standard parabolic subalgebra of $\mfrak{g}$ associated to the subset $\Sigma=\Pi\setminus \{\alpha_1\}$ of $\Pi$. We have the corresponding triangular decomposition
\begin{align*}
  \mfrak{g} = \widebar{\mfrak{u}} \oplus \mfrak{l} \oplus \mfrak{u}
\end{align*}
of $\mfrak{g}$ with $\mfrak{p}=\mfrak{l} \oplus \widebar{\mfrak{u}}$. It follows immediately that
\begin{align*}
  \Delta_+^\mfrak{l} = \{\veps_i-\veps_j;\, 2\leq i < j \leq n+1\} \qquad \text{and} \qquad \Delta_+^\mfrak{u} = \{\veps_1-\veps_i;\, 2\leq i \leq n+1\}.
\end{align*}
Furthermore, we introduce the element
\begin{align*}
  h_c = \begin{pmatrix}
    1 & 0 \\
    0 & -{1\over n}I_n
  \end{pmatrix}
\end{align*}
of $\mfrak{h}$ which spans the center $\mfrak{z}(\mfrak{l})$ of the Levi subalgebra $\mfrak{l}$. We will use the notation
\begin{align*}
  a_i^*(z) = a_{\veps_1-\veps_{i+1}}^*\!(z) \qquad \text{and} \qquad   a_i(z) = a_{\veps_1-\veps_{i+1}}\!(z)
\end{align*}
with $1 \leq i \leq n$ for the generating fields of the Weyl vertex algebra $\mcal{M}_\mfrak{u}$.
\medskip

\theorem{\label{thm:FF explicit sl at -1} For $\mfrak{g}=\mfrak{sl}_{n+1}$ with $n \geq 1$, there is a homomorphism
\begin{align*}
  \widetilde{w}_{\kappa_\mfrak{p},\mfrak{g}}^\mfrak{p} \colon \mcal{V}^{\kappa_\mfrak{p}}\!(\mfrak{g}) \rarr \mcal{M}_\mfrak{u} \otimes_\C \mcal{V}^{\kappa_\mfrak{p}+{1\over 2}\kappa_\mfrak{g}}\!(\mfrak{z}(\mfrak{l}))
\end{align*}
of $\N_0$-graded vertex algebras given by
\begin{align*}
    \widetilde{w}_{\kappa_\mfrak{p},\mfrak{g}}^\mfrak{p}(e_{\veps_1-\veps_{i+1}}\!(z)) &= -a_i(z), \\
    \widetilde{w}_{\kappa_\mfrak{p},\mfrak{g}}^\mfrak{p}(f_{\veps_1-\veps_{i+1}}\!(z)) &= \sum_{k=1}^n \normOrd{a^*_i(z)a^*_k(z)a_k(z)}+\partial_z a^*_i(z) -a^*_i(z)h_c(z)
\end{align*}
for $1 \leq i \leq n$,
\begin{align*}
    \widetilde{w}_{\kappa_\mfrak{p},\mfrak{g}}^\mfrak{p}(e_{\veps_{i+1}-\veps_{j+1}}\!(z)) = \normOrd{a^*_i(z)a_j(z)}, \qquad  \widetilde{w}_{\kappa_\mfrak{p},\mfrak{g}}^\mfrak{p}(f_{\veps_{i+1}-\veps_{j+1}}\!(z)) =  \normOrd{a^*_j(z)a_i(z)}
\end{align*}
for $1 \leq i < j \leq n$, and
\begin{align*}
    \widetilde{w}_{\kappa_\mfrak{p},\mfrak{g}}^\mfrak{p}(h_i(z)) = \begin{cases}
    -\normOrd{a^*_1(z)a_1(z)} - {\displaystyle \sum_{k=1}^n}\, \normOrd{a^*_k(z)a_k(z)} + h_c(z) & \text{if $i=1$}, \\[5mm]
    \normOrd{a^*_{i-1}(z)a_{i-1}(z)} - \normOrd{a^*_i(z)a_i(z)} & \text{if $1 <i \leq n$}
    \end{cases}
\end{align*}
for $1 \leq i \leq n$, where $\kappa_\mfrak{p} = -\kappa_0$ and $\kappa_\mfrak{g}=2(n+1)\kappa_0$.}

Let us note that the universal affine vertex algebras $\mcal{V}^{\kappa_\mfrak{p}}\!(\mfrak{g})$ of level $\kappa_\mfrak{p}$ is simple only if $n=1$ by \cite{Gorelik-Kac2007}.


\subsection{Lie algebra $\mfrak{sl}_{2n}$ with the parabolic subalgebra $\mfrak{p}_{\alpha_n}$}

Let us consider the simple Lie algebra $\mfrak{g}=\mfrak{sl}_{2n}$ with $n\geq 1$. A Cartan subalgebra $\mfrak{h}$ of $\mfrak{g}$ is given by diagonal matrices
\begin{align*}
  \mfrak{h} = \{\diag(a_1,a_2,\dots,a_{2n});\, a_1,a_2,\dots,a_{2n} \in \C,\, {\textstyle \sum_{i=1}^{2n}} a_i=0\}.
\end{align*}
For $i=1,2,\dots,2n$, we define $\veps_i \in \mfrak{h}^*$ through $\veps_i(\diag(a_1,a_2,\dots,a_{2n}))=a_i$. The root system of $\mfrak{g}$ with respect to $\mfrak{h}$ is then given by $\Delta=\{\veps_i-\veps_j;\, 1\leq i \neq j \leq 2n\}$. A positive root system in $\Delta$ is $\Delta_+=\{\veps_i - \veps_j;\, 1 \leq i < j \leq 2n\}$ with the set of simple roots $\Pi=\{\alpha_1,\alpha_2,\dots,\alpha_{2n-1}\}$, $\alpha_i= \veps_i - \veps_{i+1}$ for $i=1,2,\dots,2n-1$. The fundamental weights are $\omega_i = \smash{\sum_{j=1}^i} \veps_j$ for $i=1,2,\dots,2n-1$ and the maximal root of $\mfrak{g}$ is $\theta=\veps_1-\veps_{2n} = \alpha_1 + \alpha_2 + \dots + \alpha_{2n-1}$. Further, we define root vectors of $\mfrak{g}$ by
\begin{align*}
  e_{\veps_i-\veps_j} = E_{i,j} \quad \text{and} \quad f_{\veps_i-\veps_j} = E_{j,i}
\end{align*}
together with the corresponding coroots
\begin{align*}
  h_{\veps_i-\veps_j} = E_{i,i} - E_{j,j}
\end{align*}
for $1 \leq i < j \leq 2n$, where $E_{i,j} \in M_{2n\times 2n}(\C)$ for $1\leq i,j \leq 2n$. The normalized $\mfrak{g}$-invariant symmetric bilinear form $\kappa_0$ on $\mfrak{g}$ is given by
\begin{align*}
  \kappa_0(a,b) = \tr(ab)
\end{align*}
for $a,b \in \mfrak{g}$. Let us note that $\kappa_0$ is the trace form of the basic representation of $\mfrak{g}$.
\medskip

Let us denote by $\mfrak{p}$ the maximal opposite standard parabolic subalgebra of $\mfrak{g}$ associated to the subset $\Sigma=\Pi\setminus \{\alpha_n\}$ of $\Pi$. We have the corresponding triangular decomposition
\begin{align*}
  \mfrak{g} = \widebar{\mfrak{u}} \oplus \mfrak{l} \oplus \mfrak{u}
\end{align*}
of $\mfrak{g}$ with $\mfrak{p} = \mfrak{l} \oplus \widebar{\mfrak{u}}$. It follows immediately that
\begin{align*}
  \Delta_+^\mfrak{l} &= \{\veps_i-\veps_j;\, 1\leq i < j \leq n\} \cup \{\veps_{n+i}-\veps_{n+j};\, 1\leq i < j \leq n\}
  \intertext{and}
  \Delta_+^\mfrak{u} &= \{\veps_i-\veps_{n+j};\, 1\leq i, j \leq n\}.
\end{align*}
Furthermore, we introduce the element
\begin{align*}
  h_c = \begin{pmatrix}
    {1\over n}I_n & 0 \\
    0 & -{1\over n}I_n
  \end{pmatrix}
\end{align*}
of $\mfrak{h}$ which spans the center $\mfrak{z}(\mfrak{l})$ of the Levi subalgebra $\mfrak{l}$. We will use the notation
\begin{align*}
  a_{i,j}^*(z) = a_{\veps_i-\veps_{n+j}}^*\!(z) \qquad \text{and} \qquad  a_{i,j}(z) = a_{\veps_i-\veps_{n+j}}\!(z)
\end{align*}
with $1 \leq i,j \leq n$ for the generating fields of the Weyl vertex algebra $\mcal{M}_\mfrak{u}$.
\medskip

\theorem{\label{thm:FF explicit sl at -n} For $\mfrak{g}=\mfrak{sl}_{2n}$ with $n \geq 1$, there is a homomorphism
\begin{align*}
  \widetilde{w}_{\kappa_\mfrak{p},\mfrak{g}}^\mfrak{p} \colon \mcal{V}^{\kappa_\mfrak{p}}\!(\mfrak{g}) \rarr \mcal{M}_\mfrak{u} \otimes_\C \mcal{V}^{\kappa_\mfrak{p}+{1\over 2}\kappa_\mfrak{g}}\!(\mfrak{z}(\mfrak{l}))
\end{align*}
of $\N_0$-graded vertex algebras given by
\begin{align*}
  \widetilde{w}_{\kappa_\mfrak{p},\mfrak{g}}^\mfrak{p}(e_{\veps_i-\veps_{n+j}}(z)) &= -a_{i,j}(z), \\
  \widetilde{w}_{\kappa_\mfrak{p},\mfrak{g}}^\mfrak{p}(f_{\veps_i-\veps_{n+j}}(z)) &= \sum_{k,\ell =1}^n \normOrd{a^*_{k,j}(z) a^*_{i,\ell}(z)a_{k,\ell}(z)} +n\partial_za^*_{i,j}(z) - a^*_{i,j}(z)h_c(z)
\end{align*}
for $1 \leq i,j \leq n$,
\begin{align*}
  \widetilde{w}_{\kappa_\mfrak{p},\mfrak{g}}^\mfrak{p}(e_{\veps_i-\veps_j}(z)) &= -\sum_{k=1}^n \normOrd{a^*_{j,k}(z)a_{i,k}(z)}, &
  \widetilde{w}_{\kappa_\mfrak{p},\mfrak{g}}^\mfrak{p}(e_{\veps_{n+i}-\veps_{n+j}}(z)) &= \sum_{k=1}^n \normOrd{a^*_{k,i}(z)a_{k,j}(z)}, \\
  \widetilde{w}_{\kappa_\mfrak{p},\mfrak{g}}^\mfrak{p}(f_{\veps_i-\veps_j}(z)) &= -\sum_{k=1}^n \normOrd{a^*_{i,k}(z)a_{j,k}(z)}, &
  \widetilde{w}_{\kappa_\mfrak{p},\mfrak{g}}^\mfrak{p}(f_{\veps_{n+i}-\veps_{n+j}}(z)) &= \sum_{k=1}^n \normOrd{a^*_{k,j}(z)a_{k,i}(z)}
\end{align*}
for $1 \leq i < j \leq n$, and
\begin{align*}
  \widetilde{w}_{\kappa_\mfrak{p},\mfrak{g}}^\mfrak{p}(h_i(z)) =
  \begin{cases}
  -{\displaystyle \sum_{k=1}^n}\, \normOrd{a^*_{i,k}(z)a_{i,k}(z)} + {\displaystyle \sum_{k=1}^n}\, \normOrd{a^*_{i+1,k}(z)a_{i+1,k}(z)} & \text{if $1 \leq i < n$}, \\[5mm]
  -{\displaystyle \sum_{k=1}^n}\, \normOrd{a^*_{n,k}(z)a_{n,k}(z)} - {\displaystyle \sum_{k=1}^n}\, \normOrd{a^*_{k,1}(z)a_{k,1}(z)} + h_c(z) & \text{if $i=n$}, \\[5mm]
  {\displaystyle \sum_{k=1}^n}\, \normOrd{a^*_{k,i}(z)a_{k,i}(z)} - {\displaystyle \sum_{k=1}^n}\, \normOrd{a^*_{k,i+1}(z)a_{k,i+1}(z)} & \text{if $n < i \leq 2n-1$}
  \end{cases}
\end{align*}
for $1 \leq i \leq 2n-1$, where $\kappa_\mfrak{p} = -n\kappa_0$ and $\kappa_\mfrak{g}=4n\kappa_0$.}

Let us note that the universal affine vertex algebras $\mcal{V}^{\kappa_\mfrak{p}}\!(\mfrak{g})$ of level $\kappa_\mfrak{p}$ is simple only if $n=1$ by \cite{Gorelik-Kac2007}. In fact, both free field realizations of the affine vertex algebra $\mcal{V}^{-\kappa_0}(\mfrak{sl}_2)$ following from Theorem \ref{thm:FF explicit sl at -1} and Theorem \ref{thm:FF explicit sl at -n} are the same.


\subsection{Lie algebra $\mfrak{sp}_{2n}$ with the parabolic subalgebra $\mfrak{p}_{\alpha_n}$}

Let us consider the simple Lie algebra $\mfrak{g}=\mfrak{sp}_{2n}$ with $n\geq 2$ defined by
\begin{align*}
  \mfrak{sp}_{2n}= \left\{\!
  \begin{pmatrix}
    A & B \\
    C & -A^{\rm T}
  \end{pmatrix}\!;\,
  A,B,C\in M_{n,n}(\C),\, B^{\rm T}=B,\, C^{\rm T}=C\right\}\!.
\end{align*}
A Cartan subalgebra $\mfrak{h}$ of $\mfrak{g}$ is given by diagonal matrices
\begin{align*}
  \mfrak{h} = \{\diag(a_1,\dots,a_n,-a_1,\dots,-a_n);\, a_1,a_2,\dots,a_n \in \C\}.
\end{align*}
For $i=1,2,\dots,n$, we define $\veps_i \in \mfrak{h}^*$ by $\veps_i(\diag(a_1,\dots,a_n,-a_1,\dots,-a_n))=a_i$. The root system of $\mfrak{g}$ with respect to $\mfrak{h}$ is given by $\Delta=\{\pm\veps_i\pm\veps_j;\, 1\leq i < j \leq n\} \cup \{\pm 2\veps_i;\, 1 \leq i \leq n\}$. A positive root system in $\Delta$ is $\Delta_+=\{\veps_i \pm \veps_j;\, 1 \leq i < j \leq n\} \cup \{2\veps_i;\, 1 \leq i \leq n\}$ with the set of simple roots $\Pi=\{\alpha_1,\alpha_2,\dots,\alpha_n\}$, $\alpha_i= \veps_i - \veps_{i+1}$ for $i=1,2,\dots,n-1$ and $\alpha_n=2\veps_n$. The fundamental weights are $\omega_i = \smash{\sum_{j=1}^i} \veps_j$ for $i=1,2,\dots,n$ and the maximal root of $\mfrak{g}$ is $\theta = 2\veps_1 = 2\alpha_1 + \dots + 2\alpha_{n-1} + \alpha_n$. Further, we introduce root vectors of $\mfrak{g}$ by
\begin{align*}
  e_{\veps_i+\veps_j} = E_{i,n+j} + E_{j,n+i}, \qquad f_{\veps_i+\veps_j} = E_{n+i,j} + E_{n+j,i}, \\
  e_{\veps_i-\veps_j} = E_{i,j} - E_{n+j,n+i}, \qquad f_{\veps_i-\veps_j} = E_{j,i} - E_{n+i,n+j}
\end{align*}
with the corresponding coroots
\begin{align*}
  h_{\veps_i\pm \veps_j} = E_{i,i} \pm E_{j,j} - E_{n+i,n+i} \mp E_{n+j,n+j}
\end{align*}
for $1 \leq i < j \leq n$ and
\begin{align*}
  e_{2\veps_i} = E_{i,n+i}, \qquad f_{2\veps_i} = E_{n+i,i}
\end{align*}
with the corresponding coroots
\begin{align*}
  h_{2\veps_i} = E_{i,i} - E_{n+i,n+i}
\end{align*}
for $1 \leq i \leq n$, where $E_{i,j} \in M_{2n\times 2n}(\C)$ for $1\leq i,j \leq 2n$. The normalized $\mfrak{g}$-invariant symmetric bilinear form $\kappa_0$ on $\mfrak{g}$ is given by
\begin{align*}
  \kappa_0(a,b) = \tr(ab)
\end{align*}
for $a,b \in \mfrak{g}$. Let us note that $\kappa_0$ is the trace form of the basic representation of $\mfrak{g}$.
\medskip

Let us denote by $\mfrak{p}$ the maximal opposite standard parabolic subalgebra of $\mfrak{g}$ associated to the subset $\Sigma=\Pi\setminus \{\alpha_n\}$ of $\Pi$. We have the corresponding triangular decomposition
\begin{align*}
  \mfrak{g} = \widebar{\mfrak{u}} \oplus \mfrak{l} \oplus \mfrak{u}
\end{align*}
of $\mfrak{g}$ with $\mfrak{p} = \mfrak{l} \oplus \widebar{\mfrak{u}}$. It follows immediately that
\begin{align*}
  \Delta_+^\mfrak{l} = \{\veps_i-\veps_j;\, 1\leq i < j \leq n\} \qquad \text{and} \qquad \Delta_+^\mfrak{u} = \{\veps_i+\veps_j;\, 1\leq i \leq j \leq n\}.
\end{align*}
Furthermore, we introduce the element
\begin{align*}
  h_c = \begin{pmatrix}
    {1\over n}I_n & 0 \\
    0 & -{1\over n}I_n
  \end{pmatrix}
\end{align*}
of $\mfrak{h}$ which spans the center $\mfrak{z}(\mfrak{l})$ of the Levi subalgebra $\mfrak{l}$. We will use the notation
\begin{align*}
  a_{i,j}^*(z) = a_{j,i}^*(z) = a_{\veps_i+\veps_j}^*\!(z) \qquad \text{and} \qquad  a_{i,j}(z) = a_{j,i}(z) = a_{\veps_i+\veps_j}\!(z)
\end{align*}
with $1 \leq i \leq j \leq n$ for the generating fields of the Weyl vertex algebra $\mcal{M}_\mfrak{u}$.
\medskip

\theorem{\label{thm:FF explicit sp} For $\mfrak{g} = \mfrak{sp}_{2n}$ with $n \geq 2$, there is a homomorphism
\begin{align*}
  \widetilde{w}^\mfrak{p}_{\kappa_\mfrak{p},\mfrak{g}} \colon \mcal{V}^{\kappa_\mfrak{p}}\!(\mfrak{g}) \rarr \mcal{M}_\mfrak{u} \otimes_\C \mcal{V}^{\kappa_\mfrak{p} + {1 \over 2}\kappa_\mfrak{g}}\!(\mfrak{z}(\mfrak{l}))
\end{align*}
of $\N_0$-graded vertex algebras given by
\begin{align*}
  \widetilde{w}^\mfrak{p}_{\kappa_\mfrak{p},\mfrak{g}}(e_{\veps_i+\veps_j}\!(z)) &= -a_{i,j}(z), \\
  \widetilde{w}^\mfrak{p}_{\kappa_\mfrak{p},\mfrak{g}}(f_{\veps_i+\veps_j}\!(z)) &= \sum_{k,\ell =1}^n\! {\textstyle { 1+\delta_{k,\ell} \over 1+\delta_{i,j}}}\,\normOrd{a^*_{i,k}(z) a^*_{j,\ell}(z)a_{k,\ell}(z)} +(2-\delta_{i,j})({\textstyle (1+{n\over 2})}\partial_za^*_{i,j}(z) - a^*_{i,j}(z)h_c(z))
\end{align*}
for $1 \leq i \leq j \leq n$,
\begin{align*}
  \widetilde{w}^\mfrak{p}_{\kappa_\mfrak{p},\mfrak{g}}(e_{\veps_i-\veps_j}\!(z)) = -\sum_{k=1}^n (1+\delta_{i,k})a^*_{j,k}(z)a_{i,k}(z), \quad
  \widetilde{w}^\mfrak{p}_{\kappa_\mfrak{p},\mfrak{g}}(f_{\veps_i-\veps_j}\!(z)) = -\sum_{k=1}^n (1+\delta_{j,k})a^*_{i,k}(z)a_{j,k}(z)
\end{align*}
for $1 \leq i < j \leq n$, and
\begin{align*}
  \widetilde{w}^\mfrak{p}_{\kappa_\mfrak{p},\mfrak{g}}(h_i(z)) & = -\sum_{k=1}^n (1+\delta_{i,k}) \normOrd{a^*_{i,k}(z)a_{i,k}(z)} + \sum_{k=1}^n (1+\delta_{i+1,k}) \normOrd{a^*_{i+1,k}(z)a_{i+1,k}(z)}, \\
  \widetilde{w}^\mfrak{p}_{\kappa_\mfrak{p},\mfrak{g}}(h_n(z)) & = -\sum_{k=1}^n (1+\delta_{n,k}) \normOrd{a^*_{n,k}(z)a_{n,k}(z)} + h_c(z),
\end{align*}
for $1 \leq i \leq n-1$, where $\kappa_\mfrak{p} = -(1+{n \over 2})\kappa_0$ and $\kappa_\mfrak{g}=2(n+1)\kappa_0$.}
\vspace{-2mm}


\subsection{Lie algebra $\mfrak{so}_{2n+1}$ with the parabolic subalgebra $\mfrak{p}_{\alpha_1}$}

Let us consider the simple Lie algebra $\mfrak{g}=\mfrak{so}_{2n+1}$ with $n\geq 2$ defined by
\begin{align*}
  \mfrak{so}_{2n+1}= \left\{\!
  \begin{pmatrix}
    A & B & u\\
    C & -A^{\rm T} & v \\
    -v^{\rm T} & -u^{\rm T} & 0
  \end{pmatrix}\!;\,u,v \in \C^n,\,
  A,B,C\in M_{n,n}(\C),\, B^{\rm T}=-B,\, C^{\rm T}=-C\right\}\!.
\end{align*}
A Cartan subalgebra $\mfrak{h}$ of $\mfrak{g}$ is given by diagonal matrices
\begin{align*}
  \mfrak{h} = \{\diag(a_1,\dots,a_n,-a_1,\dots,-a_n,0);\, a_1,a_2,\dots,a_n \in \C\}.
\end{align*}
For $i=1,2,\dots,n$, we define $\veps_i \in \mfrak{h}^*$ through $\veps_i(\diag(a_1,\dots,a_n,-a_1,\dots,-a_n,0))=a_i$. The root system of $\mfrak{g}$ with respect to $\mfrak{h}$ is then given by $\Delta=\{\pm\veps_i\pm\veps_j;\, 1\leq i < j \leq n\} \cup \{\pm \veps_i;\, 1 \leq i \leq n\}$. A positive root system in $\Delta$ is $\Delta_+=\{\veps_i \pm \veps_j;\, 1 \leq i < j \leq n\} \cup \{\veps_i;\, 1 \leq i \leq n\}$ with the set of simple roots $\Pi=\{\alpha_1,\alpha_2,\dots,\alpha_n\}$, $\alpha_i= \veps_i - \veps_{i+1}$ for $i=1,2,\dots,n-1$ and $\alpha_n=\veps_n$. The fundamental weights are $\omega_i = \smash{\sum_{j=1}^i} \veps_j$ for $i=1,2,\dots,n-1$, $\omega_n={1\over 2}\smash{\sum_{j=1}^n} \veps_j$ and the maximal root of $\mfrak{g}$ is $\theta = \veps_1+\veps_2 = \alpha_1 + 2\alpha_2 + \dots + 2\alpha_n$. Further, we introduce root vectors of $\mfrak{g}$ by
\begin{align*}
  e_{\veps_i+\veps_j} &= E_{i,n+j} - E_{j,n+i}, \qquad f_{\veps_i+\veps_j} = E_{n+j,i} - E_{n+i,j}, \\
  e_{\veps_i-\veps_j} &= E_{i,j} - E_{n+j,n+i}, \qquad f_{\veps_i-\veps_j} = E_{j,i} - E_{n+i,n+j}
\end{align*}
with the corresponding coroots
\begin{align*}
  h_{\veps_i\pm \veps_j} = E_{i,i} \pm E_{j,j} - E_{n+i,n+i} \mp E_{n+j,n+j}
\end{align*}
for $1 \leq i < j \leq n$ and
\begin{align*}
  e_{\veps_i} = E_{i,2n+1} - E_{2n+1,i+n}, \qquad f_{\veps_i} = E_{2n+1,i} - E_{i+n,2n+1}
\end{align*}
with the corresponding coroots
\begin{align*}
  h_{\veps_i} = 2E_{i,i} - 2E_{n+i,n+i}
\end{align*}
for $1 \leq i \leq n$, where $E_{i,j} \in M_{2n+1\times 2n+1}(\C)$ for $1\leq i,j \leq 2n+1$. The normalized $\mfrak{g}$-invariant symmetric bilinear form $\kappa_0$ on $\mfrak{g}$ is given by
\begin{align*}
  \kappa_0(a,b) = {1 \over 2}\tr(ab)
\end{align*}
for $a,b \in \mfrak{g}$. Let us note that $\kappa_0$ is the trace form of the basic representation of $\mfrak{g}$.
\medskip

Let us denote by $\mfrak{p}$ the maximal opposite standard parabolic subalgebra of $\mfrak{g}$ associated to the subset $\Sigma=\Pi\setminus \{\alpha_1\}$ of $\Pi$. We have the corresponding triangular decomposition
\begin{align*}
  \mfrak{g} = \widebar{\mfrak{u}} \oplus \mfrak{l} \oplus \mfrak{u}
\end{align*}
of $\mfrak{g}$ with $\mfrak{p} = \mfrak{l} \oplus \widebar{\mfrak{u}}$. It follows immediately that
\begin{align*}
  \Delta_+^\mfrak{l} = \{\veps_i\pm\veps_j;\, 2\leq i < j \leq n\} \cup \{\veps_i;\, 2 \leq i \leq n\} \quad \text{and} \quad \Delta_+^\mfrak{u} = \{\veps_1 \pm\veps_i;\, 2\leq i \leq n\} \cup \{\veps_1\}.
\end{align*}
Furthermore, we introduce the element
\begin{align*}
  h_c = E_{1,1} - E_{n+1,n+1}
\end{align*}
of $\mfrak{h}$ which spans the center $\mfrak{z}(\mfrak{l})$ of the Levi subalgebra $\mfrak{l}$. We will use the notation
\begin{align*}
\begin{aligned}
  a_i^*(z) = \begin{cases}
    a_{\veps_1-\veps_{i+1}}^*\!(z) & \text{if $1 \leq i < n$}, \\
    a_{\veps_1}^*\!(z) & \text{if $i=n$}, \\
    a_{\veps_1+\veps_{i-n+1}}^*\!(z) & \text{if $n < i \leq 2n-1$},
  \end{cases}
\end{aligned}  \qquad
\begin{aligned}
  a_i(z) = \begin{cases}
    a_{\veps_1-\veps_{i+1}}\!(z) & \text{if $1 \leq i < n$}, \\
    a_{\veps_1}\!(z) & \text{if $i=n$}, \\
    a_{\veps_1+\veps_{i-n+1}}\!(z) & \text{if $n < i \leq 2n-1$}
  \end{cases}
\end{aligned}
\end{align*}
with $1 \leq i \leq 2n-1$ for the generating fields of the Weyl vertex algebra $\mcal{M}_\mfrak{u}$.
\medskip

\theorem{\label{thm:FF explicit so at -2 odd} For $\mfrak{g} = \mfrak{so}_{2n+1}$ with $n \geq 2$, there is a homomorphism
\begin{align*}
  \widetilde{w}^\mfrak{p}_{\kappa_\mfrak{p},\mfrak{g}} \colon \mcal{V}^{\kappa_\mfrak{p}}\!(\mfrak{g}) \rarr \mcal{M}_\mfrak{u} \otimes_\C \mcal{V}^{\kappa_\mfrak{p} + {1 \over 2}\kappa_\mfrak{g}}\!(\mfrak{z}(\mfrak{l}))
\end{align*}
of $\N_0$-graded vertex algebras given by
\begin{align*}
  \widetilde{w}^\mfrak{p}_{\kappa_\mfrak{p},\mfrak{g}}(e_{\veps_1-\veps_{1+i}}\!(z)) &= -a_i(z), \\
  \widetilde{w}^\mfrak{p}_{\kappa_\mfrak{p},\mfrak{g}}(f_{\veps_1-\veps_{1+i}}\!(z)) &= \sum_{k=1}^{2n-1} \normOrd{a_i^*(z)a_k^*(z)a_k(z)} - \sum_{k=1}^{n-1} \normOrd{a_k^*(z)a_{n+k}^*(z)a_{n+i}(z)} - {\textstyle {1 \over 2}} \normOrd{a_n^*(z)^2a_{n+i}(z)} \\
  &\quad + 2\partial_z a_i^*(z) - a_i^*(z)h_c(z), \\
  \widetilde{w}^\mfrak{p}_{\kappa_\mfrak{p},\mfrak{g}}(e_{\veps_1+\veps_{1+i}}\!(z)) &= -a_{n+i}(z), \\
  \widetilde{w}^\mfrak{p}_{\kappa_\mfrak{p},\mfrak{g}}(f_{\veps_1+\veps_{1+i}}\!(z)) &= \sum_{k=1}^{2n-1} \normOrd{a_{n+i}^*(z)a_k^*(z)a_k(z)} - \sum_{k=1}^{n-1} \normOrd{a_k^*(z)a_{n+k}^*(z)a_i(z)} - {\textstyle {1 \over 2}} \normOrd{a_n^*(z)^2a_i(z)} \\
  &\quad +2\partial_z a_{n+i}^*(z) - a_{n+i}^*(z)h_c(z), \\
  \widetilde{w}^\mfrak{p}_{\kappa_\mfrak{p},\mfrak{g}}(e_{\veps_1}\!(z)) &= -a_n(z), \\
  \widetilde{w}^\mfrak{p}_{\kappa_\mfrak{p},\mfrak{g}}(f_{\veps_1}\!(z)) &= \sum_{k=1}^{2n-1} \normOrd{a_n^*(z)a_k^*(z)a_k(z)} - \sum_{k=1}^{n-1} \normOrd{a_k^*(z)a_{n+k}^*(z)a_n(z)} - {\textstyle {1 \over 2}} \normOrd{a_n^*(z)^2a_n(z)} \\
  &\quad +2\partial_z a_n^*(z) - a_n^*(z)h_c(z)
\end{align*}
for $1 \leq i \leq n-1$,
\begin{align*}
  \widetilde{w}^\mfrak{p}_{\kappa_\mfrak{p},\mfrak{g}}(e_{\veps_{1+i}-\veps_{1+j}}\!(z)) &= \normOrd{a_i^*(z)a_j(z)} - \normOrd{a_{n+j}^*(z)a_{n+i}(z)}, \\
  \widetilde{w}^\mfrak{p}_{\kappa_\mfrak{p},\mfrak{g}}(f_{\veps_{1+i}-\veps_{1+j}}\!(z)) &= \normOrd{a_j^*(z)a_i(z)} - \normOrd{a_{n+i}^*(z)a_{n+j}(z)}, \\
  \widetilde{w}^\mfrak{p}_{\kappa_\mfrak{p},\mfrak{g}}(e_{\veps_{1+i}+\veps_{1+j}}\!(z)) &= \normOrd{a_i^*(z)a_{n+j}(z)} - \normOrd{a_j^*(z)a_{n+i}(z)}, \\
  \widetilde{w}^\mfrak{p}_{\kappa_\mfrak{p},\mfrak{g}}(f_{\veps_{1+i}+\veps_{1+j}}\!(z)) &= \normOrd{a_{n+j}^*(z)a_i(z)} - \normOrd{a_{n+i}^*(z)a_j(z)}
\end{align*}
for $1 \leq i < j \leq n-1$,
\begin{align*}
  \widetilde{w}^\mfrak{p}_{\kappa_\mfrak{p},\mfrak{g}}(e_{\veps_{1+i}}\!(z)) &= \normOrd{a_i^*(z)a_n(z)} -\normOrd{a_n^*(z)a_{n+i}(z)}, \\
  \widetilde{w}^\mfrak{p}_{\kappa_\mfrak{p},\mfrak{g}}(f_{\veps_{1+i}}\!(z)) &= -\normOrd{a_{n+i}^*(z)a_n(z)} +\normOrd{a_n^*(z)a_i(z)}
\end{align*}
for $1 \leq i \leq n-1$, and
\begin{align*}
  \widetilde{w}^\mfrak{p}_{\kappa_\mfrak{p},\mfrak{g}}(h_i(z)) =
  \begin{cases}
    -\normOrd{a_1^*(z)a_1(z)} + \normOrd{a_{n+1}^*(z)a_{n+1}(z)} -{\displaystyle \sum_{k=1}^{2n-1}} \normOrd{a_k^*(z)a_k(z)} + h_c(z) & \text{if $i=1$}, \\[5mm]
    \begin{gathered}[c]
    \normOrd{a_{i-1}^*(z)a_{i-1}(z)} - \normOrd{a_i^*(z)a_i(z)} - \normOrd{a_{n+i-1}^*(z)a_{n+i-1}(z)}
    \\ + \normOrd{a_{n+i}^*(z)a_{n+i}(z)}
    \end{gathered} & \text{if $1 < i < n$}, \\[4mm]
    2\normOrd{a_{n-1}^*(z)a_{n-1}(z)} - 2\normOrd{a_{2n-1}^*(z)a_{2n-1}(z)} & \text{if $i=n$}
  \end{cases}
\end{align*}
for $1 \leq i \leq n$, where $\kappa_\mfrak{p} = -2\kappa_0$ and $\kappa_\mfrak{g}= 2(2n-1)\kappa_0$.}
\vspace{-2mm}


\subsection{Lie algebra $\mfrak{so}_{2n}$ with the parabolic subalgebras $\mfrak{p}_{\alpha_1}$ and $\mfrak{p}_{\alpha_n}$}

Let us consider the simple Lie algebra $\mfrak{g}=\mfrak{so}_{2n}$ with $n\geq 3$ defined by
\begin{align*}
  \mfrak{so}_{2n}= \left\{\!
  \begin{pmatrix}
    A & B \\
    C & -A^{\rm T}
  \end{pmatrix}\!;\,
  A,B,C\in M_{n,n}(\C),\, B^{\rm T}=-B,\, C^{\rm T}=-C\right\}\!.
\end{align*}
A Cartan subalgebra $\mfrak{h}$ of $\mfrak{g}$ is given by diagonal matrices
\begin{align*}
  \mfrak{h} = \{\diag(a_1,\dots,a_n,-a_1,\dots,-a_n);\, a_1,a_2,\dots,a_n \in \C\}.
\end{align*}
For $i=1,2,\dots,n$, we define $\veps_i \in \mfrak{h}^*$ by $\veps_i(\diag(a_1,\dots,a_n,-a_1,\dots,-a_n))=a_i$. The root system of $\mfrak{g}$ with respect to $\mfrak{h}$ is then given by $\Delta=\{\pm\veps_i\pm\veps_j;\, 1\leq i < j \leq n\}$. A positive root system in $\Delta$ is $\Delta_+=\{\veps_i \pm \veps_j;\, 1 \leq i < j \leq n\}$ with the set of simple roots $\Pi=\{\alpha_1,\alpha_2,\dots,\alpha_n\}$, $\alpha_i= \veps_i - \veps_{i+1}$ for $i=1,2,\dots,n-1$ and $\alpha_n=\veps_{n-1}+\veps_n$. Besides, the fundamental weights are $\omega_i = \smash{\sum_{j=1}^i} \veps_j$ for $i=1,2,\dots,n-2$, $\omega_{n-1}={1\over 2}\big(\smash{\sum_{j=1}^{n-1}} \veps_j - \veps_n\big)$, $\omega_n={1\over 2}\big(\smash{\sum_{j=1}^{n-1}} \veps_j + \veps_n\big)$ and the maximal root of $\mfrak{g}$ is $\theta = \veps_1 + \veps_2 = \alpha_1 + 2\alpha_2 + \dots + 2\alpha_{n-2} + \alpha_{n-1} + \alpha_n$. Further, we introduce root vectors of $\mfrak{g}$ by
\begin{align*}
  e_{\veps_i+\veps_j} &= E_{i,j+n} - E_{j,i+n}, \qquad f_{\veps_i+\veps_j} = E_{j+n,i} - E_{i+n,j}, \\
  e_{\veps_i-\veps_j} &= E_{i,j} - E_{j+n,i+n}, \qquad f_{\veps_i-\veps_j} = E_{j,i} - E_{i+n,j+n}
\end{align*}
for $1 \leq i < j \leq n$ together with the corresponding coroots given through
\begin{align*}
  h_{\veps_i\pm \veps_j} = E_{i,i} \pm E_{j,j} - E_{i+n,i+n} \mp E_{j+n,j+n}
\end{align*}
for $1 \leq i < j \leq n$, where $E_{i,j} \in M_{2n\times 2n}(\C)$ for $1\leq i,j \leq 2n$. Moreover, the normalized $\mfrak{g}$-invariant symmetric bilinear form $\kappa_0$ on $\mfrak{g}$ is given by
\begin{align*}
  \kappa_0(a,b) = {1 \over 2}\tr(ab)
\end{align*}
for $a,b \in \mfrak{g}$. Let us note that $\kappa_0$ is the trace form of the basic representation of $\mfrak{g}$.
\medskip

Let us denote by $\mfrak{p}$ the maximal opposite standard parabolic subalgebra of $\mfrak{g}$ associated to the subset $\Sigma=\Pi\setminus \{\alpha_1\}$ of $\Pi$. We have the corresponding triangular decomposition
\begin{align*}
  \mfrak{g} = \widebar{\mfrak{u}} \oplus \mfrak{l} \oplus \mfrak{u}
\end{align*}
of $\mfrak{g}$ with $\mfrak{p} = \mfrak{l} \oplus \widebar{\mfrak{u}}$. It follows immediately that
\begin{align*}
  \Delta_+^\mfrak{l} = \{\veps_i\pm\veps_j;\, 2\leq i < j \leq n\} \qquad \text{and} \qquad \Delta_+^\mfrak{u} = \{\veps_1 \pm\veps_i;\, 2\leq i \leq n\}.
\end{align*}
Furthermore, we introduce the element
\begin{align*}
  h_c = E_{1,1} - E_{n+1,n+1}
\end{align*}
of $\mfrak{h}$ which spans the center $\mfrak{z}(\mfrak{l})$ of the Levi subalgebra $\mfrak{l}$. We will use the notation
\begin{align*}
\begin{aligned}
  a_i^*(z) = \begin{cases}
    a_{\veps_1-\veps_{i+1}}^*\!(z) & \text{if $1 \leq i < n$}, \\
    a_{\veps_1+\veps_{i-n+2}}^*\!(z) & \text{if $n \leq i \leq 2n-2$},
  \end{cases}
\end{aligned}  \qquad
\begin{aligned}
  a_i(z) = \begin{cases}
    a_{\veps_1-\veps_{i+1}}\!(z) & \text{if $1 \leq i < n$}, \\
    a_{\veps_1+\veps_{i-n+2}}\!(z) & \text{if $n \leq i \leq 2n-2$}
  \end{cases}
\end{aligned}
\end{align*}
with $1 \leq i \leq 2n-2$ for the generating fields of the Weyl vertex algebra $\mcal{M}_\mfrak{u}$.
\medskip

\theorem{\label{thm:FF explicit so at -2 even} For $\mfrak{g} = \mfrak{so}_{2n}$ with $n \geq 3$, there is a homomorphism
\begin{align*}
  \widetilde{w}^\mfrak{p}_{\kappa_\mfrak{p},\mfrak{g}} \colon \mcal{V}^{\kappa_\mfrak{p}}\!(\mfrak{g}) \rarr \mcal{M}_\mfrak{u} \otimes_\C \mcal{V}^{\kappa_\mfrak{p} + {1 \over 2}\kappa_\mfrak{g}}\!(\mfrak{z}(\mfrak{l}))
\end{align*}
of $\N_0$-graded vertex algebras given by
\begin{align*}
  \widetilde{w}^\mfrak{p}_{\kappa_\mfrak{p},\mfrak{g}}(e_{\veps_1-\veps_{1+i}}\!(z)) &= -a_i(z), \\
  \widetilde{w}^\mfrak{p}_{\kappa_\mfrak{p},\mfrak{g}}(f_{\veps_1-\veps_{1+i}}\!(z)) &= \sum_{k=1}^{2n-2} \normOrd{a_i^*(z)a_k^*(z)a_k(z)} - \sum_{k=1}^{n-1} \normOrd{a_k^*(z)a_{n-1+k}^*(z)a_{n-1+i}(z)} \\
  &\quad + 2\partial_z a_i^*(z) - a_i^*(z)h_c(z), \\
  \widetilde{w}^\mfrak{p}_{\kappa_\mfrak{p},\mfrak{g}}(e_{\veps_1+\veps_{1+i}}\!(z)) &= -a_{n-1+i}(z), \\
  \widetilde{w}^\mfrak{p}_{\kappa_\mfrak{p},\mfrak{g}}(f_{\veps_1+\veps_{1+i}}\!(z)) &= \sum_{k=1}^{2n-2} \normOrd{a_{n-1+i}^*(z)a_k^*(z)a_k(z)} - \sum_{k=1}^{n-1} \normOrd{a_k^*(z)a_{n-1+k}^*(z)a_i(z)} \\
  &\quad +2\partial_z a_{n-1+i}^*(z) - a_{n-1+i}^*(z)h_c(z)
\end{align*}
for $1 \leq i \leq n-1$,
\begin{align*}
  \widetilde{w}^\mfrak{p}_{\kappa_\mfrak{p},\mfrak{g}}(e_{\veps_{1+i}-\veps_{1+j}}\!(z)) &= \normOrd{a_i^*(z)a_j(z)} - \normOrd{a_{n-1+j}^*(z)a_{n-1+i}(z)}, \\
  \widetilde{w}^\mfrak{p}_{\kappa_\mfrak{p},\mfrak{g}}(f_{\veps_{1+i}-\veps_{1+j}}\!(z)) &= \normOrd{a_j^*(z)a_i(z)} - \normOrd{a_{n-1+i}^*(z)a_{n-1+j}(z)}, \\
  \widetilde{w}^\mfrak{p}_{\kappa_\mfrak{p},\mfrak{g}}(e_{\veps_{1+i}+\veps_{1+j}}\!(z)) &= \normOrd{a_i^*(z)a_{n-1+j}(z)} - \normOrd{a_j^*(z)a_{n-1+i}(z)}, \\
  \widetilde{w}^\mfrak{p}_{\kappa_\mfrak{p},\mfrak{g}}(f_{\veps_{1+i}+\veps_{1+j}}\!(z)) &= \normOrd{a_{n-1+j}^*(z)a_i(z)} - \normOrd{a_{n-1+i}^*(z)a_j(z)}
\end{align*}
for $1 \leq i < j \leq n-1$, and
\begin{align*}
  \widetilde{w}^\mfrak{p}_{\kappa_\mfrak{p},\mfrak{g}}(h_i(z)) =
  \begin{cases}
    -\normOrd{a_1^*(z)a_1(z)} + \normOrd{a_n^*(z)a_n(z)} -{\displaystyle \sum_{k=1}^{2n-2}} \normOrd{a_k^*(z)a_k(z)} + h_c(z) & \text{if $i=1$}, \\[5mm]
    \begin{gathered}[c]
    \normOrd{a_{i-1}^*(z)a_{i-1}(z)} - \normOrd{a_i^*(z)a_i(z)} - \normOrd{a_{n-1+i-1}^*(z)a_{n-1+i-1}(z)}
    \\ + \normOrd{a_{n-1+i}^*(z)a_{n-1+i}(z)}
    \end{gathered} & \text{if $1 < i < n$}, \\[4mm]
    \begin{gathered}[c]
    \normOrd{a_{n-2}^*(z)a_{n-2}(z)} + \normOrd{a_{n-1}^*(z)a_{n-1}(z)} - \normOrd{a_{2n-3}^*(z)a_{2n-3}(z)}
    \\ - \normOrd{a_{2n-2}^*(z)a_{2n-2}(z)}
    \end{gathered} & \text{if $i=n$}
  \end{cases}
\end{align*}
for $1 \leq i \leq n$, where $\kappa_\mfrak{p} = -2\kappa_0$ and $\kappa_\mfrak{g}=4(n-1)\kappa_0$.}

Let us denote by $\mfrak{p}$ the maximal opposite standard parabolic subalgebra of $\mfrak{g}$ associated to the subset $\Sigma=\Pi\setminus \{\alpha_n\}$ of $\Pi$. We have the corresponding triangular decomposition
\begin{align*}
  \mfrak{g} = \widebar{\mfrak{u}} \oplus \mfrak{l} \oplus \mfrak{u}
\end{align*}
of $\mfrak{g}$ with $\mfrak{p} = \mfrak{l} \oplus \widebar{\mfrak{u}}$. It follows immediately that
\begin{align*}
  \Delta_+^\mfrak{l} = \{\veps_i-\veps_j;\, 1\leq i < j \leq n\} \qquad \text{and} \qquad \Delta_+^\mfrak{u} = \{\veps_i+\veps_j;\, 1\leq i < j \leq n\}.
\end{align*}
Furthermore, we introduce the element
\begin{align*}
  h_c = \begin{pmatrix}
    {1\over n}I_n & 0 \\
    0 & -{1\over n}I_n
  \end{pmatrix}
\end{align*}
of $\mfrak{h}$ which spans the center $\mfrak{z}(\mfrak{l})$ of the Levi subalgebra $\mfrak{l}$. We will use the notation
\begin{align*}
  a_{i,j}^*(z) = -a_{j,i}^*(z) = a_{\veps_i+\veps_j}^*\!(z) \qquad \text{and} \qquad  a_{i,j}(z) = -a_{j,i}(z) = a_{\veps_i+\veps_j}\!(z)
\end{align*}
with $1 \leq i < j \leq n$ for the generating fields of the Weyl vertex algebra $\mcal{M}_\mfrak{u}$.
\medskip

\theorem{\label{thm:FF explicit so} For $\mfrak{g} = \mfrak{so}_{2n}$ with $n \geq 3$, there is a homomorphism
\begin{align*}
  \widetilde{w}^\mfrak{p}_{\kappa_\mfrak{p},\mfrak{g}} \colon \mcal{V}^{\kappa_\mfrak{p}}\!(\mfrak{g}) \rarr \mcal{M}_\mfrak{u} \otimes_\C \mcal{V}^{\kappa_\mfrak{p} + {1 \over 2}\kappa_\mfrak{g}}\!(\mfrak{z}(\mfrak{l}))
\end{align*}
of $\N_0$-graded vertex algebras given by
\begin{align*}
  \widetilde{w}^\mfrak{p}_{\kappa_\mfrak{p},\mfrak{g}}(e_{\veps_i+\veps_j}\!(z)) &= -a_{i,j}(z), \\
  \widetilde{w}^\mfrak{p}_{\kappa_\mfrak{p},\mfrak{g}}(f_{\veps_i+\veps_j}\!(z)) &= \sum_{k,\ell =1}^n \normOrd{a^*_{i,k}(z) a^*_{\ell,j}(z)a_{\ell,k}(z)} +(n-2)\partial_za^*_{i,j}(z) - 2a^*_{i,j}(z)h_c(z)
\end{align*}
for $1 \leq i < j \leq n$,
\begin{align*}
  \widetilde{w}^\mfrak{p}_{\kappa_\mfrak{p},\mfrak{g}}(e_{\veps_i-\veps_j}\!(z)) = -\sum_{k=1}^n a^*_{j,k}(z)a_{i,k}(z), \qquad
  \widetilde{w}^\mfrak{p}_{\kappa_\mfrak{p},\mfrak{g}}(f_{\veps_i-\veps_j}\!(z)) = -\sum_{k=1}^n a^*_{i,k}(z)a_{j,k}(z)
\end{align*}
for $1 \leq i < j \leq n$, and
\begin{align*}
  \widetilde{w}^\mfrak{p}_{\kappa_\mfrak{p},\mfrak{g}}(h_i(z)) =
  \begin{cases}
    -{\displaystyle \sum_{k=1}^n}\, \normOrd{a^*_{i,k}(z)a_{i,k}(z)} + {\displaystyle \sum_{k=1}^n}\, \normOrd{a^*_{i+1,k}(z)a_{i+1,k}(z)} & \text{if $1 \leq i < n$}, \\[5mm]
    -{\displaystyle \sum_{k=1}^n}\, \normOrd{a^*_{n-1,k}(z)a_{n-1,k}(z)} - {\displaystyle \sum_{k=1}^n}\, \normOrd{a^*_{n,k}(z)a_{n,k}(z)} + 2h_c(z) & \text{if $i=n$}
    \end{cases}
\end{align*}
for $1 \leq i \leq n$, where $\kappa_\mfrak{p} = -(n-2)\kappa_0$ and $\kappa_\mfrak{g}=4(n-1)\kappa_0$.}
\vspace{-2mm}


\section*{Acknowledgments}

T.\,A.\ is partially supported by JSPS KAKENHI Grant Numbers 21H04993 and 19KK0065. V.\,F.\ is partially supported by the NSFC of China, Grant number 12350710178.



\newcommand{\etalchar}[1]{$^{#1}$}
\providecommand{\bysame}{\leavevmode\hbox to3em{\hrulefill}\thinspace}
\providecommand{\MR}{\relax\ifhmode\unskip\space\fi MR }
\providecommand{\MRhref}[2]{%
  \href{http://www.ams.org/mathscinet-getitem?mr=#1}{#2}
}
\providecommand{\href}[2]{#2}

\end{document}